\documentclass[11pt]{article}
\usepackage[utf8]{inputenc}
\usepackage[a4paper,margin=2cm]{geometry}
\usepackage[english]{babel}
\usepackage{latexsym,amsmath,enumerate,graphics,enumerate,amsthm,tikz,hyperref,float}
\usepackage[mathscr]{euscript}
\usepackage[affil-it]{authblk}
\usepackage{enumerate,amsthm,dsfont,pstricks}
\usepackage{latexsym,amsmath,amssymb,amscd,wrapfig,graphicx}
\usepackage{hyperref,cleveref}
\usepackage{circuitikz}
\usepackage{bbm}

\newtheorem{theorem}{Theorem}[section]
\newtheorem{lemma}[theorem]{Lemma}
\newtheorem{proposition}[theorem]{Proposition}
\newtheorem{corollary}[theorem]{Corollary}

\theoremstyle{definition}
\newtheorem{definition}[theorem]{Definition}
\newtheorem{example}[theorem]{Example}

\newtheorem*{acknowledgement}{Acknowledgement}
\theoremstyle{remark}
\newtheorem{remark}[theorem]{Remark}

\newtheorem*{theorem*}{{\bf Theorem}}

\newtheorem*{assumption*}{{\bf Assumption}}

\let\phi=\varphi

\def\R{\mathbb{R}}

\def\eps{\varepsilon}

\newcommand{\comment}[1]{}
\newcommand{\norm}[1]{\left\Vert #1 \right\Vert}

\newcommand{\one}{\mathbbm{1}}

\numberwithin{equation}{section}

\let\epsilon=\varepsilon

\makeatletter
\def\@maketitle{%
  \newpage
  \null
  \vskip 2em%
  \begin{center}%
  \let \footnote \thanks
    {\Large\bfseries \@title \par}%
    \vskip 1.5em%
    {\normalsize
      \lineskip .5em%
      \begin{tabular}[t]{c}%
        \@author
      \end{tabular}\par}%
    \vskip 1em%
    {\normalsize \@date}%
  \end{center}%
  \par
  \vskip 1.5em}
\makeatother

\begin{document}

\title{\sc \huge Order theoretical structures in atomic JBW-algebras: disjointness, bands, and centres}

\author{Onno van Gaans%
\thanks{Email: \texttt{vangaans@math.leidenuniv.nl}}}
\affil{Mathematical Institute, Leiden University, 2300 RA Leiden,
The Netherlands}

\author{Anke Kalauch%
\thanks{Email: \texttt{anke.kalauch@tu-dresden.de}}}
\affil{Department of Mathematics, Technische Universit\"{a}t Dresden, 01062 Dresden, Germany}

\author{Mark Roelands%
\thanks{Email: \texttt{m.roelands@math.leidenuniv.nl}}}
\affil{Mathematical Institute, Leiden University, 2300 RA Leiden,
The Netherlands}

\maketitle
\date{}

\begin{abstract} 
Every atomic JBW-algebra is known to be a direct sum of JBW-algebra factors of type I. Extending Kadison's anti-lattice theorem, we show that each of these factors is a disjointness free anti-lattice. We characterise disjointness, bands, and disjointness preserving bijections with disjointness preserving inverses in direct sums of disjointness free anti-lattices and, therefore, in atomic JBW-algebras. We show that in unital JB-algebras the algebraic centre and the order theoretical centre are isomorphic. Moreover, the order theoretical centre is a Riesz space of multiplication operators. A survey of JBW-algebra factors of type I is included.
\end{abstract}



\section{Introduction}

Jordan algebras equipped with their cones of squares are interesting instances of partially ordered vector spaces that are not lattices, in general. A prominent example is the Jordan algebra $\mathrm{B}(H)_\mathrm{sa}$ consisting of all self-adjoint operators on some complex Hilbert space $H$ with the Jordan product given by 
\begin{equation}\label{eq.jordanproduct}
A\circ B= \frac{1}{2}(AB+BA).
\end{equation} 
Kadison \cite{Kad1951b} has shown that this space is actually an anti-lattice, which means that the supremum of two elements exists only if they are comparable. With the notion of disjointness in partially ordered vector spaces \cite{KalvanGaa2019}, a partially ordered vector space is an anti-lattice if and only if there are no non-trivial disjoint positive elements \cite[Theorem 14]{KalLemGaa2014}. In the space $\mathrm{B}(H)_\mathrm{sa}$, it turns out that there are even no disjoint elements at all. We call such a partially ordered vector space \emph{disjointness free}. The space $\mathrm{B}(H)_\mathrm{sa}$ is one of the possible factors in the algebraic direct sum that represents atomic JBW-algebras. In this paper, we study all atomic JBW-algebras that are factors and show that all of them are disjointness free anti-lattices. By \cite[Theorem~3.39 and Proposition~3.45]{Alfsen}, every atomic JBW-algebra that is a factor is isomorphic as JBW-algebra to a member of one of the following classes of JBW-algebras,
\begin{itemize}
\item[(i)] the self-adjoint bounded operators $\mathrm{B}(H)_\mathrm{sa}$ on a real or complex Hilbert space $H$ of dimension $d\ge 3$, or $\mathrm{B}(\mathcal{H}_q)$ where $\mathcal{H}_q$ is a quaternionic Hilbert space of dimension $d\ge 3$, endowed with the product \eqref{eq.jordanproduct},
\item[(ii)] the spin factors $H\oplus\mathbb{R}$, where $H$ is a real Hilbert space of dimension at least 2, with the multiplication defined in \eqref{eq.spinproduct},
\item[(iii)] the $3\times 3$ self-adjoint matrices $\mathrm{M}_3(\mathbb{O})_\mathrm{sa}$ with entries from the octonions $\mathbb{O}$, endowed with the product \eqref{eq.jordanproduct}.
\end{itemize}
For general atomic JBW-algebras, there is a representation theorem as follows. See \cite[Proposition~3.45]{Alfsen}
\begin{theorem}\label{the.atomicdecomp}
Every atomic JBW-algebra equals the algebraic direct sum of atomic JBW-algebras that are factors, that is, of factors that are isomorphic as JBW-algebras to those listed in (i)--(iii).
\end{theorem}
The factors listed in (i)--(iii) are exactly the factors among all JBW-algebras that are of so called `type I', up to JBW-algebra isomorphism.

\Cref{the.atomicdecomp} leads to the question what can be said about disjointness and related notions in such direct sums. The algebraic direct sum in \Cref{the.atomicdecomp} is in fact an order direct sum of order unit spaces. We will characterise disjointness and bands in order direct sums of order unit spaces that are disjointness free anti-lattices. We will apply this characterisation to describe which disjointness preserving bijections have a disjointness preserving inverse, proceeding corresponding research in Banach lattices and finite-dimensional pre-Riesz spaces  \cite{HuiPag1993,KalLemGaa2019}. These results apply to atomic JBW-algebras.

In the theory of Jordan algebras, there is a notion of an algebraic centre, whereas, in the theory of operators on partially ordered vector spaces, there is a notion of an order theoretical centre. The algebraic centre of a Jordan algebra consists of all elements where the corresponding left multiplication operator commutes with all other left multiplication operators. The order theoretical centre of a partially ordered vector space consists of all operators that are in an order interval whose end points are multiples of the identity. We study the natural question how these two notions of centre are related. For a unital JB-algebra, we show that the algebraic centre and the order theoretical centre are isomorphic as JB-algebras. 

The structure of the paper is as follows. There are two preliminary sections with the basic relevant notions from the theory of partially ordered vector spaces and JB-algebras. 
We need quite a few details on the factors of atomic JBW-algebras, as listed above in (i)--(iii). These results are known, but not easily collected from the different sources in the literature. 
Therefore, a survey on this subject is included in Appendix~\ref{JBappendix}.

In Section \ref{sec.antilattices}, we show that every factor of an atomic JBW-algebra is a disjointness free anti-lattice. In Section \ref{sec.orderdirectsums}, we develop basic theory on  direct sums of pre-Riesz spaces and order direct sums of order unit spaces. Disjointness and bands in order direct sums of order unit spaces that are disjointness free anti-lattices are characterised in Section \ref{sec.disjointnessinsums}. As a  consequence, we obtain a characterisation of disjointness and bands in atomic JBW-algebras. This is used in Section \ref{sec.inverses} to show that disjointness preserving linear bijections with disjointness preserving inverses are exactly the bijections that permute the factors in the direct sum. In Section \ref{sec.centres}, we show that the algebraic centre and the order theoretical centre of a unital JB-algebra are isomorphic as JB-algebras. Consequently, the order theoretical centre is a Riesz space. We introduce in Section \ref{sec.moreonordercentre} a class of order unit spaces, including all finite-dimensional ones, whose order theoretical centre is isomorphic to $\mathbb{R}^n$ for some $n$.

\section{Preliminaries on partially ordered vector spaces}\label{sec.prelimpov}

Let $X$ be  a real vector space containing a  cone $K$, i.e., $K$ is convex, $\lambda K\subseteq K$ for every $\lambda \geq 0$, and $K\cap -K=\{0\}$. The cone $K$ induces a partial order $\leq $ in $X$ by $x\leq y$ if $y-x\in K$. 
We call $(X,K)$ a \emph{partially ordered vector space}. We say that $(X,K)$ is \emph{directed} if $X=K-K$.
The space $(X,K)$ is called \emph{Archimedean} if, for every $x,y\in X$ with $nx\leq y$ for all $n\in \mathbb{N}$, we have $x\leq 0$. A partially ordered vector space $X$ is called \emph{monotone complete} if for any increasing net $(x_i)_i$ in $X$ that is bounded from above the supremum exists in $X$. 

If for every $x,y\in X$ the supremum of $\{x,y\}$ exists, then $X$ is called a \emph{vector lattice} or a \emph{Riesz space}. For further terminology on vector lattices, see \cite{AliBur1985}. We say that $(X,K)$ is an \emph{anti-lattice} if for every $x,y\in X$ the supremum of $\{x,y\}$ exists only if $x$ and $y$ are comparable, that is, $x\le y$ or $x\ge y$. Trivially, if $(X,K)$ is totally ordered, then $X,K)$ is an anti-lattice. Hence $\mathbb{R}$ is both a lattice and an anti-lattice.

A linear subspace $D$ of $X$ is \emph{order dense} in $X$ if, for every $x\in X$, we have 
\[
x=\inf\{d\in D\colon d\geq x\},
\]
and a subspace $Y$ of a partially ordered vector space $X$ is is called \emph{majorizing} in $X$ if for every $x \in X$ there is a $y\in Y$ such that $x \le y$. A linear map $T\colon X\to Y$, where $X$ and $Y$ are partially ordered vector spaces, is called \emph{positive} if for every $x\in X$ with $x\ge 0$ we have $Tx\ge 0$ and $T$ is called \emph{bipositive} if $x\ge 0$ is equivalent to $Tx\ge 0$. 

A partially ordered vector space $X$ is called a \emph{pre-Riesz space} if there is a Riesz space $Y$ and a bipositive linear map $i\colon X\to Y$ such that $i[X]$ is order dense in $Y$. We call $(Y,i)$ a \emph{vector lattice cover} of $X$. An intrinsic definition of pre-Riesz spaces is given by van Haandel in \cite{Haa1993}, see also \cite[Section 2.2]{KalvanGaa2019}.
Note that every directed Archimedean partially ordered vector space is pre-Riesz, and that every pre-Riesz space is directed. Clearly, every Riesz space is pre-Riesz. 
If $(Y,i)$ is a vector lattice cover of a pre-Riesz space $X$ such that no proper Riesz subspace of $Y$ contains $i[X]$, then we call $(Y,i)$ a \emph{Riesz completion} of $X$ and is denoted by $X^\rho$. Such a space is unique up to isomorphism (for details see, e.g., \cite[Section 2.4]{KalvanGaa2019}). 

For $A\subseteq X$, denote 
\[
A^{\mathrm{u}}:=\{x\in X\colon x\ge a \mbox{ for all $a\in A$}\}\quad \mbox{and}\quad  
A^{\mathrm{l}}:=\{x\in X\colon x\le a \mbox{ for all $a \in A$}\}. 
\]
Riesz* homomorphisms are defined in \cite[Definition 5.1 and Corollary 5.4(iv)]{Haa1993} and Riesz homomorphisms %and complete Riesz homomorphisms 
in  \cite{BusRoo1993}.  
\begin{definition}
	Let $X$ and $Y$ be directed partially ordered vector spaces. A linear map $T\colon X\to Y$ is called 
	\begin{itemize}
		\item[(i)] 
		a \emph{Riesz* homomorphism} 
		if,
		for every non-empty finite subset $F$ of $X$, one has 
		\[T\left[F^\mathrm{ul}\right]\subseteq T[F]^{\mathrm{ul}},\]
		\item[(ii)] 
		a \emph{Riesz homomorphism} if, for every $x,y\in X$, one has
		\[
		T\left[\{x,y\}^\mathrm{u}\right]^\mathrm{l}=T[\{x,y\}]^\mathrm{ul}.
		\]
%		\item[-] 
%		a \emph{complete Riesz homomorphism} if, for every nonempty set $A\subset X$, we have
%		\[\inf A=0\ \Longrightarrow \ \inf T[A]=0.\]
	\end{itemize}
\end{definition}
If $X$ and $Y$ are pre-Riesz spaces, then %every complete Riesz homomorphism is a Riesz homomorphism, 
every Riesz homomorphism is a Riesz* homomorphism, and every Riesz* homomorphism is positive, see \cite[Theorem 2.3.19]{KalvanGaa2019}. 
If $X$ and $Y$ are vector lattices, then the notions of a Riesz homomorphism and a Riesz* homomorphism both coincide with the notion of a Riesz homomorphism from vector lattice theory, see, e.g., \cite[Lemma 2.3.2]{KalvanGaa2019}. %Moreover, in this case, $T$ is a complete Riesz homomorphism if and only if $T$ is an order continuous Riesz homomorphism, see \cite[Proposition 1.4.5]{KalvanGaa2019}. 

The following Lipecki-Luxemburg-Schep theorem can be found, e.g., in \cite[Theorem 2.1.17]{KalvanGaa2019}.
\begin{theorem}\label{thm:LLS}
	Let $Y$ be a Riesz space, let $Z$ be a Dedekind complete Riesz space, and let $D$ be a majorizing Riesz subspace of $Y$. If $h\colon D\to Z$ is a Riesz homomorphism, then there exists a Riesz homomorphism $H\colon Y\to Z$ that extends $h$.
\end{theorem}

The subsequent theorem is due to van Haandel, see, e.g.,
\cite[Theorem 2.4.11]{KalvanGaa2019}. In this section, $\circ$ denotes composition.
\begin{theorem}\label{the:vanHaandel}
	Let $X_1$ and $X_2$ be pre-Riesz spaces and let $(Y_1, i_1)$ and $(Y_2,i_2)$ be vector lattice covers, respectively. Let $h\colon X_1\to X_2$ be a linear map. 
	\begin{itemize}	
		\item[(i)] If there exists a Riesz homomorphism $\hat{h}\colon Y_1\to Y_2$ such that $\hat{h}\circ i_1=i_2\circ h$, then $h$ is a Riesz* homomorphism.
		\item[(ii)] If $(Y_1, i_1)$ is the Riesz completion of $X_1$ and $h$ is a Riesz* homomorphism, then there exists a unique Riesz homomorphism $\hat{h}\colon Y_1\to Y_2$ with $\hat{h}\circ i_1=i_2\circ h$. 
	\end{itemize}		 
\end{theorem}

An element $u\in K$ is said to be an \emph{order unit} if for every $x \in X$ there is a $\lambda>0$ such that $-\lambda u \le x \le \lambda u$. If $(X,K)$ is an Archimedean partially ordered vector space with order unit $u$, it can be equipped with the \emph{order unit norm} which is defined by $$\|x\|_u:=\inf\{\lambda>0\colon -\lambda u\leq x\leq \lambda u\}$$ for $x\in X$, see, e.g., \cite[Section 1.5.3]{KalvanGaa2019}. In this case, the triple $(X,K,u)$ is called an \emph{order unit space}. Every order unit space is a pre-Riesz space. In the setting of order unit spaces, we recall characterisations of functionals that are Riesz homomorphisms or Riesz* homomorphisms, respectively, and construct a vector lattice cover with the pointwise partial ordering.
The functional representation of $X$ is given by means of the \emph{state space}, $\Sigma_X$, which is the weak* compact convex set 
\begin{equation}%\label{equ:Sigma}
\Sigma_X:=\{\varphi\in X^*\colon \varphi[K]\subseteq [0,\infty),\ \varphi(u)=1\}
\end{equation}
by the Banach-Alaoglu theorem, and the set $\Lambda_X$ of the  extreme points of $\Sigma_X$, which exist by the Krein-Milman theorem. The weak* closure $\overline{\Lambda}_X$ of $\Lambda_X$ in $\Sigma_X$ is a compact Hausdorff space, and the map \begin{equation}
\label{equ:func_repr_OUS}
\Phi_X\colon X\to \mathrm{C}(\overline{\Lambda}_X),\quad x\mapsto (\varphi\mapsto\varphi(x)),\end{equation} is a bipositive linear map, and hence injective (for details, see, e.g., \cite[Section 2.5]{KalvanGaa2019}). %In \cite{KalLemGaa2014},
%it is shown that
Moreover, $(\mathrm{C}(\overline{\Lambda}_X),\Phi_X)$ is a vector lattice cover of $X$ 
\cite[Theorem 2.5.9]{KalvanGaa2019}.
We recall the statement in \cite[Proposition 2.5.5]{KalvanGaa2019}.
\begin{proposition}\label{pro:RieszhomOUS}
	Let $(X,K,u)$ be an order unit space and let $\varphi\in\Sigma_X$. 
	\begin{itemize}
		\item[(i)] One has $\varphi\in\Lambda_X$ if and only if $\varphi$ is a Riesz homomorphism.
		\item[(ii)] One has $\varphi\in\overline{\Lambda}_X$ if and only if $\varphi$ is a Riesz* homomorphism.
	\end{itemize}
\end{proposition} 

Recall that two elements $x$ and $y$ in a Riesz space are said to be disjoint if $|x|\wedge|y|=0$. This notion is generalised to pre-Riesz spaces as follows. Two elements $x$ and $y$ in a pre-Riesz space $(X,K)$ are called \emph{disjoint}, denoted $x\perp y$, if $\{x+y,x-y\}^{\mathrm{u}}=\{x-y,-x+y\}^{\mathrm{u}}$. The disjoint complement of a set $A\subseteq X$ is denoted by $A^{\mathrm{d}}$. If $(Y,i)$  is  a vector lattice cover of $X$, then $x\perp y$ if and only if $i(x)\perp i(y)$, see, e.g., \cite[Proposition 4.1.4]{KalvanGaa2019}. Anti-lattices can be characterised by means of disjointness. A pre-Riesz space $(X,K)$ is an anti-lattice if and only if there are no non-trivial positive disjoint elements in $X$, see \cite[Theorem 14]{KalLemGaa2014}. We call $X$ \emph{disjointness free} if there are no non-trivial disjoint elements in $X$. Clearly, every disjointness free partially ordered vector space is an anti-lattice. In \cite{KalLemGaa2014}, an example of an anti-lattice that is not disjointness free can be found.

%Let $X$ and $V$ be pre-Riesz spaces. A linear map $T\colon X\to V$ is called 
%\emph{disjointness preserving} if, for every $x,y\in X$ with $x\perp y$, one has $T(x)\perp T(y)$. If $T$ is a  Riesz* homomorphism, then $T$ is a positive disjointness preserving operator, see \cite[Theorem 5.1.12]{KalvanGaa2019}.
A set $B\subseteq X$ is called a \emph{band} if $B=B^{\mathrm{dd}}$. Bands in pre-Riesz spaces are linear subspaces, see \cite[Proposition 4.1.5(ii)]{KalvanGaa2019}. Many examples are given in \cite[Chapter 4]{KalvanGaa2019}. It is straightforward to verify that the intersection of two bands is a band. In contrast to vector lattices, there may exist bands that are not directed in pre-Riesz spaces. 

A projection $P$ in $V$ is called an \emph{order projection} if both $P$ and $I-P$ are positive operators, where $I$ denotes the identity operator. If $V$ is a pre-Riesz space, then \cite[Proposition 3.1]{Glueck2021} yields that a projection $P$ in $V$ is an order projection if and only if $P$ is a \emph{band projection}, which means that the range and kernel of $P$ both are bands in $V$. The range of a band projection is called a \emph{projection band}. If $P$ and $Q$ are two band projections in $V$, then $PQ$ is a band projection in $V$, see \cite[Proposition 3.6]{Glueck2021}.

Direct sums will play a crucial role in later sections. Let $\mathcal{I}$ be a non-empty set and let $((V_i,C_i,u_i))_{i\in \mathcal{I}}$ be a collection of order unit spaces. We define the \emph{order direct sum} to be the vector space 
\begin{equation}\label{eq.Visorderdirectsum}
\bigoplus_{i\in \mathcal{I}} V_i:=\left\{ i\mapsto v_i\colon \mathcal{I}\to\bigcup_{i\in \mathcal{I}} V_i\colon\, v_i\in V_i\mbox{ for every }i\in \mathcal{I}\mbox{ and } \sup_{i\in \mathcal{I}} \|v_i\|_{u_i}<\infty\right\}
\end{equation}
with the cone $\{v\in \bigoplus_{i\in \mathcal{I}} V_i\colon\, v(i)\in C_i\mbox{ for every }i\in \mathcal{I}\}$. Then $\bigoplus_{i\in \mathcal{I}} V_i$ is an Archimedean directed partially ordered vector space with order unit $i\mapsto u_i$, which we denote by $u$. Note that for every $v\in \bigoplus_{i\in \mathcal{I}} V_i$ we have that 
\begin{equation}\label{eq.normdirectsum}
\|v\|_u=\sup_{\alpha\in \mathcal{I}} \|v(i)\|_{u_i}.
\end{equation} 
Let $\mathcal{J}$ be a non-empty subset of $\mathcal{I}$. We define $\Phi_\mathcal{J}\colon \bigoplus_{j\in \mathcal{J}} V_j\to \bigoplus_{i\in \mathcal{I}} V_i$ by $\Phi_\mathcal{J}(w):=v$, where $v_i=w_i$ for every $i\in \mathcal{J}$ and $v_i=0$ otherwise. Clearly, $\Phi_\mathcal{J}$ is a bipositive linear map. If $\mathcal{J}=\{j\}$, then we write $\Phi_j$ instead of $\Phi_\mathcal{J}$.

If $((W_i,K_i,w_i))_{i\in \mathcal{I}}$ is another family of order unit spaces and for every $i\in \mathcal{I}$ we have a linear map $T_i\colon V_i\to W_i$ such that for every $v\in \bigoplus_{i\in \mathcal{I}} V_i$ the map $i\mapsto T_i v(i)$ from $\mathcal{I}$ to $\bigcup_{i\in \mathcal{I}} W_i$ belongs to $\bigoplus_{i\in \mathcal{I}} W_i$, then we denote this map by $\bigoplus_{i\in \mathcal{I}} T_i$. 

Let $(V,C,u)$ be an order unit space with $V\neq \{0\}$. Then $V$ is called \emph{irreducible} if for every collection $((V_i,C_i,u_i))_{i\in \mathcal{I}}$ of order unit spaces such that $V$ is isomorphic to $\bigoplus_{i\in \mathcal{I}} V_i$ as order unit spaces, there exists $i\in \mathcal{I}$ such that $V_j=\{0\}$ for all $j\in \mathcal{I}\setminus\{i\}$. Otherwise, $V$ is called \emph{reducible}. If there exists an order projection $P$ in $V$ with $P\neq 0$ and $P\neq I$, then 
\[
(P[V],P[C],Pu)\quad \mbox{and}\quad ((I-P)[V], (I-P)[C],(I-P)u) 
\]
are non-trivial order unit spaces and $V$ is as order unit space isomorphic to the order direct sum $P[V]\oplus (I-P)[V]$. Hence, $V$ is reducible. 

Let $(V,C)$ be a directed partially ordered vector space. The \emph{order theoretical centre} of $V$ is the set 
\[\mathrm{E}(V):=\{T\colon V\to V\colon\, T\mbox{ is linear and there exists }\lambda\ge 0\mbox{ such that }-\lambda I\le T\le \lambda I\},\]
which is a partially ordered vector space with order unit $I$. If $V$ is Archimedean, then so is $\mathrm{E}(V)$. In that case, Buck shows in \cite{Buck} that the restriction of each $T\in \mathrm{E}(V)$ to any subspace of $V$ that is an order unit space corresponds to a multiplication operator on the functional representation of that subspace. As a consequence, he obtains that $\mathrm{E}(V)$ is commutative under composition. In our analysis, we also need properties of the operator norm on $\mathrm{E}(V)$, stated in \Cref{cor.normsoncentresame}(i),(ii), and (iv) below. As it is little extra work, we will reprove Buck's result for order unit spaces and make the treatment of $\mathrm{E}(V)$ self-contained. 

A norm $\left\|\cdot\right\|$ on $V$ is called \emph{semimonotone} if there exists a constant $\kappa$ such that for all $v,w\in C$ with $v\le w$ we have $\|v\|\le \kappa\|w\|$. The norm is called \emph{regular} if $\|v\|=\inf\{\|w\|\colon\, -w\le v\le w\}$ for every $v\in V$. If $\left\|\cdot\right\|$ is a semimonotone norm on $V$ such that $C$ is closed and $V$ is complete, then the norm is equivalent to a regular norm and then every element of $\mathrm{E}(V)$ is a bounded operator; see \cite[Corollary 3.4.13 and Lemma 5.4.1]{KalvanGaa2019}. In that case, we obtain
\begin{equation}\label{eq.ordercentrebounded}
\mathrm{E}(V)=\{T\in \mathrm{B}(V)\colon\, \mbox{ there exists }\lambda\ge 0\mbox{ such that } -\lambda I\le T\le \lambda I\},
\end{equation}
where $\mathrm{B}(V)$ denotes the vector space of all bounded linear operators on $V$. Note that for a general order unit space, the order theoretical centre $\mathrm{E}(V)$ can be equipped with both the operator norm and the order unit norm. In \Cref{cor.normsoncentresame} below, we show that these norms coincide.

For a linear subspace $X\subseteq \mathrm{C}(\Omega)$, where $\Omega$ is a compact Hausdorff space, we denote the point evaluation at a point $w\in \Omega$ by $\delta_w$, that is, $\delta_w(x)=x(w)$ for every $x\in X$. For a function $f\in \mathrm{C}(\Omega)$, we define the corresponding multiplication operator $M_f$ on $\mathrm{C}(\Omega)$ by
\[M_f(g)(w)=f(w)g(w),\ w\in \Omega,\ g\in \mathrm{C}(\Omega).\]
Note that $\|M_f\|=\|f\|_\infty$. We show that the operators in the order theoretical centre of an order unit space correspond to multiplication operators on the functional representation.

\begin{proposition}\label{pro.centralismultiplication}
Let $\Omega$ be a compact Hausdorff space and let $X\subseteq \mathrm{C}(\Omega)$ be a linear subspace containing the constant one function $\one$. Assume that
\[\Omega_0:=\{w\in \Omega\colon\, \delta_w\colon X\to \mathbb{R}\mbox{ is a Riesz homomorphism}\}\]
is dense in $\Omega$. Let $T\colon X\to X$ be a linear map with $0\le T\le I$. Then for every $x\in X$ and for every $w\in \Omega$ we have
\[(Tx)(w)=(T\one)(w) x(w).\]
\end{proposition}
\begin{proof}
Let $w\in \Omega_0$. Define $\varphi\colon X\to\mathbb{R}$ by
\[\varphi(x):=(Tx)(w),\ x\in X.\]
Then $\varphi\colon X\to\mathbb{R}$ is linear and for every $x\in X$ with $x\ge 0$ we have $0\le Tx\le Ix=x$, so
\[\varphi(x)=(Tx)(w)\ge 0\mbox{ and } \varphi(x)=(Tx)(w)\le x(w)=\delta_w(x).\]
Hence $0\le \varphi\le\delta_w$. Due to \Cref{pro:RieszhomOUS}, we have $\delta_w\in\Lambda_X$, hence $\delta_w$ is an extreme point of $\Sigma_X$. Then $\delta_w$ is extremal in the dual cone of $X$ by \cite[Lemma 1.5.19]{KalvanGaa2019}. Therefore, there exists $\lambda\in [0,1]$ such that $\varphi=\lambda\delta_w$. In particular, $\lambda=\lambda\delta_w(\one)=\varphi(\one)=(T\one)(w)$. Thus, for every $x\in X$ we have
\[(Tx)(w)=\varphi(x)=\lambda \delta_w(x)=(T\one)(w)x(w).\]
Next, let $x\in X$. As $Tx$, $T\one$, and $x$ are continuous on $\Omega$ and $\Omega_0$ is dense in $\Omega$, it follows that $(Tx)(w)=(T\one)(w)x(w)$ for every $w\in \Omega$.
\end{proof}

The next result and \Cref{cor.normsoncentresame}(iii) are due to Buck in \cite{Buck}.

\begin{theorem}\label{the.centreandmultiplication}
Let $(V,C,u)$ be an order unit space and let $\Phi\colon V\to \mathrm{C}(\overline{\Lambda})$ be its functional representation. 
\begin{itemize}
\item[(i)] We have
$\mathrm{E}(V)=\{T\colon V\to V\colon\, \Phi\circ T= M_{\Phi(Tu)}\circ \Phi\}$.
\item[(ii)] For every $S,T\in \mathrm{E}(V)$, we have that $\Phi\circ(S\circ T)=M_{fg}\circ \Phi$, where $f=\Phi(Su)$ and $g=\Phi(Tu)$. 
\end{itemize}
\end{theorem}
\begin{proof}
$(i)$ Let $T\in \mathrm{E}(V)$. Let $\alpha>0$ be such that $-I\le \alpha T\le I$ and denote $S:= I-\frac{1}{2}\alpha T$. Then $0\le S\le I$. Since $\Phi$ is an order isomorphism from $V$ onto the subspace $X:=\Phi[V]$ of $\mathrm{C}(\overline{\Lambda})$, we have that $\Lambda$ equals 
\begin{align*}
&\{w\in \overline{\Lambda}\colon\, \delta_w\colon X\to\mathbb{R}\mbox{ is a Riesz homomorphism}\} =\{w\in\overline{\Lambda}\colon\, w\colon V\to\mathbb{R}\mbox{ is a Riesz homomorphism}\},
\end{align*}
which is dense in $\overline{\Lambda}$. \Cref{pro.centralismultiplication} yields for every $x\in X$ and every $w\in\overline{\Lambda}$ that 
\[(\Phi\circ S\circ \Phi^{-1}x)(w)=(\Phi\circ S\circ \Phi^{-1}\one)(w) x(w),\]
so that for every $v\in V$ we have
\[((\Phi\circ S) v)(w)=((\Phi\circ S)u)(w)(\Phi v)(w).\]
Hence, $((\Phi-\frac{1}{2}\alpha \Phi\circ T)v)(w)=(\one-\frac{1}{2}\alpha\Phi\circ Tu)(w)(\Phi v)(w)$, which yields that $\Phi\circ Tv=M_{\Phi(Tu)}( \Phi(v))$.

Conversely, let $T\colon V\to V$ be such that $\Phi\circ T=M_{\Phi(Tu)}\circ\Phi$. There is $\alpha>0$ with $-\alpha \one\le \Phi(Tu)\le \alpha \one$. Then for every $v\in C$ we have $-\alpha\Phi(v)\le \Phi(Tv)\le \alpha\Phi (v)$, so that $-\alpha I\le T\le \alpha I$.

$(ii)$ We have $\Phi\circ (S\circ T) = (\Phi\circ S)\circ T=(M_f\circ\Phi)\circ T= M_f\circ(\Phi\circ T) = M_f\circ(M_g\circ \Phi) = (M_f\circ M_g)\circ\Phi= M_{fg}\circ\Phi$.
\end{proof}

\begin{corollary}\label{cor.normsoncentresame}
Let $(V,C,u)$ be an order unit space. 
\begin{itemize}
\item[(i)]  $\mathrm{E}(V)$ is a subspace of $\mathrm{B}(V)$ and the operator norm and the order unit norm induced by $I$ coincide.
\item[(ii)] $\mathrm{E}(V)$ is a closed subspace of $\mathrm{B}(V)$.
\item[(iii)] $\mathrm{E}(V)$ with composition is a commutative associative algebra.
\item[(iv)] For every $S,T\in \mathrm{E}(V)$ we have 
\[\|ST\|\le \|S\| \|T\|,\ \|T^2\|= \|T\|^2, \mbox{ and } \|T^2\|\le \|S^2+T^2\|.\]
\end{itemize}
\end{corollary}
\begin{proof}
$(i)$ We use the functional representation of $V$ and \Cref{the.centreandmultiplication}. Let $T\in \mathrm{E}(V)$ , denote $f:=\Phi(Tu)$, and put $X:=\Phi[V]$. Then $\Phi\circ T= M_{f}\circ \Phi$. Since $\one=\Phi(u)\in \Phi[V]$, both the operator norm and the order unit norm of $M_f$ on $X$ equal $\|f\|_\infty$. Indeed, for every $x\in X$ we have $\|M_fx\|_\infty\le \|f\|_\infty \|x\|_\infty$ and $\|M_f\one\|_\infty=\|f\|_\infty$, hence $\|M_f\|=\|f\|_\infty$. For every $\alpha\ge 0$ and every $x\in X$ with $x\ge 0$ we have that $-\alpha Ix\le M_fx\le \alpha I x$ if and only if $-\alpha x(w)\le f(w)x(w)\le \alpha x(w)$ for all $w\in \overline{\Lambda}$, hence $-\alpha I\le M_f\le \alpha I$ if and only if $-\alpha\one\le f\le \alpha\one$. Thus, $\|M_f\|_I=\|f\|_\infty$.

As $\Phi\colon V\to X$ is an order isomorphism, for all $x \in C$ we have $-\lambda x \le Tx \le \lambda x$ if and only if $-\lambda \Phi(x) \le M_f\Phi(x)\le \lambda\Phi(x)$. Hence $\|T\|_I=\|M_f\|_I$, where $\|T\|_I$ is the order unit norm in $\mathrm{E}(V)$ and $\|M_f\|_I$ is the order unit norm in $\mathrm{E}(\mathrm{C}(\Omega))$. Also, $\Phi\colon V\to X$ is an isometry, so $T$ is bounded and $\|T\|=\|M_f\|$. Thus, $\|T\|_I=\|M_f\|_I=\|f\|_\infty=\|M_f\|=\|T\|$.

$(ii)$ Let $(T_n)_n$ be a sequence in $\mathrm{E}(V)$ and let $T\in \mathrm{B}(V)$ be such that $\|T_n-T\|\to 0$. For every $n\in\mathbb{N}$ denote $f_n:=\Phi(T_nu)$. By $(i)$, we have $\|f_n-f_m\|_\infty=\|M_{f_n}-M_{f_m}\|=\|T_n-T_m\|$ for every $n,m\in\mathbb{N}$, hence $(f_n)_n$ is a Cauchy sequence in $\mathrm{C}(\overline{\Lambda})$. Let $f\in \mathrm{C}(\overline{\Lambda})$ be the limit of $(f_n)_n$. Next we show that $\Phi\circ T=M_f\circ \Phi$. Indeed, let $v\in V$. As $\Phi\circ T_n=M_{f_n}\circ \Phi$ for every $n\in\mathbb{N}$, we have 
\begin{align*}
\|(\Phi\circ T)v-(M_f\circ \Phi)v\|=& \|(\Phi\circ T)v-(\Phi\circ T_n)v+(M_{f_n}\circ \Phi)v-(M_f\circ \Phi)v\|\\
&\le \|\Phi\|\|T_n-T\|\|v\|+\|f_n-f\|_\infty\|v\|\to 0,
\end{align*}
hence $\Phi\circ T=M_f\circ \Phi$. Hence $\Phi(Tu)=f$ and $T \in \mathrm{E}(V)$ by \Cref{the.centreandmultiplication}$(i)$. 

$(iii)$ Let $S,T\in \mathrm{E}(V)$ and let $\alpha,\beta\in\mathbb{R}$ be such that $-\alpha I\le S\le \alpha I$ and $-\beta I\le T\le \beta I$. As $\beta I+T$ is positive, we obtain $S(\beta I +T)\le \alpha (\beta I +T)$, hence 
\[ST\le \alpha \beta I+\alpha T -\beta S\le 3\alpha\beta I.\]
Similarly, $ST\ge -3\alpha\beta I$, hence $ST\in \mathrm{E}(V)$. 

Let $f=\Phi(Su)$ and $g=\Phi(Tu)$. By \Cref{the.centreandmultiplication}$(ii)$, we have 
\[\Phi\circ(ST)=M_{fg}\circ\Phi=M_{gf}\circ \Phi=\Phi\circ(TS),\]
hence $ST=TS$.

$(iv)$ Clearly, the operator norm is submultiplicative. With the aid of \Cref{the.centreandmultiplication}$(ii)$, we obtain \[\|T^2\|=\|M_{g^2}\|=\|g^2\|_\infty=\|g\|_\infty^2 = \|T\|^2.\] Similarly, $\|T^2\|=\|g^2\|_\infty\le \|f^2+g^2\|_\infty=\|M_{f^2+g^2}\|=\|M_f^2+M_g^2\|=\|S^2+T^2\|$, since $(M_f^2+M_g^2)\circ \Phi= \Phi\circ (S^2+T^2)$.
\end{proof}

\begin{remark}\label{rem.coneQ}
Let $(V,C,u)$ be an order unit space and $Q=\{T^2\colon\, T\in \mathrm{E}(V)\}$. By \Cref{cor.normsoncentresame}$(iii)$, we have $Q\subseteq \mathrm{E}(V)$. Moreover, $Q\subseteq\mathrm{E}(V)_+$. Indeed, let $T\in \mathrm{E}(V)$ and $f=\Phi(Tu)$. As in \Cref{cor.normsoncentresame}$(iv)$, we have $\Phi\circ T^2= M_{f^2}\circ\Phi$. Since $\Phi$ is an order isomorphism, it follows that $T^2$ is positive. In general, the sets $Q$ and $\mathrm{E}(V)_+$ differ. For an example, let $V:=\mathrm{Pol}[0,1]$ be the subspace of $\mathrm{C}([0,1])$ consisting of all polynomials. Then $V$ is an order unit space with the constant function $\one$ as order unit. Since $V$ is order dense in $\mathrm{C}([0,1])$, the space $\mathrm{C}([0,1])$ (with identity as embedding map) is the functional representation of $V$. Further, $V$ is an associative subalgebra of $\mathrm{C}([0,1])$. Let $f(t)=t^2+1$. The multiplication operator $M_f$ maps $V$ into $V$. Also, $0\le M_f\le 2 I$, hence $M_f\in \mathrm{E}(V)_+$. Suppose that there exists $T\in \mathrm{E}(V)$ such that $T^2=M_f$. Let $g=T\one$. Then, by \Cref{the.centreandmultiplication}, we have $T=M_g$, so that $M_{g^2}=T^2=M_f$, hence $g^2=f$ and thus $g$ is not a polynomial. Then $T\one=M_g\one$ is not in $V$, which yields a contradiction. Note that we also obtain that $Q$ is not a cone, as $M_\one,M_{t\mapsto t^2}\in Q$. 
\end{remark}

\begin{corollary}\label{cor.ordercentreofsubalgebra}
Let $(V,C,u)$ be an order unit space and let $\Phi\colon V\to \mathrm{C}(\overline{\Lambda})$ be its functional representation. If $\Phi[V]$ is a subalgebra of $\mathrm{C}(\overline{\Lambda})$, then $\mathrm{E}(V)=\{\Phi^{-1}\circ M_f\circ \Phi\colon\, f\in \Phi[V]\}$.
\end{corollary}

%\begin{example}
%Let $V=C^1[0,1]$ with pointwise order and order unit $\one$. The inclusion of $V$ in $C[0,1]$ is the functional representation of $V$. As $V$ is a subalgebra of $C[0,1]$, Corollary \ref{cor.ordercentreofsubalgebra} yields that $E(V)=\{M_f\colon\, f\in C^1[0,1]\}$. In particular, $E(V)$ is not a lattice.
%\end{example}

\section{Preliminaries on JB-algebras}

A \emph{Jordan algebra} $(A, \circ)$ is a commutative, not necessarily associative algebra such that
\[
x \circ (y \circ x^2) = (x \circ y) \circ x^2 \mbox{\quad  for all }x,y \in A.
\]
A \emph{JB-algebra} $A$ is a normed, complete Jordan algebra over the scalar field $\mathbb{R}$ satisfying
\begin{align*}
\norm{x \circ y} &\leq \norm{x}\norm{y}, \\
\norm{x^2} &= \norm{x}^2, \\
\norm{x^2} &\leq \norm{x^2 + y^2}
\end{align*}
for all $x,y \in A$. If $A$ is finite-dimensional and there is an inner product on $A$ such that $\langle x\circ y,z \rangle = \langle x,y\circ z \rangle$ for all $x,y,z \in A$, then $A$ is a so called \emph{Euclidean Jordan algebra}, see \cite[Chapter III]{Faraut-Koranyi}. As mentioned in the introduction, an important example of a JB-algebra is the set of self-adjoint elements of a $C^*$-algebra equipped with the Jordan product $x \circ y := \frac{1}{2}(xy + yx)$.

The elements $x,y \in A$ are said to \emph{operator commute} if $x \circ (y \circ z) = y \circ (x \circ z)$ for all $z \in A$. An element $x\in A$ is said to be \emph{central} if it operator commutes with all elements of $A$. The \emph{algebraic centre} of $A$, denoted by $\mathrm{Z}(A)$, consists of all elements that operator commute with all elements of $A$, and it is an associative subalgebra of $A$. In the remainder of this paper, it will be assumed that all JB-algebras have an algebraic unit $e$ and then $e\in \mathrm{Z}(A)$. The following representation theorem can be found in \cite[Theorem~3.2.2]{Hanche-Olsen-Stormer}. Here $\mathrm{C}(\Omega)$ is endowed with supremum-norm and pointwise multiplication.

\begin{theorem}\label{the.representationJB}
Every associative unital JB-algebra is isometrically isomorphic as a JB-algebra to $\mathrm{C}(\Omega)$ for some compact Hausdorff space $\Omega$.
\end{theorem}

\begin{corollary}\label{cor.algcentreisCspace}
The algebraic centre $\mathrm{Z}(A)$ of a unital JB-algebra is isometrically isomorphic as a JB-algebra to $\mathrm{C}(\Omega)$ for some compact Hausdorff space $\Omega$.
\end{corollary}

The \emph{spectrum} $\sigma(x)$ of $x \in A$ is defined to be the set of $\lambda \in \R$ such that $x - \lambda e$ is not invertible in JB$(x,e)$, the JB-subalgebra of $A$ generated by $x$ and $e$, see \cite[Section~3.2.3]{Hanche-Olsen-Stormer}. Furthermore, there is a continuous functional calculus, which means that there exists an isometric JB-algebra isomorphism from  JB$(x,e)$ onto $\mathrm{C}(\sigma(x))$, see \cite[Corollary~1.19]{Alfsen}. The cone of elements with non-negative spectrum is denoted by $A_+$, and equals the set of squares by the functional calculus, and its interior $A_+^\circ$ consists of all elements with strictly positive spectrum. This cone turns $A$ into an order unit space with order unit $e$, that is,
\[ \norm{x} = \inf \{ \lambda > 0: -\lambda e \leq x \leq \lambda e \}. \]
The \emph{Jordan triple product} $\{ \cdot, \cdot, \cdot \}$ is defined as
\[ \{x,y,z \} := (x \circ y) \circ z + (z \circ y) \circ x - (x \circ z) \circ y \]
for $x,y,z \in A$. The linear map $U_x \colon A \to A$ defined by $U_x y := \{x,y,x\}$ will play an important role and is called the \emph{quadratic representation} of $x$. It is always a positive map by \cite[Theorem~1.25]{Alfsen}. In case $x$ is invertible, it follows that $U_x$ is an automorphism of the cone $A_+$ and its inverse is $U_{x^{-1}}$ by \cite[Lemma~1.23]{Alfsen} and \cite[Theorem~1.25]{Alfsen}. A {\it state} $\varphi$ of $A$ is a positive linear functional on $A$ such that $\varphi(e)=1$. The set of states on $A$ is called the \emph{state space} of $A$. The extreme points of the state space are referred to as \emph{pure states} on $A$ (cf. \cite[A~17]{Alfsen}). In the notation introduced for pre-Riesz spaces in the previous section, the state space is denoted by $\Sigma_A$ and the pure states by $\Lambda_A$, which would be the Riesz homomorphisms from the functional representation of $A$ given in \eqref{equ:func_repr_OUS}.  

A \emph{JBW-algebra} $M$ is the Jordan analogue of a von Neumann algebra: it is a monotone complete JB-algebra with unit $e$ and a separating set of normal states, or equivalently, a JB-algebra that is a dual space. A state $\varphi$ on $M$ is  said to be {\it normal} if for any bounded increasing net $(x_i)_{i}$ with supremum $x$ we have $\varphi(x_i)\rightarrow\varphi(x)$. The (convex) set of normal states on $M$ is called the {\it normal state space} of $M$. The topology on $M$ defined by the duality of $M$ and the normal state space of $M$ is called the {\em $\sigma$-weak topology}. That is, we say a net $(x_i)_{i}$ converges $\sigma$-weakly to $x$ if $\varphi(x_i)\to \varphi(x)$ for all normal states $\varphi$ on $M$. The Jordan multiplication on a JBW-algebra is separately $\sigma$-weakly continuous in each variable and jointly $\sigma$-weakly continuous on bounded sets by \cite[Proposition~2.4]{Alfsen} and \cite[Proposition~2.5]{Alfsen}. Furthermore, for any $x$ the corresponding quadratic representation $U_x$ is $\sigma$-weakly continuous by \cite[Proposition~2.4]{Alfsen}. If $A$ is a JB-algebra, then one can extend the Jordan product uniquely to its bidual $A^{**}$ turning $A^{**}$ into a JBW-algebra, see \cite[Corollary~2.50]{Alfsen}. 

An element $p$ in a JBW-algebra $M$ is a {\it projection} if $p^2=p$. For a projection $p\in M$, the {\it orthogonal complement}, $e-p$, will be denoted by $p^\perp$ and a projection $q$ is {\it orthogonal} to $p$ precisely when $q\le p^\perp$, see \cite[Proposition~2.18]{Alfsen}. In each JBW-algebra $M$, the spectral theorem \cite[Theorem~2.20]{Alfsen} holds, which implies, in particular, that the linear span of projections is norm dense in $M$, see \cite[Proposition~4.2.3]{Hanche-Olsen-Stormer}. 
%\begin{comment}
%The collection of projections forms a complete orthomodular lattice by %\cite[Proposition~2.25]{Alfsen}, which means in particular that every set of projections has %a supremum. We remark that these sets of projections need not have a supremum in $M$. 
%\end{comment}

Let $(V_i)_{i\in \mathcal{I}}$ be a family of JBW-algebras with units $u_i$. The \emph{algebraic direct sum} of $(V_i)_{i\in \mathcal{I}}$ is the vector space given by \eqref{eq.Visorderdirectsum} endowed with the norm given by \eqref{eq.normdirectsum} and componentwise multiplication. According to \cite[Definition~2.42]{Alfsen}, the algebraic direct sum $\bigoplus_{i\in \mathcal{I}} V_i$ is a JBW-algebra. If $V_i$ is atomic for every $i\in\mathcal{I}$, then the direct sum $\bigoplus_{i\in \mathcal{I}} V_i$ is atomic, as well. As the positive elements are the squares, the algebraic direct sum and the order direct sum of JBW-algebras coincide. 

Any central projection $p$, i.e., $p$ is a projection that is also a central element, decomposes the JBW-algebra $M$ as an algebraic direct sum of JBW-subalgebras such that $M=U_pM\oplus U_{p^\perp} M$, see \cite[Proposition~2.41]{Alfsen}. If a JBW-algebra $M$ has trivial algebraic centre $\mathbb{R}e$, $M$ is called a \emph{factor}. A minimal element in the set of all non-zero projections of a JBW-algebra is called an \emph{atom}. A JBW-algebra in which every non-zero projection dominates an atom is called \emph{atomic}. Our investigation of atomic JBW-algebras relies on the representation given in \Cref{the.atomicdecomp} and the list (i)--(iii) (above that theorem) of atomic JBW-algebras that are factors. Those factors are discussed in detail in Appendix \ref{JBappendix}. It is shown there that they are indeed JBW-algebras and factors. Moreover the states and pure states are described. The latter are used in the subsequent section to show that the factors are disjointness free anti-lattices.

%By \cite[Proposition~3.45]{Alfsen}, every atomic JBW-algebra is a direct sum of so-called \emph{type I factors}. The type I factors have been classified and are up to JBW-algebra isomorphism either a spin factor, or $\mathrm{Mat}_3(\mathbb{O})_{sa}$, the self-adjoint $3\times 3$ matrices with octonionic entries, or of the form $B(H)_{sa}$, the bounded self-adjoint operators on $H$ where $H$ is a real, complex, or quaternionic Hilbert space of dimension $d \geq 3$,  see \cite[Theorem~3.39]{Alfsen}. These type I factors will be studied in more detail in the following subsections separately. The results stated there are known in the theory of JB-algebras, but are outlined in this paper for the convenience of the reader.   

\section{Factors of atomic JBW-algebras are anti-lattices}\label{sec.antilattices}
Recall that every unital JB-algebra $A$ is an order unit space. Therefore, its functional representation $(\mathrm{C}(\overline{\Lambda_A}),\Phi_A)$ given by \eqref{equ:func_repr_OUS} is a vector lattice cover. All elements of $\overline{\Lambda_A}$ are states of $A$, and the elements of $\Lambda_A$ are precisely the pure states. Disjointness of two elements $a,b\in A$ is equivalent to disjointness of $\Phi_A(a)$ and $\Phi_A(b)$ in $\mathrm{C}(\overline{\Lambda_A})$, which is pointwise disjointness on $\overline{\Lambda_A}$. To show that $a$ and $b$ are not disjoint, it suffices to find one element $\phi\in\overline{\Lambda_A}$ such that $\phi(a)=(\Phi_A(a))(\phi)\neq 0$ and, similarly, $\phi(b)\neq 0$.

\subsection{Disjointness in $\mathrm{B}(\mathcal{H}_q)_\mathrm{sa}$ and $\mathrm{B}(H)_\mathrm{sa}$}
 
We will use vector states to show that there are no non-zero disjoint operators in $\mathrm{B}(\mathcal{H}_q)_\mathrm{sa}$. For a normalised vector $v\in \mathcal{H}_q$ the corresponding vector state is given by $\varphi_v(T):=\langle Tv,v\rangle$, $T\in \mathrm{B}(\mathcal{H}_q)$. By \Cref{L:vector states are pure states}, vector states of $\mathrm{B}(\mathcal{H}_q)_\mathrm{sa}$ are pure states. 

%If we denote by $K$ the state space of $B(\mathcal{H}_q)_{sa}$ and by $\mathrm{ext} K$ the pure states in $K$, then the ambient Riesz space of $B(\mathcal{H}_q)_{sa}$ are the continuous functions on the weak*-closure $\Lambda$ of $\mathrm{ext} K$ in $K$, so $\Phi\colon B(\mathcal{H}_q)_{sa} \to C(\Lambda)$ with $\Phi(T)(\varphi):=\varphi(T)$ turns $B(\mathcal{H}_q)_{sa}$ into a pre-Riesz space. It follows that $T,S$ are disjoint in $B(\mathcal{H}_q)_{sa}$ if and only if $\Phi(T), \Phi(S)$ are disjoint in $C(\Lambda)$. 

\begin{proposition}\label{P: disjoint in B(H_q)}
There are no non-zero disjoint operators in $\mathrm{B}(\mathcal{H}_q)_\mathrm{sa}$.
\end{proposition}

\begin{proof}
Let $S$ and $T$ be non-zero. Then there are normalised vectors $v$ and $w$ such that $\langle Tv,v \rangle\neq 0$ and $\langle Sw,w \rangle\neq 0$ by \Cref{L:spectrum in numerical range} and the functional calculus \cite[Corollary~1.19]{Alfsen}. Consider the corresponding vector states $\varphi_v$ and $\varphi_w$.  If either $\langle Tw,w \rangle \neq 0$ or $\langle Sv,v \rangle \neq 0$, then either  
\[ \varphi_w(T)\varphi_w(S) \neq 0\mbox{ or }
 \varphi_v(T)\varphi_v(S) \neq 0. \]
On the other hand, if $\langle Tw,w \rangle = \langle Sv,v \rangle = 0 $, choose $n\ge 1$ such that $n^2\langle Tv,v\rangle + 2 \mathrm{Re}(\langle Tv,w\rangle)\neq 0$, $2\mathrm{Re}(\langle Sv,w\rangle)+\frac{1}{n^2}\langle Sw,w\rangle\neq 0$, and $nv+\frac{1}{n}w\neq 0$.
Define $q_0:=nv+\frac{1}{n}w$ and $q:=\left\|q_0\right\|^{-1}q_0$.  Then 
\[
\langle Tq_0,q_0 \rangle = n^2 \langle Tv,v \rangle + 2\operatorname{Re}(\langle Tv,w \rangle) + {\textstyle\frac{1}{n^2}}\langle Tw,w \rangle\neq 0
\]
and 
\[
\langle Sq_0,q_0 \rangle = n^2 \langle Sv,v \rangle + 2\operatorname{Re}(\langle Sv,w \rangle) + {\textstyle\frac{1}{n^2}}\langle Sw,w \rangle\neq 0.
\]
Hence, we have 
\[
\varphi_q(T)\varphi_q(S) \neq 0. 
\]
We conclude that in either case there exists a pure state $\varphi$ by \Cref{L:vector states are pure states} such that $\varphi(T)\varphi(S) \neq 0$, so $T$ and $S$ are not disjoint.
\end{proof}

The argument to show that there are no non-zero disjoint operators in $\mathrm{B}(H)_\mathrm{sa}$ where $H$ is a complex or real Hilbert space is analogous to the proof of \Cref{P: disjoint in B(H_q)}. 

\begin{proposition}\label{P:disjoint_B(H)}
Let $H$ be a real or complex Hilbert space. Then there are no non-zero disjoint operators in $\mathrm{B}(H)_\mathrm{sa}$.
\end{proposition}

This proposition is a generalisation of \cite[Proposition 16]{KalLemGaa2014}, where the space of symmetric $n\times n$-matrices with the cone of positive semi-definite matrices is considered.

\subsection{Disjointness in spin factors}
To show that there are no non-zero disjoint elements in a spin factor $H\oplus\mathbb{R}$, we will use the characterisation of the pure states given in \Cref{L:pure states in spin}.

\begin{proposition}\label{P:disjoint_in_spin}
There are no non-zero disjoint elements in a spin factor.
\end{proposition}

\begin{proof}
Let $(x,\lambda)$ and $(y,\mu)$ be non-zero elements of $H\oplus\mathbb{R}$. Then there are pure states $\varphi$ and $\psi$ such that $\varphi((x,\lambda))\neq 0$ and $\psi((y,\mu))\neq 0$. According to \Cref{L:pure states in spin}, there are unit vectors $v$ and $w$ of $H$ such that the state $\varphi$ is represented by $(v,1)$ and the state $\psi$ is represented by $(w,1)$. We have $\langle v,x \rangle + \lambda \neq 0$ and $\langle w,y \rangle + \mu \neq 0$. 

Let us first consider the case $v\neq w$ and $v\neq -w$. Then $v$ and $w$ are linearly independent.
For every $t\in (0,1)$, define $u_t= tv+(1-t)w$. Note that $u_t\neq 0$. Define $\eta_t=(u_t/\|u_t\|,1)$ for every $t\in(0,1)$. By \Cref{L:pure states in spin}, $\eta_t$ is a pure state. The equation $\eta_t((x,\lambda))=0$ holds for at most two values of $t\in(0,1)$. Indeed, $\eta_t((x,\lambda))=0$ comes down to $\langle u_t,x\rangle+\lambda\|u_t\|=0$, which yields 
\[\langle tv+(1-t)w,x\rangle = -\lambda\sqrt{t^2+(1-t)^2+2t(1-t)\langle v,w\rangle}.\]
Squaring both sides and sorting terms leads to the following quadratic equation in $t$,
\begin{align*}&\left( \langle v,x\rangle^2 +\langle w,x\rangle^2 -2\langle v,x\rangle \langle w,x\rangle -2\lambda^2+2\lambda^2\langle v,w\rangle \right)t^2\\
&\qquad +\left( -2\langle w,x\rangle^2 +2 \langle v,x\rangle\langle w,x\rangle +2\lambda^2-2\lambda^2\langle v,w\rangle \right)t +\langle w,x\rangle^2-\lambda^2=0.
\end{align*}
If the equation would be satisfied for three or more values of $t$, then its coefficients are zero, so $\lambda^2=\langle w,x\rangle^2$, $2\langle v,x\rangle\langle w,x\rangle - 2\langle w,x\rangle^2 \langle v,w\rangle =0$, and
\begin{equation}\label{eq.coefficientquadraticequation}
\langle v,x\rangle^2 - \langle w,x\rangle^2 -2\langle v,x\rangle \langle w,x\rangle + 2\langle w,x\rangle^2\langle v,w\rangle =0.
\end{equation}
Addition of the latter two equations yields $\langle v,x\rangle^2=\langle w,x\rangle^2$ and from \eqref{eq.coefficientquadraticequation} we then also obtain $$\langle v,x\rangle \langle w,x\rangle = \langle w,x\rangle^2\langle v,w\rangle.$$ Observe that $\langle w,x\rangle\neq 0$, as otherwise $\langle v,x\rangle=0$ and $\lambda=0$, whereas $\langle v,x\rangle+\lambda\neq 0$. Therefore, from $\langle v,x\rangle^2=\langle w,x\rangle^2$ we obtain either $\langle v,x\rangle=\langle w,x\rangle$ and $\langle v,w\rangle=1$, or $\langle v,x\rangle=-\langle w,x\rangle$ and $\langle v,w\rangle=-1$. As $v$ and $w$ are unit vectors, it follows from Cauchy-
Schwarz that $v=w$ or $v=-w$, which both yield a contradiction. Thus, the equation $\eta_t((x,\lambda))=0$ holds for at most two values of $t\in(0,1)$. Similarly, $\eta_t((y,\mu))=0$ holds for at most two values of $t\in(0,1)$. We conclude that there exists $t\in (0,1)$ with $\eta_t((x,\lambda))\neq 0$ and $\eta_t((y,\mu))\neq 0$. Therefore, $(x,\lambda)$ and $(y,\mu)$ are not disjoint. 

Let us now address the case $v=w$ or $v=-w$. Since $H$ is at least two dimensional, there exists a unit vector $z\in H$ which is linearly independent of $v$ and $w$. For every $s\in (0,1)$, define $w_s=w+sz$ and note that $w_s\neq 0$. Then $\psi_s=(w_s/\|w_s\|,1)$ is a pure state by \Cref{L:pure states in spin}. It follows by continuity that there exists $s\in (0,1)$ such that $\psi_s((y,\mu))=\langle w_s,y\rangle/\|w_s\|+\mu \neq 0$. The first part of the proof with $\psi$ replaced by $\psi_s$ yields that $(x,\lambda)$ and $(y,\mu)$ are not disjoint. 
\end{proof}

A finite-dimensional space with a Lorentz cone is a special case of \Cref{P:disjoint_in_spin} and for such a space the result also follows from \cite[Theorem 15]{KalLemGaa2014}.

\subsection{Disjointness in $\mathrm{M}_3(\mathbb{O})_\mathrm{sa}$}
We will use the characterisation of the pure states in \Cref{L:pure_states_in_Albert} to show that there are no non-zero disjoint elements in $\mathrm{M}_3(\mathbb{O})_\mathrm{sa}$. 

\begin{proposition}\label{P:disjoint_M_3}
There are no non-zero disjoint elements in $\mathrm{M}_3(\mathbb{O})_\mathrm{sa}$.
\end{proposition}

\begin{proof}
Let $A,B \in \mathrm{M}_3(\mathbb{O})_\mathrm{sa}$ be non-zero and distinct. By \Cref{L:pure_states_in_Albert}, there are minimal projections $P$ and $Q$ such that $\langle A,P \rangle \neq 0$ and $\langle B,Q \rangle \neq 0$. If either $\langle A,Q \rangle \neq 0$ or $\langle B,P \rangle \neq 0$, then $A$ and $B$ are not disjoint. If $\langle A,Q \rangle = \langle B,P \rangle =0$, we distinguish two cases.

As case 1, consider $P\circ Q\neq 0$. Then $\langle P,Q \rangle \neq 0$ by \cite[Exercise~III.3]{Faraut-Koranyi}. Note that $\langle P,Q \rangle \le \sqrt{\langle P,P \rangle}\sqrt{\langle Q,Q \rangle}=1$ by the Cauchy-Schwarz inequality, and $\langle P,Q \rangle <1$ as otherwise 
\[
\langle P-Q,P-Q \rangle = \langle P,P\rangle -2\langle P,Q\rangle +\langle Q,Q \rangle =0,
\]
which is impossible as $P$ and $Q$ are distinct. Hence, there exists $\theta\in \left(0,\frac{1}{2}\pi\right)$ such that $\cos^2 \theta:=\langle P,Q\rangle$. Define 
\[
P':= \begin{pmatrix}
1&0\\0&0
\end{pmatrix}, \qquad Q':=\begin{pmatrix}
\cos^2\theta & \cos \theta \sin \theta\\\cos \theta \sin \theta&\sin^2\theta
\end{pmatrix}, \qquad R':= \begin{pmatrix}
\cos^2\theta & \frac{1}{2}\cos \theta \sin \theta\\\frac{1}{2}\cos \theta \sin \theta&0
\end{pmatrix}.
\]
The Jordan algebra generated by $P$ and $Q$ (without $I_3$) is isomorphic to $\mathrm{M}_2(\mathbb{R})_\mathrm{sa}$ via the map 
\[
\alpha P+\beta Q + \gamma P\circ Q \mapsto \alpha P'+\beta Q'+\gamma R' 
\]
by \cite[Proposition~IV.1.6]{Faraut-Koranyi}. By \Cref{R:spin}, all the minimal projections in $\mathrm{M}_2(\mathbb{R})_\mathrm{sa}$ are of the form 
\[
\begin{pmatrix}
\frac{1}{2}+x_2&x_1\\x_1&\frac{1}{2}-x_2
\end{pmatrix},
\]
where $x_1^2+x_2^2=\frac{1}{4}$. Let $n\in\mathbb{N}$. With $x_1:=\sqrt{\frac{1}{n}-\frac{1}{n^2}}$ and $x_2:=\frac{1}{n}-\frac{1}{2}$, it follows that 
\[
S_n':= \begin{pmatrix}
\frac{1}{n}&\sqrt{\frac{1}{n}-\frac{1}{n^2}}\\\sqrt{\frac{1}{n}-\frac{1}{n^2}}&1-\frac{1}{n}
\end{pmatrix}
\]
is a minimal projection in $\mathrm{M}_2(\mathbb{R})_\mathrm{sa}$. A computation shows that 
\[S_n'=\alpha_n P'+\beta_n Q'+\gamma_n R',\]
where
\begin{align*}
    \alpha_n&:=\frac{1}{n} +\frac{(n-1)\cos^2\theta}{n\sin^2\theta}-\frac{2\sqrt{n-1}\cos\theta}{n\sin\theta},\\\beta_n&:=\frac{n-1}{n\sin^2\theta},\\\gamma_n&:=\frac{2\sqrt{n-1}}{n\cos\theta\sin\theta}-\frac{2(n-1)}{n\sin^2\theta}.
\end{align*}
Since $S'_n$ is a projection in $\mathrm{M}_2(\mathbb{R})_\mathrm{sa}$ with $\mathrm{trace}(S'_n)=1$, the preimage 
\[
S_n:=\alpha_n P+\beta_n Q+\gamma_nP\circ Q
\]
is a projection in $\mathrm{M}_3(\mathbb{O})_\mathrm{sa}$ with $\mathrm{trace}(S_n)=1$ as well. The spectral decomposition of $S_n$ expands $S_n$ as a linear combination of minimal projections and the spectrum of $S_n$ consists of the coefficients in this linear combination. As $S_n$ is a projection, we have $\sigma(S_n)=\{0,1\}$. Since the trace of $S_n$ equals $1$, only one term in its spectral decomposition can be non-zero. Thus, $S_n$ is a minimal projection. Suppose there is a subsequence of $(\langle A,S_n \rangle)_{n\ge 1}$ such that $\langle A,S_{n_k} \rangle = 0$ for all $k\ge 1$. Note that as $n\to\infty$,
\[\alpha_n\to \frac{\cos^2\theta}{\sin^2\theta},\quad \beta_n\to \frac{1}{\sin^2\theta},\quad\mbox{and }\gamma_n\to -\frac{2}{\sin^2\theta}.\]
By letting $k\to\infty$ in $\langle A,S_{n_k} \rangle = 0$, we find
\[\left\langle A, \frac{\cos^2\theta}{\sin^2\theta}P +\frac{1}{\sin^2\theta}Q-\frac{2}{\sin^2\theta} P\circ Q\right\rangle =0.\]
As $\langle A,Q\rangle=0$, we obtain $\langle A,P\circ Q\rangle = \frac{1}{2}\langle A, \cos^2\theta P\rangle$. Hence  $0=\langle A,S_{n_k} \rangle = (\alpha_{n_k}+\frac{1}{2}\cos^2\theta \gamma_{n_k})\langle A,P \rangle$ for all $k\ge 1$. But 
\[
\alpha_{n_k}+\frac{1}{2}\cos^2\theta \gamma_{n_k} = \frac{1}{n_{k}}-\frac{\sqrt{n_k-1}\cos\theta}{n_k\sin\theta}<0
\]
for all $n_k>\frac{1}{\cos^2\theta}$, which yields a contradiction. Hence, there is a number $N\ge 1$ such that $\langle A,S_n \rangle \neq 0$ for all $n\ge N$. 

There is a number $M\ge 1$ such that $\langle B,S_n \rangle \neq 0 $ whenever $n\ge M$. Indeed, suppose there is a subsequence such that $\langle B, S_{n_k}\rangle=0$ for all $k\ge 1$. Taking limits yields
\[\left\langle B, \frac{\cos^2\theta}{\sin^2\theta}P +\frac{1}{\cos^2\theta}Q -\frac{2}{\sin^2\theta}P\circ Q\right\rangle=0.\]
As $\langle B,P\rangle=0$, we obtain
\[\langle B,P\circ Q\rangle =\frac{\sin^2\theta}{2\cos^2\theta}\langle B,Q\rangle.\]
Hence
\[0=\langle B,S_{n_k}\rangle = \left(\beta_{n_k}+\frac{\sin^2\theta}{2\cos^2\theta}\gamma_{n_k}\right)\langle B,Q\rangle.\]
Note that
\[\beta_n+\frac{\sin^2\theta}{2\cos^2\theta}\gamma_n=\frac{n-1}{n}\left(\frac{1}{\sin^2\theta}-\frac{1}{\cos^2\theta}\right) +\frac{\sqrt{n-1}\sin\theta}{n\cos^3\theta},\]
which is non-zero for large $n$. Thus, we arrive at a contradiction. 

Therefore, for any $n\ge N,M$, we have $\langle A,S_n\rangle\langle B,S_n\rangle \neq 0$, which yields that $A$ and $B$ are not disjoint, as $\langle \cdot,S_n\rangle$ is a pure state by \Cref{L:pure_states_in_Albert}. 

Next, we consider the case 2, namely $P\circ Q =0$. We will construct a minimal projection $R$ such that $R\circ P\neq 0$ and $\langle B,R\rangle \neq 0$. Then case 1 of the proof with $Q$ replaced by $R$ yields that $A$ and $B$ are not disjoint. To construct $R$, we use that $\mathrm{M}_3(\mathbb{O})_\mathrm{sa}$ is a factor, so there is a $C\in \mathrm{M}_3(\mathbb{O})_\mathrm{sa}$ such that $P\circ C=Q\circ C=\frac{1}{2}C$ and $C^2=P+Q$ by \cite[Proposition~IV.1.4(i)]{Faraut-Koranyi} and \cite[Proposition~IV.2.4(i)]{Faraut-Koranyi}. Note that $\langle C,P \rangle = \langle C,P^2 \rangle = \langle P\circ C,P \rangle = \frac{1}{2}\langle C,P \rangle$, so $\langle C,P \rangle = 0$. Similarly, it follows that $\langle C,Q \rangle = 0$, and hence,
\[
\mathrm{trace}(C)=\mathrm{trace}\left(\textstyle\frac{1}{2}C+\frac{1}{2}C\right)=\mathrm{trace}(P\circ C)+\mathrm{trace}(Q\circ C)= \langle C,P \rangle+ \langle C,Q \rangle = 0.
\]
For numbers $\alpha,\beta$ such that $\alpha^2+\beta^2=1$, a straightforward calculation shows that 
\[
\alpha^2 P+\beta^2 Q+\alpha\beta C
\]
is a projection.  Its trace equals 1, so it is a minimal projection. Define for $n\ge 1$ the minimal projection $R_n$ by
\[
R_n:=\frac{1}{n^2}P+\left(1-\frac{1}{n^2}\right)Q+\frac{1}{n}\sqrt{1-\frac{1}{n^2}}C.
\]
Since $\langle B,P\rangle=0$ and $\langle B,Q\rangle \neq 0$, there exists $N$ such that for $n\ge N$ we have $\langle B,R_n\rangle \neq 0$. Moreover, there exists $n\ge N$ such that $R_n\circ P\neq 0$, so $R:=R_n$ is as required.
\end{proof}

\section{Order direct sums of order unit spaces}\label{sec.orderdirectsums}

By \Cref{the.atomicdecomp}, every atomic JBW-algebra is a direct sum of atomic JBW-algebra factors. In Section \ref{sec.antilattices}, we have shown that these factors are disjointness free anti-lattices. Thus, every atomic JBW-algebra is the order direct sum of order unit spaces that are disjointness free anti-lattices. In this section, we consider order direct sums of order unit spaces, we investigate their vector lattice covers, relate Riesz* homomorphisms on the direct sum with those on the components, and we determine the functional representation. For clarity of the ideas, we first consider the case of two components and then formulate the general case,  which is entirely similar.

\begin{lemma}\label{lem:1}
	Let $Y_1$ and $Y_2$ be partially ordered vector spaces and $X_1\subseteq Y_1$, $X_2\subseteq Y_2$ subspaces. If $X_i$ is order dense in $Y_i$ ($i\in\{1,2\}$), then $X_1\times X_2$ is order dense in $Y_1\times Y_2$.
\end{lemma}
\begin{proof}
	Let $y=(y_1,y_2)\in Y_1\times Y_2$ and $S:=\{(x_1,x_2)\colon\, \ x_i\in X_i,\  x_i\geq y_i,\  i\in\{1,2\}\}$. Clearly, $y$ is a lower bound of $S$. Let $z=(z_1,z_2)\in Y_1\times Y_2$ be a lower bound of $S$. For $i\in\{1,2\}$, we have $z_i\leq x_i$ for every $x_i\in X_i$ with $x_i\geq y_i$. As $X_i$ is order dense in $Y_i$, we obtain $z_i\leq y_i$. Hence $z\leq y$. Therefore $y=\inf S$.
\end{proof}
We note an immediate consequence of \Cref{lem:1}.
\begin{proposition}\label{pro:VLC}
	If $X_1$, $X_2$ are pre-Riesz spaces and $(Y_1,i_1)$, $(Y_2,i_2)$ vector lattice covers of $X_1$, $X_2$, respectively, then $(Y_1\times Y_2,i_1\times i_2)$ is a vector lattice cover of $X_1\times X_2$.
\end{proposition} 

We relate Riesz* homomorphisms on pre-Riesz spaces to Riesz* homomorphisms on their Cartesian product.
\begin{proposition}\label{pro:Riesz_star_CP}
	Let $X_1$, $X_2$ be pre-Riesz spaces and $Z$ an Archimedean pre-Riesz space.
	\begin{itemize}
		\item[(i)] If $h_1\colon X_1\to Z$ is a Riesz* homomorphism, then 
		\[h\colon X_1\times X_2\to Z,\ (x_1,x_2)\mapsto h_1(x_1),\]
		is a Riesz* homomorphism.
		\item[(ii)] Assume that $Z$ is a disjointness free anti-lattice. If $h\colon X_1\times X_2\to Z$ is a Riesz* homomorphism, then
		\begin{align*}
		&h_1\colon X_1\to Z, \ x\mapsto h(x,0),\quad
		 h_2\colon X_2\to Z, \ x\mapsto h(0,x),
		\end{align*}
		are Riesz* homomorphisms, and there is $k\in\{1,2\}$ such that for all $(x_1,x_2)\in X_1\times X_2$ we have $h(x_1,x_2)=h_k(x_k)$.
	\end{itemize}
\end{proposition}
\begin{proof}
For $j\in \{1,2\}$, let $(Y_j,i_j)$ be a vector lattice cover of $X_j$.
Let $(Z^\delta,i_Z)$ be the Dedekind completion of $Z$, see \cite[Theorem~2.1.13]{KalvanGaa2019}.	
	
$(i)$ Take for the Riesz completion $X_1^\rho$ of $X_1$ the  Riesz subspace of $Y_1$ generated by $i_1[X_1]$.  By \Cref{the:vanHaandel}, there exists a Riesz homomorphism  $h_1^\rho\colon X_1^\rho\to Z^\delta$ such that $h_1^\rho\circ i_1=i_Z\circ h_1$. As $i_1[X_1]$ is majorizing in $Y_1$, we have that $X_1^\rho$ is majorizing in $Y_1$. By \Cref{thm:LLS}, there is a Riesz homomorphism $\hat{h}_1\colon Y_1\to Z^\delta$ extending $h_1^\rho$.
Define
\[\hat{h}\colon Y_1\times Y_2\to Z^\delta,\quad (y_1,y_2)\mapsto \hat{h}_1(y_1).\]
Clearly, $\hat{h}$ is linear. Moreover, for every $(y_1,y_2)\in Y_1\times Y_2$, we have
\[\hat{h}\left(\left|(y_1,y_2)\right|\right)=\hat{h}\left(\left(|y_1|,|y_2|\right)\right)=\hat{h}_1\left(|y_1|\right)=\big|\hat{h}_1\left(y_1\right)\big|=\big|\hat{h}\left(y_1,y_2\right)\big|,\]
hence $\hat{h}$ is a Riesz homomorphism. For every $(x_1,x_2)\in X_1\times X_2$, we have
\begin{align*}
i_Z(h(x_1,x_2))&=i_Z(h_1(x_1))=h_1^\rho(i_1(x_1))=\hat{h}_1(i_1(x_1))=\hat{h}\left(i_1(x_1),i_2(x_2)\right)\\&=\hat{h}\left((i_1\times i_2)(x_1,x_2)\right).
\end{align*}
A combination of \Cref{pro:VLC} and \Cref{the:vanHaandel} yields that $h$ is a Riesz* homomorphism.

$(ii)$ Let $Y^\rho$ be the Riesz subspace of $Y_1\times Y_2$ generated by $(i_1\times i_2)[X_1\times X_2]$. According to \Cref{the:vanHaandel}, there is a Riesz homomorphism $\hat{h}\colon Y^\rho\to Z^\delta$ such that $\hat{h}\circ (i_1\times i_2)=i_Z\circ h$. As $Y^\rho$ is majorizing in $Y_1\times Y_2$, by \Cref{thm:LLS}, there is a Riesz homomorphism $H\colon Y_1\times Y_2\to Z^\delta$ such that $H$ extends $\hat{h}$. In particular, for every $y\in (i_1\times i_2)[X_1\times X_2]$, we have $H(y)=\hat{h}(y)$, therefore $H\circ(i_1\times i_2)=i_Z\circ h$. For every $(y_1,y_2)\in Y_1\times Y_2$, define
\[H_1\colon Y_1\to Z^\delta, \ y_1\mapsto H(y_1,0), \ \mbox{ and } \ H_2\colon Y_2\to Z^\delta, \ y_2\mapsto H(0,y_2).\]
$H_1$ is linear, and for every $y_1\in Y_1$ we have 
\[H_1(|y_1|)=H((|y_1|,0))=H(|(y_1,0)|)=|H(y_1,0)|=|H_1(y_1)|,\]
hence $H_1$ is a Riesz homomorphism. Moreover, for every $x\in X_1$, we have 
\begin{align*}
(H\circ i_1)(x)&=H_1\left(i_1(x)\right)=H\left((i_1(x),0)\right)=H\left((i_1(x),i_2(0))\right)\\
&=H\left((i_1\times i_2)(x,0)\right)=i_Z(h(x,0))=i_Z(h_1(x)),
\end{align*}
hence $H_1\circ i_1=i_Z\circ h_1$.
Similarly, $H_2$ is a Riesz homomorphism with $H_2\circ i_2=i_Z\circ h_2$.
It follows by \Cref{the:vanHaandel} that $h_1$ and $h_2$ are Riesz* homomorphisms.
It remains to show that there exists $k\in \{1,2\}$ such that $h(x_1,x_2)=h_k(x_k)$ for every $(x_1,x_2)\in X_1\times X_2$. 
First, observe that for every $(x_1,x_2)\in X_1\times X_2$ we have 
\[h(x_1,x_2)=h(x_1,0)+h(0,x_2)=h_1(x_1)+h_2(x_2).\]
If $h_1=0$ then we choose $k=2$. Otherwise, if $h_1\neq 0$, we show that $h_2=0$. Indeed, choose $x_1\in X_1$ such that $h_1(x_1)\neq 0$. For every $x_2\in X_2$, we have $\left(i_1(x_1),0\right)\perp\left(0,i_2(x_2)\right)$ in $Y_1\times Y_2$,
therefore $H\left(i_1(x_1),0\right)\perp H\left(0,i_2(x_2)\right)$ in $Z^\delta$. Since \[H\left(i_1(x_1),0\right)=H_1\left(i_1(x_1)\right)=i_Z(h_1(x_1))\neq 0,\]
and $Z$ is a disjointness free anti-lattice, it follows that $H\left(0,i_2(x_2)\right)=0$, hence
\[i_Z(h_2(x_2))=H_2(i_2(x_2))=H(0,i_2(x_2))=0,\]
which implies $h_2(x_2)=0$. Thus, $h_2=0$.
\end{proof}
For an order unit space $X$, we have,  by \Cref{pro:RieszhomOUS},  
\[\overline{\Lambda}_X=\left\{h\colon X\to\mathbb{R}\colon\, \ h \mbox{ is a Riesz* homomorphism}, h(u)=1\right\}.\] 
Therefore, we obtain the following consequence of \Cref{pro:Riesz_star_CP}.
\begin{proposition}\label{pro.sumvlcover}
	Let $(X_1,K_1, u_1)$ and $(X_2,K_2,u_2)$ be order unit spaces. Then $(X_1\times X_2, K_1\times K_2, (u_1,u_2))$ is an order unit space and \[\overline{\Lambda}_{X_1\times X_2} =\left\{(x_1,x_2)\mapsto f(x_1)\colon\, \ f\in \overline{\Lambda}_{X_1}\right\}\cup \left\{(x_1,x_2)\mapsto g(x_2)\colon\, \ g\in  \overline{\Lambda}_{X_2}\right\}.\]
	Moreover, if we consider the disjoint union topology on $\overline{\Lambda}_{X_1}\cup \overline{\Lambda}_{X_2}$, then the map
	$J\colon \overline{\Lambda}_{X_1}\cup \overline{\Lambda}_{X_2}\to \overline{\Lambda}_{X_1\times X_2}$ given by 
	\[(Jf)(x_1,x_2):=\begin{cases}
	f(x_1)      & \ \text{if } f\in \overline{\Lambda}_{X_1},\\
f(x_2)      & \ \text{if } f\in\overline{\Lambda}_{X_2},	
	\end{cases}\]
	for $(x_1,x_2)\in X_1\times X_2$, is a homeomorphism.
\end{proposition}

\begin{proof}
	We apply \Cref{pro:Riesz_star_CP} for $Z:=\mathbb{R}$.
	
	Let $h\in \overline{\Lambda}_{X_1\times X_2}$. Let $h_1$, $h_2$, and $k$ be as in \Cref{pro:Riesz_star_CP}$(ii)$. Then $h_k(u_k)=h(u_1,u_2)=1$, hence $h_k\in \overline{\Lambda}_{X_k}$.
	
	Conversely, let $h_1\in \overline{\Lambda}_{X_1}$ and define $h$ as in \Cref{pro:Riesz_star_CP}$(i)$. Then $h(u_1,u_2)=h_1(u_1)=1$, hence
	$h\in \overline{\Lambda}_{X_1\times X_2}$. For $h_2\in \overline{\Lambda}_{X_2}$, the proof is analogous.	
	
	It remains to show that $J$ and $J^{-1}$ are continuous. Let $(f_i)_i$ be a net in $\overline{\Lambda}_{X_1}\cup \overline{\Lambda}_{X_2}$ and let $f\in\overline{\Lambda}_{X_1}\cup \overline{\Lambda}_{X_2}$ be such that $f_i\to f$. Assume first that $f\in \overline{\Lambda}_{X_1}$. Then there is $i_0$ with $f_i\in \overline{\Lambda}_{X_1}$ for every $i\geq i_0$. Hence, for every $(x_1,x_2)\in X_1\times X_2$ and $i\geq i_0$, we have 
	\[(Jf_i)(x_1,x_2)=f_i(x_1)\to f(x_1)=(Jf)(x_1,x_2),\]
	therefore $Jf_i\to Jf$. Similarly, if $f\in \overline{\Lambda}_{X_2}$, we  obtain $Jf_i\to Jf$. Thus, $J$ is continuous. Since $\overline{\Lambda}_{X_1}\cup \overline{\Lambda}_{X_2}$ is compact and $\overline{\Lambda}_{X_1\times X_2}$ is Hausdorff, it follows that $J^{-1}$ is continuous. 
\end{proof}

\begin{corollary}
	The functional representation $\left(\mathrm{C}\left(\overline{\Lambda}_{X_1\times X_2}\right), \Phi_{X_1\times X_2}\right)$ of $X_1\times X_2$ satisfies 
	\begin{eqnarray}\nonumber
	\mathrm{C}\left(\overline{\Lambda}_{X_1\times X_2}\right)&=&\mathrm{C}\left(\overline{\Lambda}_{X_1}\right)\oplus \mathrm{C}\left(\overline{\Lambda}_{X_2}\right),\\\nonumber
	\Phi_{X_1\times X_2}(x_1,x_2)&=&
		\left(\Phi_{X_1}(x_1),
			\Phi_{X_2}(x_2)\right)
	\end{eqnarray}
	for all $x_1\in X_1$ and $x_2\in X_2$.
\end{corollary}

\begin{proposition}\label{pro.sumprojectionbands}
	Let $X_1$ and $X_2$ be pre-Riesz spaces. Then $X_1\times\{0\}$ and $\{0\}\times X_2$ are projection bands in $X_1\times X_2$.
\end{proposition}

\begin{proof}
		Let $(Y_1,i_1)$, $(Y_2,i_2)$ be  vector lattice covers of $X_1$, $X_2$, respectively. 
By \Cref{pro:VLC},
	 $(Y_1\times Y_2,i_1\times i_2)$ is a vector lattice cover of $X_1\times X_2$. 
	 
	 Let $x_1\in X_1$ and $x_2\in X_2$. Then $\left(i_1(x_1),0\right)\perp\left(0,i_2(x_2)\right)$ in $Y_1\times Y_2$, hence $(x_1,0)\perp (0,x_2)$ in $X_1\times X_2$. Thus, $X_1\times\{0\}\subseteq (\{0\}\times X_2)^\mathrm{d}$ and 
	 $\{0\}\times X_2\subseteq (X_1\times\{0\})^\mathrm{d}$.
	 
	 Let $(v_1,v_2)\in  (X_1\times\{0\})^\mathrm{d}$. Then, for every $x_1\in X_1$, we have \[\left(i_1(v_1),i_2(v_2)\right)\perp(i_1(x_1),0).\] By order denseness of $i_1[X_1]$ in $Y_1$, we obtain  $\left(i_1(v_1),i_2(v_2)\right)\perp(y_1,0)$ for every $y_1\in Y_1$. Therefore, $i_1(v_1)=0$, hence $v_1=0$. We get $(v_1,v_2)=(0,v_2)\in \{0\}\times X_2$. Consequently, $(X_1\times\{0\})^\mathrm{d}\subseteq \{0\}\times X_2$.
	 
	 We conclude $(X_1\times\{0\})^\mathrm{d}= \{0\}\times X_2$, and, similarly, 
	 $ (\{0\}\times X_2)^\mathrm{d}= X_1\times\{0\}$. Thus, $(X_1\times\{0\})^\mathrm{dd}=X_1\times\{0\}$, which means that $X_1\times\{0\}$ is a band. Analogously, $\{0\}\times X_2$ is a band. Finally, $X_1\times X_2=(X_1\times\{0\})\oplus (\{0\}\times X_2)$, hence 
	 $X_1\times\{0\}$ and  $\{0\}\times X_2$ are projection bands.
\end{proof}

Analogues of the statements of \Cref{pro.sumvlcover}, \Cref{pro:Riesz_star_CP}, and \Cref{pro.sumprojectionbands} are valid for arbitrary direct sums of order unit spaces, as we state next without proof.

\begin{proposition}
Let $((V_i,C_i,u_i))_{i\in \mathcal{I}}$ be a family of order unit spaces and, for every $i\in \mathcal{I}$, let $Y_i$ be an Archimedean partially ordered vector space and $j_i\colon V_i\to Y_i$ be a bipositive linear map such that $j_i(u_i)$ is an order unit in $Y_i$. Denote $V=\bigoplus_{i\in \mathcal{I}} V_i$.
\begin{itemize}
\item[(i)] If $j_i[V_i]$ is order dense in $Y_i$ for every $i\in \mathcal{I}$, then $\bigoplus_{i\in \mathcal{I}} j_i[V_i]$ is order dense in $\bigoplus_{i\in \mathcal{I}} Y_i$. Consequently, if $(Y_i,j_i)$ is a vector lattice cover of $V_i$ for every $i\in \mathcal{I}$, then $(\bigoplus_{i\in \mathcal{I}} Y_i,\bigoplus_{i\in \mathcal{I}} j_i)$ is a vector lattice cover of $V$. 
\item[(ii)] Let $Z$ be a disjointness free anti-lattice. For every $k\in \mathcal{I}$ and every Riesz* homomorphism  $h\colon V_k\to Z$ the map $H\colon V\to Z$ defined by
\[H(x):=h(x(k)),\quad x\in \bigoplus_{i\in \mathcal{I}} V_i,\]
is a Riesz* homomorphism. Conversely, if $H\colon V\to Z$ is a Riesz* homomorphism, then for every $k\in \mathcal{I}$ the map $h\colon V_k\to Z$ defined by
\[h_k(v):=H(\Phi_k(v)),\quad v\in V_k,\]
is a Riesz* homomorphism. Moroever, in the latter case there exists $k\in \mathcal{I}$ such that $H(v)=h_k(v)$ for every $v\in V$.
\item[(iii)] We have
\[\overline{\Lambda}_V=\bigcup_{k\in \mathcal{I}} \left\{x\mapsto f(x(k))\colon V\to \mathbb{R}\colon\, f\in \overline{\Lambda}_{V_k} \right\}.\]
Moreover, consider the disjoint union topology on $\bigcup_{k\in \mathcal{I}} \overline{\Lambda}_{V_k}$ and  the map 
$J\colon \bigcup_{i\in \mathcal{I}} \overline{\Lambda}_{V_i}\to \overline{\Lambda}_V$ given by $(J(f))(x)=f(x(i))$, where $i\in \mathcal{I}$ is such that $f\in \overline{\Lambda}_{V_i}$. Then $J$ is a homeomorphism.
\end{itemize}
\end{proposition}

%subspaces of $V$ such that for every $v\in V$ there are unique $v_\alpha\in V_\alpha$ $(\alpha\in I)$ with $v_\alpha=0$ for all but finitely many $\alpha$ and $v=\sum_{\alpha\in I} v_\alpha$, then $V$ is the \emph{direct sum} of $(V_\alpha)_{\alpha\in I}$. If for every $v\in V$ we have $v\ge 0$ if and only if $v_\alpha\ge 0$ for every $\alpha\in I$, where $v$ and $v_\alpha$ are as above, then we say that $V$ is the \emph{order direct sum} of $(V_\alpha)_{\alpha\in I}$ and we denote 
%\begin{equation}\label{eq.Visorderdirectsum}
%V=\bigoplus_{\alpha\in I} V_\alpha.
%\end{equation}

\section{Disjointness and bands in order direct sums of disjointness free anti-lattices}\label{sec.disjointnessinsums}

%Recall that every band is an ideal, so that for every $\alpha\in I$ the subspace $V_\alpha$ in the order direct sum \eqref{eq.Visorderdirectsum} is a directed ideal of $V$. By \cite[Theorem 4.3.22]{Kalauch-vanGaans-book}, this means that $V_\alpha$ is directed and for every $v\in V$ and $a,b\in V_\alpha$ with $a\le v\le b$ we have $v\in V_\alpha$.

In this section, we characterise disjointness and bands in atomic JBW-algebras. First, we show that two elements in an order direct sum of order unit spaces are disjoint if and only if they are componentwise disjoint. 

\begin{lemma}\label{lem.entrywisedisjointness}
Let $((V_i,C_i,u_i))_{i\in \mathcal{I}}$ be a collection of order unit spaces with order direct sum $(V,C,u)$. Let $v,w\in V$. Then $v$ and $w$ are disjoint in $V$ if and only if for every $i\in \mathcal{I}$ the elements $v(i)$ and $w(i)$ are disjoint in $V_i$. 
\end{lemma}
\begin{proof}
Assume that $v$ and $w$ are disjoint in $V$.  We have 
\[\{v(i)+w(i),-v(i)-w(i)\}^\mathrm{u}=\{v(i)-w(i),-v(i)+w(i)\}^\mathrm{u}.\]
Indeed, let $z\in \{v(i)+w(i),-v(i)-w(i)\}^\mathrm{u}$. Take $t:=\sup_{j\in \mathcal{I}} \|v(j)+w(j)\|_{u_j}$ and define $x\in V$ by $x(i)=z$ and $x(j)=tu_j$ for every $j\in \mathcal{I}\setminus\{i\}$. Then $x\in V$ and $x\ge v+w,-v-w$. As $v$ and $w$ are disjoint, we obtain $x\ge v-w,-v+w$. In particular, $z=x(i)\in \{v(i)-w(i),-v(i)+w(i)\}^\mathrm{u}$. The converse inclusion is proven similarly. Thus, $v(i)$ and $w(i)$ are disjoint.

Next, assume that for every $i\in \mathcal{I}$ we have that $v(i)$ and $w(i)$ are disjoint. For every $z\in\{v+w,-v-w\}^\mathrm{u}$ we have $z(i)\in\{v(i)+w(i),-v(i)-w(i)\}^\mathrm{u}$, so $z\in\{v-w,-v+w\}^\mathrm{u}$. Hence, $\{v+w,-v-w\}^\mathrm{u}\subseteq \{v-w,-v+w\}^\mathrm{u}$. The converse inclusion follows similarly, so $v$ and $w$ are disjoint.
\end{proof}

Next we observe that the components of an order direct sum of order unit spaces are projection bands that are pairwise disjoint. Let $((V_i,C_i,u_i))_{i\in \mathcal{I}}$ be a collection of order unit spaces with order direct sum $(V,C,u)$. For every $\mathcal{J}\subseteq \mathcal{I}$, we define $P_\mathcal{J}\colon V\to V$ by $P_\mathcal{J}(v)=w$, where $w_i=v_i$ for every $i\in \mathcal{J}$ and $w_i=0$ otherwise. Clearly, $P_\mathcal{J}$ is a projection with range $\Phi_\mathcal{J}\left[\bigoplus_{i\in \mathcal{J}} V_i\right]$, where $\Phi_\mathcal{J}$ is as defined below \eqref{eq.Visorderdirectsum}. 

\begin{proposition}\label{pro.componentsarebands}
Let $((V_i,C_i,u_i))_{i\in \mathcal{I}}$ be a collection of order unit spaces with order direct sum $(V,C,u)$ and let $\mathcal{J}\subseteq \mathcal{I}$. 
\begin{itemize}
\item[(i)] $P_\mathcal{J}$ is a band projection. 
\item[(ii)] $\Phi_\mathcal{J}\left[\bigoplus_{i\in \mathcal{J}} V_i\right]$ is a projection band in $V$ and is directed.
\item[(iii)] For every $x\in \Phi_\mathcal{J}\left[\bigoplus_{i\in \mathcal{J}} V_i\right]$ and $y\in \Phi_{\mathcal{I}\setminus \mathcal{J}}\left[\bigoplus_{i\in \mathcal{I}\setminus \mathcal{J}} V_i\right]$ we have that $x$ and $y$ are disjoint.
\item[(iv)] We have $\Phi_\mathcal{J}\left(u|_\mathcal{J}\right)^\mathrm{dd}=\Phi_\mathcal{J}\left[\bigoplus_{i\in \mathcal{J}} V_i\right]$.
\end{itemize}
\end{proposition}
\begin{proof}
$(i)$ As $P_\mathcal{J}$ and $I-P_\mathcal{J}=P_{\mathcal{I}\setminus \mathcal{J}}$ both are positive, we have that $P_\mathcal{J}$ is an order projection, hence $P_\mathcal{J}$ is a band projection by \cite[Proposition 3.1]{Glueck2021}.

$(ii)$ Follows directly from $(i)$ and the fact that $\Phi_\mathcal{J}\left[\bigoplus_{i\in \mathcal{J}} V_i\right]$ has an order unit.

$(iii)$ For every $i\in \mathcal{J}$ we have $y(i)=0$ and for every $i\in \mathcal{I}\setminus \mathcal{J}$ we have $x(i)=0$. It follows from \Cref{lem.entrywisedisjointness} that $x$ and $y$ are disjoint.

$(iv)$ As $u|_\mathcal{J}\in\bigoplus_{i\in \mathcal{J}} V_i$ and $\Phi_\mathcal{J}\left[ \bigoplus_{i\in \mathcal{J}} V_i\right]$ is a band in $V$, we have $\Phi_\mathcal{J}\left(u|_\mathcal{J}\right)^\mathrm{dd}\subseteq\Phi_\mathcal{J}\left[\bigoplus_{i\in \mathcal{J}} V_i\right]$.

In the proof of the converse inclusion, we call a set $S$ of a partially ordered vector space $X$ \emph{full} if for every $s,t\in S$ and $x\in X$ with $s\le x\le t$ we have $x\in S$. As $\Phi_\mathcal{J}\left(u|_\mathcal{J}\right)^\mathrm{dd}$ is a band, it is a full subspace by \cite[Theorem 4.3.13 and Lemma 4.3.5]{KalvanGaa2019}. 
%Moreover, $\Phi_\mathcal{J}\left(u|_\mathcal{J}\right)\in \Phi_\mathcal{J}\left(u|_\mathcal{J}\right)^\mathrm{dd}$. 
For every $w\in \Phi_\mathcal{J}\left[\bigoplus_{i\in \mathcal{J}} V_i\right]$ there exists $\lambda$ such that $-\lambda u\le w\le \lambda u$, hence 
$-\lambda \Phi_\mathcal{J}\left(u|_\mathcal{J}\right)\le w\le \lambda \Phi_\mathcal{J}\left(u|_\mathcal{J}\right)$ and, therefore, $w\in \Phi_\mathcal{J}\left(u|_\mathcal{J}\right)^\mathrm{dd}$.
\end{proof}

\begin{remark}\label{exa.antilatticeisirredicible}
Let $(V,C,u)$ be an order unit space which is also an anti-lattice. Then $V$ is irreducible. Indeed, %by \cite[Theorem 4.1.10]{Kalauch-vanGaans-book}, 
there are no non-trivial positive disjoint elements in $V$, hence no non-trivial directed bands. Therefore,  $V$ is irreducible.
\end{remark}

A converse of \Cref{pro.componentsarebands} is true if the components in the order direct sum are disjointness free anti-lattices.

\begin{theorem}\label{T:bands_in_atomic_JBW}
%Let $M=\bigoplus_{i\in I}M_i$ be an atomic JBW-algebra with the corresponding factor decomposition. 
Let $((V_i,C_i,u_i))_{i\in \mathcal{I}}$ be a collection of order unit spaces that are disjointness free anti-lattices with order direct sum $(V,C,u)$. 
\begin{itemize}
\item[(i)] $B\subseteq V$ is a band if and only if there exists $\mathcal{J}\subseteq \mathcal{I}$ such that $B=\Phi_\mathcal{J}\left[\bigoplus_{j\in \mathcal{J}}V_j\right]$, where it is understood that $B=\{0\}$ for $\mathcal{J}=\emptyset$.
\item[(ii)] Two non-zero $x,y\in V$ are disjoint if and only if there is $\mathcal{J}\subseteq \mathcal{I}$ with $\mathcal{J}\neq \emptyset$ and  $\mathcal{I}\setminus \mathcal{J}\neq \emptyset$ such that $x\in \Phi_\mathcal{J}\left[\bigoplus_{i\in \mathcal{J}} V_i\right]$ and $y\in \Phi_{\mathcal{I}\setminus \mathcal{J}}\left[\bigoplus_{i\in \mathcal{I}\setminus \mathcal{J}} V_i\right]$.
\end{itemize}
\end{theorem}
\begin{proof}
$(i)$ Let $B$ be a band in $V$. Define 
\[\mathcal{J}:=\{i\in \mathcal{I}\colon\, \mbox{there is a}\ v\in B\mbox{ such that }v(i)\neq 0\}.\] 
Let $v\in B$. For every $i\in \mathcal{I}\setminus \mathcal{J}$ we have $v(i)=0$, hence $v\in \Phi_\mathcal{J}(v|_\mathcal{J})\in \Phi_\mathcal{J}\left[\bigoplus_{j\in \mathcal{J}} V_j\right]$. For the converse inclusion, first observe that for every $i\in \mathcal{J}$ we have 
\[\{x(i)\colon\, x\in B\}^\mathrm{d}=\{0\}.\]
Indeed, there exists $v\in B$ with $v(i)\neq 0$. As $V_i$ is disjointness free, we have 
$\{x(i)\colon\, x\in B\}^\mathrm{d}\subseteq \{v(i)\}^\mathrm{d} =\{0\}$. 

Next, let $v\in \Phi_\mathcal{J}\left[\bigoplus_{j\in \mathcal{J}} V_j\right]$. Let $z\in B^\mathrm{d}$ and let $i\in \mathcal{I}$. If $i\in \mathcal{I}\setminus \mathcal{J}$, then $v(i)=0$, so $v(i)$ and $z(i)$ are disjoint. If $i\in \mathcal{J}$, then, by \Cref{lem.entrywisedisjointness}, $z(i)\in\{x(i)\colon\, x\in B\}^\mathrm{d}$, hence $z(i)=0$, so that $v(i)$ and $z(i)$ are disjoint. By \Cref{lem.entrywisedisjointness}, we obtain that $v$ and $z$ are disjoint. Thus, $v\in B^\mathrm{dd}=B$.

$(ii)$ Let $x,y\in V$ be disjoint. Define $\mathcal{J}:=\{i\in \mathcal{I}\colon\, x(i)\neq 0\}$. Then $x\in \Phi_\mathcal{J}\left[\bigoplus_{i\in \mathcal{J}} V_i\right]$. By \Cref{lem.entrywisedisjointness}, we have $y(i)=0$ for every $i\in \mathcal{J}$, as $V_i$ is disjointness free. Therefore, $y\in \Phi_{\mathcal{I}\setminus \mathcal{J}}\left[\bigoplus_{i\in \mathcal{I}\setminus \mathcal{J}} V_i\right]$.
%Let $B$ be a band in $M$. For each $i\in I$ let $c_i$ be the central projection such that $c_iM=M_i$.  Then $B\cap M_i$ is a band in $M_i$ for all $i\in I$ by Lemma[direct sums and bands]. But by \Cref{P: disjoint in B(H_q)}, \Cref{P:disjoint_in_spin}, \Cref{P:disjoint_B(H)}, and \Cref{P:disjoint_M_3} we have that $B\cap M_i=M_i$ or $B\cap M_i=\{0\}$. Hence, for 
%\[
%J:=\{i\in I\colon B\cap M_i=M_i\}
%\]
%it follows that $B=\bigoplus _{i\in J}M_i$. On the other hand, if $J\subset I$, then $B:=\bigoplus _{i\in J}M_i$ is a band by Theorem[Onno]. 
\end{proof}

\Cref{T:bands_in_atomic_JBW} yields a characterisation of disjointness and bands atomic JBW-algebras.

\begin{theorem}\label{T:disjointness_in_atomic_JBW}
Let $M=\bigoplus_{i\in \mathcal{I}}M_i$ be an atomic JBW-algebra with its factor decomposition given in \Cref{the.atomicdecomp}. 
\begin{itemize}
\item[(i)] $B\subseteq M$ is a band if and only if $B=\bigoplus_{j\in \mathcal{J}}M_j$ for $\mathcal{J}\subseteq \mathcal{I}$, where it is understood that $B=\{0\}$ for $\mathcal{J}=\emptyset$.
\item[(ii)] Two non-zero $x,y\in M$ are disjoint if and only if there is a $\mathcal{J}\subseteq \mathcal{I}$ with $\mathcal{J}\neq \emptyset$ and  $\mathcal{I}\setminus \mathcal{J}\neq \emptyset$ such that $x\in \bigoplus_{i\in \mathcal{J}} M_i$ and $y\in \bigoplus_{i\in \mathcal{I}\setminus \mathcal{J}} M_i$.
\end{itemize}
\end{theorem}

%\begin{proof}
%Let $x$ and $y$ be non-zero and disjoint. Then $I$ must have at least two elements by \Cref{P: disjoint in B(H_q)}, \Cref{P:disjoint_in_spin}, \Cref{P:disjoint_B(H)}, and \Cref{P:disjoint_M_3}. Let $c_i$ be the central projection in $M$ such that $c_iM=M_i$ and define $J$ to be the smallest subset of $I$ such that $x=(c_i x)_{i\in J}$. Note that $I\setminus J \neq \emptyset$ as otherwise $y\in M^{d}=\{0\}$. It follows that $y\in B^d$ where $B:=\bigoplus _{i\in J}M_i$ and as $B\cap B^d=\{0\}$, we have that $B^d\subset \bigoplus _{i\in I\setminus J}M_i$ by \Cref{T:bands_in_atomic_JBW}. Conversely, if $J\subset I$ with $J\neq \emptyset$ and  $I\setminus J\neq \emptyset$, then any $x\in \bigoplus_{i\in J} M_i$ and $y\in \bigoplus_{i\in I\setminus J} M_i$ are disjoint by \Cref{T:bands_in_atomic_JBW}.
%\end{proof}

\section{Inverses of disjointness preserving bijections}\label{sec.inverses}

In Banach lattices \cite{HuiPag1993} and in finite-dimensional pre-Riesz spaces \cite{KalLemGaa2019}, the inverse of a disjointness preserving bijection is disjointness preserving. We provide an example of a disjointness preserving bijection in an atomic JBW-algebra where the inverse is not disjointness preserving. Further, we characterise the disjointness preserving bijections with disjointness preserving inverse on order direct sums of disjointness free anti-lattices. 
\begin{example}
	Define $M$ to be the algebraic direct sum of the spin factor $\ell^2(\mathbb{N})\oplus \mathbb{R}$ and  $\ell^\infty(\mathbb{N})$. We have that $M$ is an atomic JBW-algebra. For $x:=(x_n)_{n \ge 1}\in \ell^2(\mathbb{N})$, $\alpha\in\mathbb{R}$, and $y:=(y_n)_{n \ge 1}\in\ell^\infty(\mathbb{N})$, define 
	\[T\colon M\to M,\quad (x,\alpha,y)\mapsto  \left(\, (x_2,x_3,x_4,\ldots),\alpha,(x_1,y_1,y_2,y_3,\ldots)    \,\right).\]
	Then $T$ is a linear map. Moreover, $T$ is disjointness preserving. Indeed, let $(x,\alpha,y)$, $(v,\beta,w)\in M$ be disjoint. Then $(x,\alpha)$ and $(v,\beta)$ are disjoint in $\ell^2(\mathbb{N})\oplus \mathbb{R}$ by \Cref{lem.entrywisedisjointness}. 
	By \Cref{P:disjoint_in_spin}, we get $(x,\alpha)=0$ or $(v,\beta)=0$. Without loss of generality assume $(x,\alpha)=0$. Then $(\,(x_2,x_3,x_4,\ldots),\alpha)=0$. Further, $(x_1,y_1,y_2,y_3,\ldots)=(0,y_1,y_2,y_3,\ldots)$ and $(v_1,w_1,w_2,w_3,\ldots)$ are disjoint in $\ell^\infty(\mathbb{N})$. Hence $T(x,\alpha,y)$ and $T(v,\beta,w)$ are disjoint.
	
Note that the map 
\[S\colon M\to M,\quad (x,\alpha,y)\mapsto  \left(\, (y_1,x_1,x_2,x_3,\ldots),\alpha,(y_2,y_3,y_4,\ldots)    \,\right)\]
is the inverse of $T$.
	The operator $S$ is not disjointness preserving. Indeed, \[a:=(\,(0,0,0,\ldots),1,(0,0,0,\ldots)\,) \ \mbox{ and } \ b:=(\,(0,0,0,\ldots),0,(1,0,0,0,\ldots)\,)\]
	are disjoint in $M$, but
	\[Sa=a \ \mbox{ and } \ Sb=(\,(1,0,0,0,\ldots),0,(0,0,0,\ldots)\,)\]
	are not disjoint in $M$, since $(\,(0,0,0,\ldots),1)$ and $(\,(1,0,0,0,\ldots),0)$ are both non-zero and, hence, cannot be disjoint in 
	 $\ell^2(\mathbb{N})\oplus \mathbb{R}$.
\end{example}

\begin{theorem}\label{the.characterizationofdisjointnesspreserving}
Let $((V_i,C_i,u_i))_{i\in \mathcal{I}}$ be a collection of order unit spaces that are disjointness free anti-lattices with order direct sum $(V,C,u)$. Let $T \colon V \to V$ be a disjointness preserving linear bijection. Then $T^{-1}$ is disjointness preserving if and only if there is a bijection $\sigma \colon \mathcal{I} \to \mathcal{I}$ and there are linear bijections $T_i \colon V_i \to V_{\sigma(i)}$ such that $T = \bigoplus _{i\in \mathcal{I}} T_i$. 
\end{theorem}
\begin{proof}
Assume that $T^{-1}$ is disjointness preserving. Fix $i\in \mathcal{I}$. Since $v:=T(\Phi_i(u_i))\neq 0$, there is $j\in \mathcal{I}$ such that $v(j)\neq 0$. Then $v(k)=0$ for every $k\neq j$. Indeed, the elements $\Phi_j(v(j))$ and $v-\Phi_j(v(j))$ are disjoint, hence $x:=T^{-1}(\Phi_j(v(j)))$ and $y:=T^{-1}(v-\Phi_j(v(j)))$ are disjoint. By \Cref{T:bands_in_atomic_JBW}$(ii)$, there is $\mathcal{K}\subseteq \mathcal{I}$ such that 
$$x\in \Phi_\mathcal{K}\left[\bigoplus_{k\in \mathcal{K}} V_k\right]\quad \mbox{and}\quad y\in \Phi_{\mathcal{I}\setminus \mathcal{K}}\left[\bigoplus_{i\in \mathcal{I}\setminus \mathcal{K}} V_i\right].$$ But $x+y=\Phi_i(u_i)$, so $x=0$ or $y=0$. Since $v(j)\neq 0$, we have $T(y)\neq v=T(\Phi_i(u_i))$, so $x\neq 0$, hence $y=0$. Thus, $v=\Phi_j(v(j))$ and, therefore, $v(k)=0$ for all $k\in I\setminus\{j\}$. Define $\sigma(i):=j$ and for every $w\in V_i$ define 
\begin{equation}\label{E:def of T_i}
T_i w:= T(\Phi_i(w))(j). 
\end{equation}
Clearly, $T_i\colon V_i\to V_{\sigma(i)}$ is linear. 

Next, we show that $\sigma\colon \mathcal{I}\to \mathcal{I}$ is a bijection and, for every $i\in \mathcal{I}$,  $T_i\colon V_i\to V_{\sigma(i)}$ is a bijection. Let $i\in \mathcal{I}$ and $j:=\sigma(i)$. We start by showing that $T$ maps $\Phi_i[V_i]$ into $\Phi_{j}[V_{j}]$. For every $w\in \Phi_i[V_i]$ and $z\in \Phi_j[V_j]^\mathrm{d}$ we have that $T(\Phi_i(u_i))$ and $z$ are disjoint, so $T^{-1}(z)$ and $\Phi_i(u_i)$ are disjoint. According to \Cref{pro.componentsarebands}$(iv)$, we have $\Phi_i(u_i)^\mathrm{dd}=\Phi_i[V_i]$. As $T^{-1}(z)\in \Phi_i(u_i)^\mathrm{dd}$, we obtain that $T^{-1}(z)$ and $w$ are disjoint. Therefore, $z$ and $T(w)$ are disjoint. Thus, $T(w)\in \Phi_j[V_j]^\mathrm{dd}=\Phi_j[V_j]$ by \Cref{T:bands_in_atomic_JBW}$(ii)$. Consequently, $T[\Phi_i[V_i]]\subseteq \Phi_[(V_j]$. Applying the same arguments to $T^{-1}$ instead of $T$, there exists $k\in \mathcal{I}$ such that $T^{-1}[\Phi_j[V_j]]\subseteq \Phi_k[V_k]$. Then $u_i=T^{-1}(Tu_i)\in \Phi_k(V_k)$, so $k=i$. It follows that $T_i$ is a bijection from $\Phi_i[V_i]$ onto $\Phi_j[V_{j}]$. Moreover, it follows that $\sigma$ is injective. To see that $\sigma$ is surjective, let $j\in \mathcal{I}$ and let $i\in \mathcal{I}$ be such that $T^{-1}(\Phi_j(u_j))\in \Phi_i[V_i]$. Then $T[\Phi_i[V_i]]\subseteq \Phi_j[V_j]$, hence $\sigma(i)=j$. 

It remains to show that $T=\bigoplus_{i\in \mathcal{I}} T_i$. Let $v\in V$ and $i\in \mathcal{I}$. Again, we write $j=\sigma(i)$. Since $v-\Phi_i(v(i))\in \Phi_i[V_i]^\mathrm{d}$ and $T$ is disjointness preserving, we obtain
\[T(v-\Phi_i(v(i)))\in (T[\Phi_i[V_i])^\mathrm{d}=\Phi_j[V_j]^\mathrm{d},\]
hence $(T(v-\Phi_i(v(i))))(j)=0$. Therefore, by \eqref{E:def of T_i},
\[(Tv)(j)=T(\Phi_i(v(i)))(j)=T_i v(i),\]
which shows that $T=\bigoplus_{i\in \mathcal{I}} T_i$.

For a proof of the converse implication, assume that there are a bijection $\sigma\colon \mathcal{I}\to \mathcal{I}$ and linear bijections $T_i \colon V_i\to V_{\sigma(i)}$ such that $T=\bigoplus_{i\in \mathcal{I}} T_i$. Let $v,w\in V$ be disjoint and denote $x:=T^{-1}(v)$ and $y:=T^{-1}(w)$. By \Cref{T:bands_in_atomic_JBW}$(ii)$, there exists $\mathcal{J}\subseteq \mathcal{I}$ such that $v\in \Phi_\mathcal{J}\left[\bigoplus_{i\in \mathcal{J}} V_i\right]$ and $w\in \Phi_{\mathcal{I}\setminus \mathcal{J}}\left[\bigoplus_{i\in \mathcal{I}\setminus \mathcal{J}} V_i\right]$. Denote $\mathcal{K}:=\sigma^{-1}[\mathcal{J}]$. Then $x\in \Phi_\mathcal{K}\left[\bigoplus_{i\in \mathcal{K}} V_i\right]$. Indeed, for every $i\in \mathcal{I}$, we have
\[\left(T\left(\Phi_\mathcal{K}(x|_\mathcal{K})\right)\right)(\sigma(i))=T_i\left(\Phi_\mathcal{K}(x|_\mathcal{K})(i)\right).\]
If $i\in \mathcal{K}$, then $\Phi_\mathcal{K}(x|_\mathcal{K})(i)=x(i)$, so 
\[T_i\left(\Phi_\mathcal{K}(x|_\mathcal{K})(i)\right) = T_i (x(i))=(Tx)(\sigma(i))=v(\sigma(i)).\]
If $i\in \mathcal{I}\setminus \mathcal{K}$, then $\Phi_\mathcal{K}(x|_\mathcal{K})(i)=0$, so
\[T_i\left(\Phi_\mathcal{K}(x|_\mathcal{K})(i)\right)=0=v(\sigma(i)),\]
since $\sigma(i)\in \sigma[\mathcal{I}\setminus \mathcal{K}]=\mathcal{I}\setminus\sigma[\mathcal{K}]=\mathcal{I}\setminus \mathcal{J}$. Hence,
\[T\left(\Phi_\mathcal{K}(x|_\mathcal{K})\right)(\sigma(i))=v(\sigma(i))= (Tx)(\sigma(i))\]
for every $i\in \mathcal{I}$. Thus, $\Phi_\mathcal{K}(x|_\mathcal{K})=x$, which yields that $x\in \Phi_\mathcal{K}\left[\bigoplus_{i\in \mathcal{K}} V_i\right]$. Similarly, it follows that $y\in \Phi_{\mathcal{I}\setminus \mathcal{K}}\left[\bigoplus_{i\in \mathcal{I}\setminus \mathcal{K}} V_i\right]$, so that $x$ and $y$ are disjoint. Consequently, $T^{-1}$ is disjointness preserving.
\end{proof}

\begin{corollary}
Let $M=\bigoplus_{i\in \mathcal{I}}M_i$ be an atomic JBW-algebra with the corresponding factor decomposition as in \Cref{the.atomicdecomp}, and $T \colon M \to M$ be a disjointness preserving linear bijection. Then $T^{-1}$ is disjointness preserving if and only if there is a bijection $\sigma \colon \mathcal{I} \to \mathcal{I}$ and there are linear bijections $T_i \colon M_i \to M_{\sigma(i)}$ such that $T = \bigoplus _{i\in \mathcal{I}} T_i$. 
\end{corollary}
Since every disjointness preserving bijection on a finite-dimensional order unit space has a disjointness preserving inverse \cite[Theorem 3.4]{KalLemGaa2019}, \Cref{the.characterizationofdisjointnesspreserving} also has the following consequence.

\begin{corollary}
Let $((V_i,C_i,u_i))_{i\in\{1,\ldots,n\}}$ be a collection of finite-dimensional order unit spaces that are disjointness free anti-lattices with order direct sum $(V,C,u)$. Let $T \colon V \to V$ be a disjointness preserving linear bijection. Then there is a bijection $\sigma \colon \{1,\ldots,n\} \to \{1,\ldots,n\}$ and, for every $i\in\{1,\ldots,n\}$ there is a linear bijection $T_i \colon V_i \to V_{\sigma(i)}$ such that $T = \bigoplus _{i\in\{1,\ldots,n\}} T_i$. 
\end{corollary}

\section{The algebraic and order theoretical centre of unital JB-algebras}\label{sec.centres}
Recall that the order theoretical centre of an order unit space $(V,C,u)$ is the partially ordered vector space
\[
\mathrm{E}(V) := \left\{ T \in \mathrm{B}(V) \colon -\lambda I \le T \le \lambda I\ \mbox{for some $\lambda >0$} \right\}.
\]
The operator norm on $\mathrm{E}(V)$ coincides with the order unit norm generated by the order unit $I$, see \Cref{cor.normsoncentresame}. In the theory of Riesz spaces, it is known that the order theoretical centre of a Riesz space is again a Riesz space, see \cite[Theorem~3.30]{Abramovich-Aliprantis}. It turns out that also for complete order unit spaces the order theoretical centre is a Riesz space.

\begin{proposition}\label{pro.EVJB}
Let $(V,C,u)$ be a complete order unit space. Then $\mathrm{E}(V)$ with composition is an associative unital JB-algebra, $\mathrm{E}(V)$ is isomorphic as a JB-algebra to a space of continuous functions on a compact Hausdorff space, and $\mathrm{E}(V)$ is a Riesz space.
\end{proposition}
\begin{proof}
As $V$ is complete, we have that $\mathrm{E}(V)$ is a Banach space by \Cref{cor.normsoncentresame}$(ii)$. With \Cref{cor.normsoncentresame}$(iii)$ and $(iv)$, it follows that $\mathrm{E}(V)$ is an associative JB-algebra. Clearly, $I\in \mathrm{E}(V)$. According to \Cref{the.representationJB}, $\mathrm{E}(V)$ is therefore isomorphic as a JB-algebra to a space of continuous functions $\mathrm{C}(\Omega)$. Since the positive cone of a JB-algebra consists of the squares of the algebra, $\mathrm{E}(V)$ and $\mathrm{C}(\Omega)$ are then also isomorphic as partially ordered vector spaces. Thus, $\mathrm{E}(V)$ is a Riesz space. 
\end{proof}

The order theoretical centre of an order unit space that is not complete need not be a Riesz space, as the next example shows.

\begin{example}
Let $V:=\mathrm{C}^1([0,1])$ with pointwise order and order unit $\one$. The inclusion of $V$ in $\mathrm{C}([0,1])$ is the functional representation of $V$. As $V$ is a subalgebra of $\mathrm{C}([0,1])$, \Cref{cor.ordercentreofsubalgebra} yields the identity $\mathrm{E}(V)=\{M_f\colon\, f\in \mathrm{C}^1([0,1])\}$. In particular, $\mathrm{E}(V)$ is not a Riesz space as the map $f\mapsto M_f$ is an order isomorphism.
\end{example}

%\Cref{pro.EVJB} yields that the set $Q$ of squares in $E(V)$ as in \Cref{rem.coneQ} is a cone, provided $V$ is complete. 

Let $A$ be a unital JB-algebra. By \Cref{pro.EVJB}, the order theoretical centre $\mathrm{E}(A)$ is isomorphic as a JB-algebra to a space of continuous functions on a compact Hausdorff space. According to \Cref{cor.algcentreisCspace}, the same is true for the algebraic centre $\mathrm{Z}(A)$. In this section, we show that the algebraic centre $\mathrm{Z}(A)$ and the order theoretical centre $\mathrm{E}(A)$ are isometrically isomorphic as JB-algebras. 

In what follows, we will denote the multiplication operator by an element $x \in A$ on $A$ by $T_x$, that is, $T_xy := x\circ y$ for all $y\in A$. The isomorphism between $\mathrm{Z}(A)$ and $\mathrm{E}(A)$ will be given by the map $z\mapsto T_z$. Thus, $\mathrm{E}(A)$ consists of multiplication operators. This extends a result from Banach lattice theory, where it is known that the order theoretical centre of a space of continuous functions on a compact Hausdorff space consist of all multiplication operators, see \cite[Theorem 3.32]{Abramovich-Aliprantis}. Let us collect some elementary properties.

\begin{lemma}\label{lem.easypropertiesmuop}
Let $A$ be a unital JB-algebra and let $f\colon \mathrm{Z}(A)\to \mathrm{B}(A)$ be given by $f(z):=T_z$ for all $z\in \mathrm{Z}(A)$. Then $f$ is linear, multiplicative, injective, and $f$ maps the algebraic unit to the identity operator. Moreover, $\|T_z\|=\|z\|$ for all $z\in \mathrm{Z}(A)$.
\end{lemma}
\begin{proof}
It is clear that $f$ is linear and that it maps the algebraic unit to the identity operator. Since $T_ae=T_be$ implies $a=b$ for all $a,b\in A$, we get that $f$ is injective. As mentioned in the proof of \cite[Proposition~1.52]{Alfsen}, we have that $T_{z \circ w} = T_z \circ T_w$ for any $z,w \in \mathrm{Z}(A)$, so $f$ is an algebra homomorphism. 

If $x\in A$ is such that $\|x\| \le 1$, then $\|T_zx\| = \|z\circ x\| \le \|z\|$, so $\|T_z\| \le \|z\|$. Conversely, $\|T_z\| \ge \|T_z e\| = \|z\|$. 
\end{proof}

The idea to show that $z\mapsto T_z$ maps $\mathrm{Z}(A)$ onto $\mathrm{E}(A)$ is to consider the JBW-algebra case first and investigate the order interval $[0,e]$ in the JBW-algebra and the operator interval $[0,I]$ in the order theoretical centre of the JBW-algebra. 
%The compactness of the operator interval combined with the fact that its extreme points are order projections yields the desired isomorphism. 
We extend the result to JB-algebras by passing to the bidual.

Let $M$ be a JBW-algebra. Recall that the operator interval $[0,I]$ consists of bounded linear operators $T\colon M \to M$ such that $0 \le T \le I$, see \eqref{eq.ordercentrebounded}.

The predual $M_*$ of $M$ generates the \emph{$\sigma$-weak operator topology} on $\mathrm{B}(M)$ by letting $T_i \to T$ if and only if $\varphi(T_ix) \to \varphi(Tx)$ for all $x \in M$ and all normal states $\varphi \in M_*$. Note that the interval $[0,I]$ is closed for the $\sigma$-weak operator topology. Indeed, if $T_i \to T$ for the $\sigma$-weak operator topology with $T_i$ in $[0,I]$ for all $i$, then, for any $x \in M_+$ and any normal state $\varphi$, we have $\varphi((I-T_i)x) \ge 0$ for all $i$, so $\varphi((I-T)x) \ge 0$ and $T \le I$ by \cite[Corollary~2.17]{Alfsen}. Similarly, we find that $T \ge 0$. It turns out that the operator interval $[0,I]$ is actually compact for the $\sigma$-weak operator topology, which is essentially  \cite[Remark~2.10(b)]{Choi-Kim}. We provide the details.

\begin{lemma}\label{L:operintiscompact}
Let $M$ be a JBW-algebra. The operator interval $[0,I]$ is compact for the $\sigma$-weak operator topology.
\end{lemma}
\begin{proof}
If we identify $M$ with $(M_*)^*$, then an operator $T \in \mathrm{B}(M)$ can be thought of as $Tx(\varphi) := \varphi(Tx)$ for all $x \in M$ and all $\varphi \in M_*$ by \cite[Corollary~2.17(2.11)]{Alfsen}. Define the linear map 
\[
\Phi \colon (M\hat{\otimes}M_*)^* \to \mathrm{B}(M) \qquad \mbox{($M\hat{\otimes}M_*$ is the projective tensor product)}
\]
by $\Phi(\psi)(x)(\varphi):= \psi(x \otimes \varphi)$ for all $x \in M$. We show that $\Phi$ is an isometry. Indeed, let $\varphi \in M_*$ be such that $\|\varphi\| \le 1$. It follows that $|\Phi(\psi)(x)(\varphi)| = |\psi(x \otimes \varphi)| \le \|\psi\|\|x\|$, so $\|\Phi(\psi)(x)\| \le \|\psi\|\|x\|$ and, therefore, $\|\Phi(\psi)\| \le \|\psi\|$. Let $\eta \in M \hat{\otimes} M_*$ and let 
\[
\eta = \sum_{k=1}^\infty x_k \otimes \varphi_k
\]
be a representation of $\eta$. For $\psi \in (M\hat{\otimes}M_*)^*$, it follows by continuity of $\psi$ that
\[
|\psi(\eta)| = \left|\sum_{k=1}^\infty \psi(x_k \otimes \varphi_k)\right| \le \sum_{k=1}^\infty |\psi(x_k \otimes \varphi_k)| \le \sum_{k=1}^\infty \|\Phi(\psi)(x_k)\|\|\varphi_k\| \le \|\Phi(\psi)\|\sum_{k=1}^\infty\|x_k\|\|\varphi_k\|.
\]
Taking the infimum over all such representations of $\eta$ yields $|\psi(\eta)| \le \|\Phi(\psi)\|\|\eta\|_\pi$ where $\|\cdot\|_\pi$ denotes the projective norm on $M\hat{\otimes}M_*$. Hence, $\|\psi\| \le \|\Phi(\psi)\|$ and we conclude that $\Phi$ is an isometry. 

Next we show that $\Phi$ is surjective. Indeed, let $T \in \mathrm{B}(M)$ and define the bilinear map $\vartheta \colon M \times M_* \to \mathbb{R}$ by $\vartheta (x,\varphi) := \varphi(Tx)$. By the universal property of the projective tensor product, there is a unique bounded linear map $\psi \colon M\hat{\otimes} M_* \to \mathbb{R}$ with $\|\psi\| = \|\vartheta\|$ such that $\psi(x \otimes \varphi) = \vartheta (x,\varphi) = \varphi(Tx)$. Hence, $\Phi(\psi) = T$. We conclude that $\Phi$ is an isometric isomorphism. 

Equip $\mathrm{B}(M)$ with the $\sigma$-weak operator topology and the norm dual of the projective tensor product $M \hat{\otimes} M_*$ with the weak* topology. If $(\psi_i)_i$ is a net in $M \hat{\otimes} M_*$ that converges weak* to $\psi$, then 
\[
\varphi(\Phi(\psi_i)(x)) = \Phi(\psi_i)(x)(\varphi) = \psi_i(x \otimes \varphi) \to \psi (x \otimes \varphi) = \varphi(\Phi(\psi)(x)),
\]
so $\Phi(\psi_i)$ converges to $\Phi(\psi)$ in the $\sigma$-weak operator topology. Hence, $\Phi$ is continuous with respect to the weak* topology and the $\sigma$-weak operator topology. It follows that $\Phi^{-1}\left[[0,I]\right]$ is a norm bounded and weak* closed set which is weak* compact by the Banach-Alaoglu theorem. Therefore, $[0,I]$ is compact for the $\sigma$-weak operator topology as it is the continuous image of a compact set.
\end{proof}

\begin{lemma}\label{L: U_p central}
Let $p$ be a projection in a JBW-algebra $M$. Then $U_p$ is in the operator interval $[0,I]$ if and only if $p$ is central.
\end{lemma}

\begin{proof}
Let $p\in M$ be a projection. Suppose $U_p$ is in $[0,I]$. Then we have $U_p x \le x$ for all positive $x \in M$, so $p$ operator commutes with all elements in $M_+$ by \cite[Lemma~1.48]{Alfsen}. As $M_+$ generates $M$, it follows that $p$ must be a central projection. Conversely, if $p$ is a central projection in $M$, then $U_p x \le x$ for all $x \in M_+$ by \cite[Lemma~1.48]{Alfsen}, so $U_p$ is in the operator interval $[0,I]$.
\end{proof}

The following lemma characterises the extreme points of the positive ball $[0,e]$ in a JBW-algebra. This is \cite[Proposition~1.40]{Alfsen}, but the proof there is not correct. We will provide an alternative argument here. 

\begin{lemma}\label{L: extreme [0e]}
Let $A$ be a JB-algebra with unit $e$. Then the extreme points of $[0,e]$ are precisely the projections in $A$. 
\end{lemma}

\begin{proof}
Let $p$ be a projection and $x,y\in [0,e]$ be such that $p = t x + (1-t)y$ for some $0 < t < 1$. Since $t x \le p$ and $(1-t)y \le p$, it follows that $U_{p^\perp}x = U_{p^\perp}y = 0$. Thus, $U_p x = x$, and $U_py = y$ by \cite[Proposition~1.38]{Alfsen}. Hence, $x = U_px \le U_pe = p = tx + (1-t)y$ and $x \le y$. Similarly, $y \le p$ implies that $y \le x$ and we find that $p$ is an extreme point of $[0,e]$.  

Let $x$ be an extreme point of $[0,e]$. Then $2x-x^2, x^2 \in [0,e]$ by the functional calculus \cite[Corollary~1.19]{Alfsen}, so that $x=\frac{1}{2}(2x-x^2)+\frac{1}{2}x^2$, hence $x=x^2$ and $x$ must be a projection.
\end{proof}

Next, we investigate the extreme points of the order interval $[0,I]$ in the space of bounded linear operators on a JBW-algebra $M$. We have that $0$ and $I$ are extreme points of $[0,I]$. Indeed, if $0 = t S + (1-t) T$ with $S,T\in [0,I]$ and $t\in (0,1)$, then $t S, (1-t)T \le 0$, so that $S = T = 0$. Similarly, if $I = t S + (1-t) T$ with $S,T\in [0,I]$ and $t\in (0,1)$, then $0 = t(I-S) +(1-t)(I-T)$ and both $I-S$ and $I-T$ are in $[0,I]$, so $I- S = I- T = 0$. Hence $S = T = I$. Therefore, $\{0,I\}$ is a subset of the set of extreme points of $[0,I]$.  We will see below that every extreme point in $[0,I]$ comes from a central projection. First, we need some more terminology concerning projection  operators on $M$.

A positive projection $P \colon M \to M$ with $\|P\|=1$ is called \emph{complemented} if there is a positive projection $Q \colon M\to M$ with $\|Q\| = 1$ and 
\[
M_+ \cap \ker P = M_+ \cap \mathrm{ran}\,Q.
\]
In this case, $P$ and $Q$ are said to be \emph{complementary}. For a positive projection $P\colon M\to M$, the dual operator $P^*\colon M^*\to M^*$ is a positive projection.

\begin{lemma}\label{lem.priortoextreme}
Let $M$ be a JBW-algebra and let $T\colon M\to M$ be a bounded linear operator. 
\begin{itemize}
\item[(i)] If $T\in [0,I]$, then $T^*[M_*]\subseteq M_*$. 
\item[(ii)] If $T^*[M_*]\subseteq M_*$, then $T$ is $\sigma$-weakly continuous.
\end{itemize}
\end{lemma}
\begin{proof}
$(i)$ Let $\varphi$ be a normal state. If $(x_i)_i$ decreases to 0, then $0 \le Tx_i \le x_i$, so $(Tx_i)_i$ decreases to 0, and $T^*\varphi(x_i) = \varphi(Tx_i) \to 0$. Hence $T^*\varphi$ is a normal state.

$(ii)$ Let $(x_i)_i$ and $x$ be such that $\varphi(x_i)\to \varphi(x)$ for every $\varphi\in M_*$. For each $\varphi\in M_*$, we have $T^*\varphi\in M_*$, so $\varphi(Tx_i)=T^*\varphi(x_i)\to T^*\varphi(x)=\varphi(Tx)$.
\end{proof}

For a bounded linear operator $T\colon M\to M$ with $T^*[M_*]\subseteq M_*$, let $T_*$ denote the restriction of $T^*$ to $M_*$. A $\sigma$-weakly continuous positive projection $P \colon M \to M$ is called \emph{bicomplemented} if there is a $\sigma$-weakly continuous positive projection $Q \colon M \to M$ such that $P$ and $Q$ are complementary, $P^*[M_*]\subseteq M_*$, $Q^*[M_*]\subseteq M_*$, and the projections $P_* \colon M_* \to M_*$ and $Q_* \colon M_* \to M_*$ are complementary. If $c$ is a central projection, then $U_c=T_c$; see \cite[(1.65) on p.\ 29]{Alfsen}.

\begin{proposition}\label{P: extreme in [0I]}
Let $M$ be a JBW-algebra. A bounded linear operator $T$ on $M$ is an extreme point of the interval $[0,I]$ if and only if it is of the form $T_p$ for some central projection $p$ in $M$.
\end{proposition}

\begin{proof}
Let $T \in [0,I]$ be an extreme point. Since $0\le T^2 \le T$ and $0\le 2T-T^2 \le I$, we have $T = \frac{1}{2}(2T-T^2) + \frac{1}{2}T^2$. Hence, $T=T^2$. Note that $T$ and $I-T$ are complementary. By \Cref{lem.priortoextreme}, the adjoint operators $T^*$ and $(I-T)^*$ map normal states to normal states and are $\sigma$-weakly continuous. Moreover, $T_*$ is in $[0_*,I_*]$ with $T_*^2 = T_*$, and $T_*$ is complemented by $I_* - T_*$. Hence $T$ is bicomplemented. By \cite[Theorem~2.83]{Alfsen}, it therefore follows that $T = U_p$ for some projection $p$ in $M$. By \Cref{L: U_p central}, we obtain that $p$ must be a central projection, so that 
$T = T_p$. 

Conversely, let $p\in M$ be a central projection. According to \Cref{L: U_p central}, we have  $T_p=U_p\in [0,I]$. Let $S,V\in [0,I]$ be such that $T_c = \frac{1}{2}S + \frac{1}{2}V$. Define $x := Se$, $y:=Ve$, $S_0:=S-T_x$, and $V_0:=V-T_y$. Note that $S_0e=V_0e = 0$. Evaluating $T_p$ at $e$ yields $p = \frac{1}{2} x + \frac{1}{2}y$. By \Cref{L: extreme [0e]}, we obtain $x=y=p$, so $S_0 = S-T_x = 2T_p-V -T_x=T_y-V= -V_0$. Further, as $T_p^2=T_p$, we have $T_p +T_pS_0 = T_p+ T_pS-T_p^2 =  T_pS \le T_p$, so $T_p S_0 \le 0$. Thus, $-T_pS_0$ is a positive operator which vanishes at the order unit $e$, so that $-T_pS_0=0$. But now, by \cite[Proposition~1.47]{Alfsen}, $S_0 = T_p S_0 + T_{p^\perp}S_0 = T_{p^\perp}S_0$ , hence $T_{p^\perp}S = T_{p^\perp}S_0 = S_0 \ge 0$. As $S_0e=0$, we get $S_0 = 0$. Hence $T_p$ is an extreme point of $[0,I]$. 
\end{proof}

Let $M$ be a JBW-algebra with algebraic centre $\mathrm{Z}(M)$. By  \cite[Proposition~2.36]{Alfsen}, the algebraic centre $\mathrm{Z}(M)$ is a JBW-subalgebra of $M$ with algebraic unit $e$. Denote $[0,e]_{\mathrm{Z}(M)} := [0,e] \cap \mathrm{Z}(M)$. 
%Indeed, since $x \in Z(M)$ if and only if $\mathrm{JB}(x,e) \subset Z(M)$, the spectrum of $x$ in $Z(M)$ is the same as the spectrum of $x$ in $M$. Hence $Z(M)_+ = Z(M) \cap M_+$ and so $[0,e]_{Z(M)} = [0,e] \cap Z(M)$ by \cite[Corollary~1.22]{Alfsen}. Furthermore, the centre $Z(M)$ is $\sigma$-weakly closed in $M$. To see this, suppose that $(x_i)_i$ is a net in $Z(M)$ that converges $\sigma$-weakly to some $x \in M$. Let $s$ be a symmetry in $M$, that is, $s^2 = e$. Then $U_s$ is $\sigma$-weakly continuous on $M$ by \cite[Proposition~2.4]{Alfsen} and so $U_sx_i \to U_sx$, but by \cite[Lemma~2.35]{Alfsen} we have that $x_i = U_sx_i$, so $U_sx = x$. Hence $x \in Z(M)$ by applying \cite[Lemma~2.35]{Alfsen} again. 
Since the unit ball $[-e,e]$ is $\sigma$-weakly compact in $M$ by the Banach-Alaoglu theorem, the interval $[0,e]$ is $\sigma$-weakly compact as well, as it is the image of $[-e,e]$ under the homeomorphism $x \mapsto \frac{1}{2}(x+e)$. Since Jordan multiplication is $\sigma$-weakly continuous in each separate variable, $\mathrm{Z}(M)$ is $\sigma$-weakly closed in $M$. We conclude that $[0,e]_{\mathrm{Z}(M)}$ is therefore $\sigma$-weakly compact in $M$.   

\begin{proposition}\label{T: characterization [0I]}
Let $M$ be a JBW-algebra with algebraic centre $\mathrm{Z}(M)$. The map $f \colon [0,e]_{\mathrm{Z}(M)} \to [0,I]$ defined by $f(z) := T_z$ is a homeomorphism with respect to the $\sigma$-weak topology on $M$ and the $\sigma$-weak operator topology on $[0,I]$.
\end{proposition}

\begin{proof}
We start by an observation that we need three times. If $(z_i)_i$ and $z$ in $M$ are such that $z_i$ converges $\sigma$-weakly to $z$, then $T_{z_i}$ converges to $T_z$ in the $\sigma$-weak operator topology. Indeed, Jordan multiplication is $\sigma$-weakly continuous in each separate variable, so $(z_i \circ x)_i$ converges $\sigma$-weakly to $z\circ x$. Hence, for every normal state $\varphi$, we have $\varphi(T_{z_i}x) = \varphi (z_i \circ x) \to \varphi(z\circ x) = \varphi(T_zx)$. 

We will show next that $f$ indeed maps into $[0,I]$. By 
\Cref{L: extreme [0e]}, every extreme point $z$ of $[0,e]_{\mathrm{Z}(M)}$ is a projection. As $z$ is also central, we have by \Cref{P: extreme in [0I]} that $T_z\in [0,I]$. Since $[0,I]$ is convex, for every convex combination $z$ of extreme points of $[0,e]_{\mathrm{Z}(M)}$, we have $T_z\in [0,I]$. The set $[0,e]_{\mathrm{Z}(M)}$ is convex and $\sigma$-weakly compact, so the Krein-Milman theorem yields that, for any $z\in [0,e]_{\mathrm{Z}(M)}$, there is a net $(z_i)_i$ that converges $\sigma$-weakly to $z$, where each $z_i$ is a convex combination of extreme points of $[0,e]_{\mathrm{Z}(M)}$. Then $T_{z_i}\in [0,I]$ for every $i$ and $T_{z_i}$ converges to $T_z$ in the $\sigma$-weak operator topology. By  \Cref{L:operintiscompact}, it follows that $T_z\in [0,I]$.  

From the observation at the beginning of the proof, it is clear that $f$ is continuous. To see that $f$ is injective, observe that $T_xe=T_ye$ implies $x=y$ for every $x,y\in M$. 

Next, we show that $f$ is surjective. For every extreme point $T$ of $[0,I]$, \Cref{P: extreme in [0I]} yields that $T=T_z$ for some central projection $z$ in $M$. Then $z\in [0,e]_{\mathrm{Z}(M)}$. As $[0,e]_{\mathrm{Z}(M)}$ is convex, for every convex combination $T$ of extreme points of $[0,I]$, there exists $z\in [0,e]_{\mathrm{Z}(M)}$ such that $T=T_z$. Let $T\in [0,I]$. According to \Cref{L:operintiscompact}, the convex set $[0,I]$ is compact. By the Krein-Milman theorem, there is a net $(T_i)_i$ in $[0,I]$ that converges to $T$ in the $\sigma$-weak operator topology, where each $T_i$ is a convex combination of extreme points of $[0,I]$. For each $i$, there is $z_i\in [0,e]_{\mathrm{Z}(M)}$ such that $T_i=T_{z_i}$.  Since $[0,e]_{\mathrm{Z}(M)}$ is $\sigma$-weakly compact, there is a subnet $(w_j)_j$ of $(z_i)_i$ that converges $\sigma$-weakly to an element $z$ of $[0,e]_{\mathrm{Z}(M)}$. But then $T_{w_j} \to T_z$. As $(T_{w_j})_j$ is a subnet of $(T_i)_i$, we obtain $T = T_z$. Thus, $f$ is surjective.

Since $[0,e]_{\mathrm{Z}(M)}$ is $\sigma$-weakly compact and $[0,I]$ is Hausdorff for the $\sigma$-weak operator topology, it follows that $f$ is a homeomorphism.  
\end{proof}

We are now in a position to show that the algebraic centre and the order theoretic centre of a JBW-algebra are isometrically isomorphic as JBW-algebras.

\begin{theorem}\label{T: centre is centre}
Let $M$ be a JBW-algebra. Consider its algebraic centre $\mathrm{Z}(M)$ and its order theoretical centre $\mathrm{E}(M)$ equipped with the order unit norm induced by $I$. 
\begin{itemize}
	\item[(i)] The map $f \colon \mathrm{Z}(M) \to \mathrm{E}(M)$ defined by $f(z):=T_z$ is a multiplicative isometric isomorphism.
	\item[(ii)] The order unit norm induced by $I$ and the operator norm coincide on $\mathrm{E}(M)$.
	\item[(iii)] $f$ is a homeomorphism if we equip $\mathrm{Z}(M)$ with the $\sigma$-weak topology and $\mathrm{E}(M)$ with the $\sigma$-weak operator topology.   
	\end{itemize}
\end{theorem}

\begin{proof} 
$(i)$ Due to \Cref{lem.easypropertiesmuop}, it remains to show that $f$ maps into and onto $\mathrm{E}(M)$, and that it is an isometry with respect to the order unit norm on $\mathrm{E}(M)$. We first check that it maps into $\mathrm{E}(M)$. If $z \in \mathrm{Z}(M)$ is non-zero, then $w := \frac{1}{2}(e+\|z\|^{-1}z) \in [0,e]_{\mathrm{Z}(M)}$. By \Cref{T: characterization [0I]}, we obtain $T_w \in [0,I]$,  hence $-\|z\|I \le T_z \le \|z\|I$, which means that $T_z\in \mathrm{E}(M)$.

To see that $f$ is surjective, let $T\in \mathrm{E}(M)$. For $\lambda > 0$ such that $-\lambda I \le T \le \lambda I$, we have  $\frac{1}{2}I +\frac{1}{2\lambda}T \in [0,I]$, so, by \Cref{T: characterization [0I]}, there is a $z_\lambda \in [0,e]_{\mathrm{Z}(M)}$ such that $\frac{1}{2}I +\frac{1}{2\lambda}T=T_{z_\lambda}$. Then $T = T_{\lambda(2z_\lambda-e)}$, which shows that  $f$ is surjective.

We show that $f$ is an isometry. From $-\|z\|I \le T_z \le \|z\|I$, it follows that $\|T_z\|_I \le \|z\|$, where $\left\|\cdot\right\|_I$ denotes the the order unit norm on $\mathrm{E}(M)$ induced by $I$. On the other hand, $-\|T_z\|_I I \le T_z \le \|T_z\|_I I$. By evaluating at $e$,  we find that $-\|T_z\|_I e\le z\le \|T_z\|_Ie$, so that $\|T_z\|_I\ge \|z\|$. 

$(ii)$ By $(i)$ and \Cref{lem.easypropertiesmuop}, we get $\|T_z\| = \|z\| = \|T_z\|_I$.

$(iii)$ Let $(z_i)_i$ be a net in $\mathrm{Z}(M)$ that converges $\sigma$-weakly to $z$. Then $z_i\circ x$ converges $\sigma$-weakly to $z\circ x$ for any $x \in M$ as Jordan multiplication is separately $\sigma$-weakly continuous.
% by \cite[Proposition~2.4]{Alfsen}, 
For any normal state $\varphi$, it follows that $\varphi(T_{z_i}x) \to \varphi(T_zx)$ and so $T_{z_i}$ converges to $T_z$ in the $\sigma$-weak operator topology. On the other hand, if $T_{z_i}$ converges to $T_z$ in the $\sigma$-weak operator topology, then, for any normal state $\varphi$, we find that $\varphi(z_i) = \varphi(T_{z_i}e) \to \varphi(T_ze) = \varphi(z)$, so $z_i \to z$ in the $\sigma$-weak topology. 
\end{proof}

The main step to obtain a result for the general case of JB-algebras is to modify \Cref{T: characterization [0I]} in the following way.

\begin{proposition}\label{T: characterization [0I] JB}
Let $A$ be a unital JB-algebra with algebraic centre $\mathrm{Z}(A)$. The map $f \colon [0,e]_{\mathrm{Z}(A)} \to [0,I]$ defined by $f(z) := T_z$ is a homeomorphism with respect to the weak topology on $A$ and the weak operator topology on $[0,I]$.
\end{proposition}

\begin{proof}
According to \cite[Corollary~2.50]{Alfsen}, the bidual $A^{**}$ of the JB-algebra $A$ is a JBW-algebra. Moreover, $A\cap \mathrm{Z}(A^{**})=\mathrm{Z}(A)$. Indeed, every element of $A\subseteq A^{**}$ that operator commutes with every element of $A^{**}$ also operator commutes with every element of $A$, so $A\cap \mathrm{Z}(A^{**})\subseteq\mathrm{Z}(A)$. To see that $\mathrm{Z}(A)\subseteq A\cap \mathrm{Z}(A^{**})$, let $z\in \mathrm{Z}(A)$ and let $x,y\in A^{**}$ be arbitrary. By Goldstine's theorem, $A$ is weak* dense in $A^{**}$, so there are nets $(x_i)_{i}$ and $(y_j)_{j}$ in $A$ that weak* converge to $x$ and $y$, respectively. Then $x_i\circ(z\circ y_j)=T_{x_i}T_zy_j=T_zT_{x_i}y_j=z\circ(x_i\circ y_j)$ for every $i$ and $j$. Due to \cite[Corollary~2.50]{Alfsen}, Jordan multiplication is separately weak* continuous on $A^{**}$. Hence, $x_i\circ(z\circ y)=z\circ(x_i\circ y)$ for all $i$ and, thus, $x\circ(z\circ y)=z\circ(x\circ y)$. The latter means that $z$ and $x$ operator commute in $A^{**}$, so that $z\in A\cap \mathrm{Z}(A^{**})$. 

We proceed by showing that $f$ maps indeed into $[0,I]$. Let $z \in [0,e]_{\mathrm{Z}(A)}$ and let $x\in A_+$. By the functional calculus \cite[Corollary~1.19]{Alfsen}, $x$ has a positive square root $x^{1/2}$ and the quadratic representation $U_{x^{1/2}}$ of $x^{1/2}$ is a positive operator. As $z$ is central, we have 
\begin{equation}\label{eq.gescharrelmetU}
U_{x^{1/2}}z=\{x^{1/2},z,x^{1/2}\}=2T_{x^{1/2}}T_z x^{1/2}-T_z T_{x^{1/2}}x^{1/2}=T_z T_{x^{1/2}}x^{1/2}=T_zx.
\end{equation}
It follows that $T_z x = U_{x^{1/2}}z \ge 0$. Also, $T_zx=U_{x^{1/2}}z \le U_{x^{1/2}}e = x$. Hence, $0 \le T_z \le I$. 

The map $f$ is injective, since $T_xe=T_ye$ implies $x=y$. To see that $f$ is surjective, let $T\in [0,I]$. Then the double adjoint operator $T^{**}$ of $T$ satisfies $0^{**} \le T^{**} \le I^{**}$, where $0^{**}$ and $I^{**}$ denote the zero operator and the identity operator on $A^{**}$, respectively. By \Cref{T: characterization [0I]}, there exists $z\in [0,e]_{\mathrm{Z}(A^{**})}$ such that $T^{**}x = z \circ x$ for all $x \in A^{**}$. Note that left multiplication by $z$ leaves $\mathrm{Z}(A^{**})$ invariant and that $T^{**}$ leaves $A$ invariant. Thus, $T^{**}$ leaves the intersection $A \cap \mathrm{Z}(A^{**}) = \mathrm{Z}(A)$ invariant. In particular, $z = T^{**} e \in \mathrm{Z}(A)$, since $e$ is also the identity of $A^{**}$ by \cite[Corollary~2.50]{Alfsen}. It follows that $T=T^{**}|_A=T_z$, so that $T = f(z)$. 

Next we show that $f$ is continuous. Let $(z_i)_i$ be a net in $[0,e]_{\mathrm{Z}(A)}$ that converges weakly to $z\in [0,e]_{\mathrm{Z}(A)}$. Then for any state $\varphi$ on $A$ and $x \in A_+$, it follows with the aid of \eqref{eq.gescharrelmetU} that 
\[
\left|\varphi(T_{z_i}x) - \varphi(T_zx)\right| = \left| \varphi(T_{z_i - z} x) \right| = \left| U^*_{x^{1/2}}\varphi(z_i - z) \right| \to 0,
\]
hence $T_{z_i} \to T_z$ for the weak operator topology. If $T_{z_i} \to T_z$ with respect to the weak operator topology in $[0,I]$, then for any state $\varphi$ on $A$, we have 
\[
\varphi(z_i) = \varphi(T_{z_i}e) \to \varphi(T_z e) = \varphi(z),
\]
so $z_i \to z$ weakly. We conclude that $f$ is a homeomorphism with respect to the weak topology on $A$ and the weak operator topology on $[0,I]$. 
\end{proof}

In the proof of \Cref{T: centre is centre}, we can now replace \Cref{T: characterization [0I]} by \Cref{T: characterization [0I] JB} and, thus, we obtain the following result on the algebraic centre and order theoretical centre of unital JB-algebras.

\begin{theorem}\label{T: centre is centre JB}
Let $A$ be a unital JB-algebra. Consider its algebraic centre $\mathrm{Z}(A)$ and its order theoretical centre $\mathrm{E}(A)$ equipped with the order unit norm induced by $I$. 
\begin{itemize}
	\item[(i)] The map $f\colon \mathrm{Z}(A) \to \mathrm{E}(A)$ defined by $f(z):=T_z$ is a multiplicative isometric isomorphism.
	\item[(ii)] The order unit norm induced by $I$ and the operator norm coincide on $\mathrm{E}(A)$.
	\item[(iii)] $f$ is a homeomorphism if we equip $\mathrm{Z}(A)$ with the weak topology and $\mathrm{E}(A)$ with the weak operator topology.   
	\end{itemize}
\end{theorem}

\section{The order theoretical centre of order unit spaces whose cone is positively spanned by extreme vectors} \label{sec.moreonordercentre}

In this section, we further investigate a certain class of complete order unit spaces (including all finite-dimensional ones) for which the order theoretical centre is $\mathbb{R}^n$. The order theoretical centre, just as in the case of (atomic) JBW-algebras, contains information about the decomposability of the order unit space in terms of order direct sums. In particular, it turns out that the order theoretical centre is isomorphic to $\mathbb{R}$ when the order unit space is irreducible or an anti-lattice. 

A non-zero element $p\in C$ is said to be \emph{extreme} if, for every $x\in C$ with $x\le p$, there exists $\lambda\ge 0$ such that $x=\lambda p$. The set of all extreme elements of $C$ is denoted by $\mathrm{ext}(C)$. We will consider order unit spaces with the property that $\mathrm{Span}_+\mathrm{ext}(C)=C$. The natural cone of $\mathrm{C}([0,1])$ does not have extreme elements, hence does not satisfy this property.

\begin{lemma}\label{lem.finitedimhastheextremeproperty}
Let $(V,C,u)$ be an order unit space. If $V$ is finite-dimensional, then $\mathrm{Span}_+\mathrm{ext}(C)=C$.
\end{lemma}
\begin{proof}
By \cite[Corollary 5.4.11]{KalvanGaa2019}, the dual space $V^*$ of $V$ is directed and then \cite[Proposition 1.5.13]{KalvanGaa2019} yields that the cone $C^*$ of $V^*$ has a non-empty interior. Then $C^{**}$ has a compact base by \cite[Theorem 1.5.21]{KalvanGaa2019} and, according to \cite[Lemma 2.6.8]{KalvanGaa2019}, we have that $C^{**}$ is naturally isomorphic to $C$. Thus, the cone $C$ has a compact base $S$. Due to \cite[Lemma 1.5.19]{KalvanGaa2019}, the extreme elements of $C$ correspond to the extreme points of $S$. Since $S$ is a compact convex set in a finite-dimensional space, Minkowski's theorem yields that $S$ equals the convex hull of its extreme points. It follows that $\mathrm{Span}_+\mathrm{ext}(C)=C$. 
\end{proof}

We consider a direct sum of a collection of order unit spaces $((V_i,C_i,u_i))_{i\in \mathcal{I}}$ and use the notations introduced below \eqref{eq.Visorderdirectsum}. Denote by $P_i$ the order projection onto $\Phi_i[V_i]$.
 
\begin{lemma}\label{lem.itsplits}
Let $((V_i,C_i,u_i))_{i\in \mathcal{I}}$ be a collection of order unit spaces with order direct sum $(V,C,u)$.
\begin{itemize}
\item[(i)] Let $W$ be a projection band in $V$ with $P$ the order projection onto $W$ such that $(W,C\cap W,Pu)$ is an irreducible order unit space. Then there exists $i\in \mathcal{I}$ such that $W\subseteq \Phi_i[V_i]$. 
\item[(ii)] If $p\in\mathrm{ext}(C)$, then there exists $i\in \mathcal{I}$ such that $p=P_i p$ and $p(i)\in \mathrm{ext}(C_i)$.
%\item[(ii)] If $P$ is an order projection, then $P(\Phi_\alpha(V_\alpha))\subseteq \Phi_\alpha(V_\alpha)$ for every $\alpha\in I$. Moreover, for every $\alpha\in I$ and every $v\in V$ we have $(Pv)(\alpha)=  (P\Phi_\alpha(v(\alpha)))(\alpha)$. 
%$P(V)=\bigoplus_{\alpha\in I} P_\alpha (P(\Phi_\alpha(V_\alpha)))$ is an order direct sum.
\item[(iii)] If $\mathrm{Span}_+\mathrm{ext}(C) = C$, then $\mathrm{Span}_+\mathrm{ext}(C_i) = C_i$ for every $i\in \mathcal{I}$.
\end{itemize}
\end{lemma}
\begin{proof}
$(i)$ For every $i\in \mathcal{I}$, we have that $P=P_i P+(I-P_i)P$. Since $W$ is irreducible and $P$ is the identity on $W$, we obtain that $P_i P=0$ or $(I-P_i)P=0$. There exists $i\in \mathcal{I}$ such that $P_i P\neq 0$ and then $P=P_i P$ and, therefore, $W\subseteq \Phi_i[V_i]$.

$(ii)$ For every $i\in \mathcal{I}$, we have $\Phi_i(p(i))\le p$, so there is $\lambda_i\ge 0$ with $\Phi_i(p(i))=\lambda_i p$. Using that $p\neq 0$, we choose $i\in \mathcal{I}$ with $p(i)\neq 0$. Then $\lambda_i\neq 0$. For all $j\neq i$ we have that $\Phi_j(p(j))$ and $\Phi_i(p(i))$ are disjoint, so that $\lambda_j p$ and $\lambda_i p$ are disjoint, which yields that $\lambda_j=0$. Therefore, $p(j)=0$. Hence $p=P_i p$. 

If $x\in C_i$ is such that $x\le p(i)$, then there is $\lambda\ge 0$ such that $\Phi_i(x)=\lambda p$ as $p$ is extreme in $C$. Hence, $x=(\Phi_i(x))(i)=\lambda p(i)$. Thus, $p(i)$ is extreme in $C_i$.

%(ii) Let $\alpha\in I$. To see that $P(\Phi_\alpha(V_\alpha))\subset \Phi_\alpha(V_\alpha)$, it suffices to consider $v\in C_\alpha$, since $V_\alpha$ is directed. As $P$ and $I-P$ are positive, we have $0\le P(\Phi_\alpha(v))\le \Phi_\alpha(v)$ and $0\le (I-P)(\Phi_\alpha(v))\le \Phi_\alpha(v)$. From the latter we get $(I-P)(\Phi_\alpha(v))(\beta)=0$ for all $\beta\in I\setminus\{\alpha\}$, hence $P(\Phi_\alpha(v))(\beta)=(\Phi_\alpha(v))(\beta)$. 

%Since $\Phi_\alpha(V_\alpha)$ is a directed ideal in $V$, it follows that $P(\Phi_\alpha(v))\in \Phi_\alpha(V_\alpha)$. 

%If $v\in V$ and $v=\sum_{\alpha\in I} v_\alpha$ with $v_\alpha\in V_\alpha$ and $v_\alpha=0$ for all but finitely many $\alpha$, then $Pv=\sum_{\alpha\in I} Pv_\alpha$ with $Pv_\alpha=0$ for all but finitely many $\alpha$. 

%Let $\alpha\in I$ and $v\in V$. 

%If for every $\alpha\in I$ we have $P_\alpha(P(\Phi_\alpha(v_\alpha)))\ge 0$, then ??????. If $Pv\ge 0$, then there are $w_\alpha\in V_\alpha$ with $w_\alpha=0$ for all but finitely many $\alpha$ such that $Pv=\sum_{\alpha\in I} w_\alpha$ and $w_\alpha\ge 0$. Then $Pv=P(Pv)=\sum_{\alpha\in I} P w_\alpha$ and $Pw_\alpha\ge 0$ for all $\alpha\in I$.

$(iii)$ Let $v\in C_i$. By assumption, there are $p_1,\ldots, p_n\in\mathrm{ext}(C)$ and $\lambda_1,\ldots,\lambda_n\ge 0$ such that $\Phi_i(v)=\sum_{k=1}^n \lambda_k p_k$. By $(ii)$, for every $k\in \{1,\ldots,n\}$, there exists $i_k\in \mathcal{I}$ such that $p_k=P_{i_k}p_k$ and $p_k(i_k)\in \mathrm{ext}(C_{i_k})$. Then $v=(\Phi_i(v))(i)= \sum_{k=1}^n \lambda_k p_k(i)$ and, for every $k\in\{1,\ldots,n\}$ with $i_k\neq i$, we have $p_k(i)=0$. Hence, $v\in\mathrm{Span}_+\mathrm{ext}(C_i)$. 
\end{proof}

\begin{lemma}\label{L:finte decomposition}
If $(V,C,u)$ is an order unit space such that $\mathrm{Span}_+\mathrm{ext}(C) = C$, then it is the order direct sum of finitely many irreducible order unit spaces. Moreover, this decomposition is unique up to possibly reordering the indices. 
\end{lemma}

\begin{proof}
By assumption, we can write the order unit $u = \lambda_1p_1+\dots+\lambda_np_n$, where $\lambda_k>0$ and $p_k \in \mathrm{ext}(C)$. Suppose that $V$ is reducible and that $V=\bigoplus_{i\in \mathcal{I}} V_i$ is an order direct sum of order unit spaces $((V_i,C_i,u_i))_{i\in \mathcal{I}}$. By \Cref{lem.itsplits}$(ii)$, for every $k\in \{1,\ldots,n\}$, there exists $i_k\in \mathcal{I}$ such that $p_k=P_{i_k}p_k$ and $p_k(i_k)\in \mathrm{ext}(C_{i_k})$. Let $\mathcal{J}=\left\{i_k\colon\, k\in\{1,\ldots,n\}\right\}$ and let $P$ be the order projection onto the band $\Phi_\mathcal{J}\left[\bigoplus_{i\in \mathcal{J}} V_{i}\right]$. Then $I-P$ is a positive order projection as well and the operator norm satisfies $\|I-P\|=\|(I-P)u\|=0$, so $P = I$. Hence, $V=\Phi_\mathcal{J}\left[\bigoplus_{i\in \mathcal{J}} V_{i}\right]$. Note that, for every $i\in \mathcal{I}\setminus \mathcal{J}$ and every $x\in V_i$, we have that $u$ and $\Phi_i(x)$ are disjoint. Hence, $V_i=\{0\}$ for every $i\in \mathcal{I}\setminus \mathcal{J}$. It follows that at most $n$ of the summands $V_i$ are non-zero. 
 
 Thus, each decomposition of $V$ into a direct sum of order unit spaces has at most $n$ non-trivial summands. By splitting up reducible summands inductively, it follows that $V$ equals the direct sum of finitely many irreducible order unit spaces.  

Suppose that $((V_i,C_i,u_i))_{i\in\{1,\ldots,m\}}$ and $((W_j,C_j,u_j))_{j\in\{1,\ldots,n\}}$ are irreducible order unit spaces such that $V$ is isomorphic to $\bigoplus_{i=1}^m V_i$ and to $\bigoplus_{j=1}^n W_j$. Then there is a bipositive surjective linear map $\Psi\colon \bigoplus_{i=1}^m V_i\to \bigoplus_{j=1}^n W_j$. 

Fix $k\in\{1,\ldots,m\}$. Since $\Phi_k[V_k]$ is a projection band in $\bigoplus_{i=1}^m V_i$, we have that $\Psi[\Phi_k[V_k]]$ is a projection band in $\bigoplus_{j=1}^n W_j$.  Let $P$ be the order projection in $\bigoplus_{j=1}^n W_j$ onto $\Psi[\Phi_k[V_k]]$. Then, by \Cref{lem.itsplits}$(i)$, there exists $j\in\{1,\ldots,n\}$ such that $\Psi[\Phi_k[V_k]]\subseteq \tilde{\Phi}_j[W_j]$, where $\tilde{\Phi}_j$ denotes the natural embedding of $W_j$ into $\bigoplus_{j=1}^n W_j$. 

Similarly, there exists $k'\in\{1,\ldots,n\}$ such that $\Psi^{-1}[\tilde{\Phi}_j[W_j]]\subseteq \Phi_{k'}[V_{k'}]$. Then
\[\Psi[\Phi_kV_k]]=P[\Psi[\Phi_k[V_k]]]\subseteq P[\tilde{\Phi}_j[W_j]]=P[\Psi\Psi^{-1}[\tilde{\Phi}_j[W_j]]]\subseteq P[\Psi[\Phi_{k'}[V_{k'}]]], \]
hence $P[\Psi[\Phi_{k'}[V_{k'}]]]\neq\{0\}$. It follows that $k'=k$. Indeed, if $k'\neq k$, then $\Phi_{k'}[V_{k'}]$ and $\Phi_k[V_k]$ are disjoint, hence $\Psi[\Phi_{k'}[V_{k'}]]$ and $\Psi[\Phi_k[V_k]]$ are disjoint, so that $P[\Psi[\Phi_{k'}[V_{k'}]]]=\{0\}$, which is a contradiction. Then
\[\Psi[\Phi_k[V_k]]\subseteq \tilde{\Phi}_j[W_j]=\Psi\Psi^{-1}\tilde{\Phi}_j[W_j]\subseteq \Psi\Phi_k[V_k],\]
so $\Psi[\Phi_k[V_k]]=\tilde{\Phi}_j [W_j]$. For distinct elements $k$, we obtain distinct elements $j$ this way, hence $n\ge m$. We conclude that $n=m$ and that the decompositions as order direct sums are unique up to possibly reordering the indices.
\end{proof}

We proceed with preparations for the proof of \Cref{T:order centre of pos spanned V}. In that proof, to a positive element of the order theoretical centre, we will associate a positive bijection with a positive inverse. Hereby we will use the following lemma, which relies on completeness of the space. 

\begin{lemma}\label{L:centre is Banach}
The order theoretical centre $\mathrm{E}(V)$ of a complete order unit space $(V,C,u)$ equipped with the order unit norm $\left\|\cdot\right\|_I$ is a Banach space and if $T \in [0,I]$ is such that $\|T\|_I<1$, then $I-T$ is invertible with inverse $\sum_{k=0}^\infty T^k$.
\end{lemma}

\begin{proof}
Suppose that $(T_n)_{n\ge 1}$ is a Cauchy sequence of operators in $\mathrm{E}(V)$ for $\left\|\cdot\right\|_I$. Let $\varepsilon>0$. Then there is an $N\ge 1$ such that $-\eps I \le T_n-T_m \le \eps I$ whenever $n,m\ge N$. Hence, for $x\in C$, it follows that $-\eps\|x\|u \le (T_n-T_m)x \le \eps \|x\|u$, so $\|(T_n-T_m)x\|\le \eps\|x\|$ whenever $n,m\ge N$. This implies that $(T_nx)_{n\ge 1}$ is a Cauchy sequence in $V$. Consequently, we can define an additive positively homogeneous map $\tau \colon C \to V$ by $\tau(x):=\lim_{n\to\infty}T_nx$, which can be extended to a linear map $T\colon V\to V$ via $Tx=T(y-z):=\tau(y)-\tau(z)$ see \cite[Theorem 1.2.5]{KalvanGaa2019}, where $x=y-z$ is a difference of positive elements $y$ and $z$. Since $(T_n)_{n\ge 1}$ is  bounded with respect to $\left\|\cdot\right\|_I$, by say $M$, it follows that $-MI\le T_n \le MI$ for all $n\ge 1$, and, for any $x\in C$, we therefore have $-Mx \le T_nx \le Mx$. Taking the limit as $n \to \infty$ yields $-MI\le T\le MI$, so $T \in \mathrm{E}(V)$. For every $n,m\ge N$ and every $x\in C$, we have $-\varepsilon x\le (T_n-T_m)x\le \varepsilon x$, so $-\varepsilon x\le (T_n-T)x\le \varepsilon x$, hence $-\varepsilon I\le T_n-T\le \varepsilon I$. Therefore, $\left\|T_n-T\right\|_I\le \varepsilon$. Thus, $\mathrm{E}(V)$ is a Banach space.

Note that, for any positive $S,T\in \mathrm{E}(V)$, it follows that $\|ST\|_I\le\|S\|_I\|T\|_I$. So, if $T \in [0,I]$ is such that $\|T\|_I<1$, then $\sum_{k=1}^\infty T^k$ converges in the Banach space $(\mathrm{E}(V),\left\|\cdot\right\|_I)$ to an operator $R$. We have that $(I-T)\sum_{k=0}^N T^k\to (I-T)R$ as $N\to\infty$. Since  
$
I-(I-T)\sum_{k=0}^N T^k = T^{N+1}
$
for every $N$ and $T^{N+1}\to 0$ as $N\to\infty$, we obtain $I-(I-T)R=0$, so that $(I-T)R=I$. Similarly, $R(I-T)=I$. It follows that $(I-T)$ is invertible in $\mathrm{E}(V)$ with inverse $R$.
\end{proof}

\begin{lemma}\label{lem.wegetlambda}
Let $(V,C,u)$ be an order unit space and let $S\colon V\to V$ be a positive linear bijection with a positive inverse such that $S\le I$. For every $p\in \mathrm{ext}(C)$, there exists $\lambda>0$ such that $Sp=\lambda p$.
\end{lemma}
\begin{proof}
Since $S \in [0,I]$, for every $p \in \mathrm{ext}(C)$ we have $Sp\le p$, hence there is $\lambda >0$ with $Sp = \lambda p$. 
\end{proof}

In the next proposition, we deal with irreducible spaces and proceed with the general case in the subsequent theorem.

\begin{proposition}\label{P: irreducible trivila centre}
Let $(V,C,u)$ be a complete irreducible order unit space such that $\mathrm{Span}_+ \mathrm{ext}(C) = C$. Then the order theoretical centre of $V$ is isomorphic to $\mathbb{R}$ as partially ordered vector spaces. % $\langle I\rangle \cong \mathbb{R}$.
\end{proposition}

\begin{proof}
Let $T \in [0,I]$. Then $\|\frac{1}{2}T\|_I < 1$, so $S:=I-\frac{1}{2}T$ is invertible with positive inverse by \Cref{L:centre is Banach}.
% for some  where $q \in \mathrm{ext}(C)$ as well since $S$ must map extreme vectors of the cone to extreme vectors of the cone by Lemma~\ref{L:totally ordered extreme}. As $\lambda q = Sp \le p$, it follows that $\lambda q = \mu p$ for some $\mu>0$, but then $q=p$ since both $p$ and $q$ are normalised. 
%Let $\mathcal{B}$ be a maximal linearly independent set in $\mathrm{ext}(C)$ so that $\mathrm{Span}\,\mathcal{B} = V$, and 
%Let $q \in \mathrm{ext} (C)$ be such that $q = \lambda_1 p_1 + \dots +\lambda_np_n$ where the $p_k$ may be chosen to be linearly independent, and $\lambda_k>0$. Then there is a $\lambda>0$ such that $\lambda q= Sq=\lambda_1\mu_1p_1+\dots+\lambda_n\mu_np_n$, and hence, it follows that $\lambda = \mu_k$ for all $1\le k\le n$. 
For every $\lambda>0$, define the set
\[
K_\lambda := \left\{q \in \mathrm{ext}(C) \colon Sq = \lambda q\right\}.
\]
By \Cref{lem.wegetlambda}, $\mathrm{ext}(C) = \bigcup_\lambda K_\lambda$ and $K_\lambda$ and $K_\mu$ are disjoint whenever $\lambda \neq \mu$. We wish to show that, actually, there exists $\lambda>0$ such that $\mathrm{ext}(C)=K_\lambda$. Choose $\lambda>0$ such that $K_\lambda\neq\emptyset$. 
%Suppose that there are at least two distinct $K_\lambda$ and $K_\mu$ in this union. 
By assumption, every element of $C$ is a positive linear combination of elements of $\mathrm{ext}(C)$. Define $\pi\colon C \to \mathrm{Span}_+ K_\lambda$ by 
\[
\pi(\alpha_1 p_1+\dots+\alpha_np_n):=\sum_{k\in M} \alpha_k p_k,
\]
where $\alpha_k\ge 0$, $p_k\in\mathrm{ext}(C)$ for all $k$, and $M:=\{k\in\{1,\ldots,n\}\colon\, p_k\in K_\lambda\}$. We will first show that $\pi$ is well defined. Let $\alpha_1 p_1+\cdots+\alpha_np_n=\beta_1p_1+\cdots+\beta_np_n$, where $\alpha_k,\beta_k \ge 0$ and $p_k\in\mathrm{ext}(C)$ for all $k\in\{1,\ldots,n\}$. Note that we allow coefficients to be $0$, so that we may indeed assume that the same vectors $p_k$ appear in both positive linear combinations. Denote $N:=\{1,\ldots,n\}\setminus M$. Without loss of generality, we assume that all $p_k$ with $k\in N$ are linearly independent, where we have to allow that $\alpha_k,\beta_k\in\mathbb{R}$ for $k\in N$. For $k\in N$, let $\lambda_k>0$ be such that $p_k\in K_{\lambda_k}$. We have
\begin{equation}\label{eq.firstlinearrelation}
	\sum_{k\in M} \alpha_k p_k -\sum_{k\in M} \beta_kp_k = \sum_{k\in N} \beta_k p_k - \sum_{k\in N} \alpha_k p_k.
\end{equation}
If we apply $S$ to this equality and divide by $\lambda$ it follows that
\begin{equation}\label{eq.secondlinearrelation}
	\sum_{k\in M} (\alpha_k-\beta_k) p_k  = \sum_{k\in N} (\beta_k-\alpha_k) {\textstyle\frac{\lambda_k}{\lambda}}p_k,
	\end{equation}
and by subtracting \eqref{eq.firstlinearrelation} from \eqref{eq.secondlinearrelation} we obtain
\[\sum_{k\in N} (\beta_k-\alpha_k)\left({\textstyle\frac{\lambda_k}{\lambda}} -1\right) p_k=0.\]
Since all $p_k$ with $k\in N$ are assumed to be linearly independent and $\lambda_k\neq\lambda$ for all $k\in N$, it follows that $\beta_k-\alpha_k=0$ for all $k\in N$. Thus, by \eqref{eq.firstlinearrelation}, $\sum_{k\in M}\alpha_k p_k=\sum_{k\in M}\beta_k p_k$ and, therefore, $\pi$ is well-defined.

%Suppose $x \in C$ can be written as $x_\lambda + x_0$ and simultaneously as $y_\lambda + y_0$ where none of the extreme vectors in the linear combination of $x_0$ and $y_0$ contain elements of $K_\lambda$, and both $x_\lambda$, $y_\lambda$ are in $\mathrm{Span}_+ K_\lambda$. Suppose that $x_\lambda - y_\lambda$ is non-zero. Then we may assume without loss of generality that $x_\lambda - y_\lambda = \mu_1 q_1+\dots+\mu_nq_n$ where the $q_k$ are linearly independent. By applying $S$ and then dividing by $\lambda$, we see that
%\[
%x_\lambda - y_\lambda = \mu_1\lambda^{-1}\lambda_1q_1+\dots+\mu_n\lambda^{-1}\lambda_n q_n,
%\]
%so $\lambda = \lambda_k$ for all $1\le k\le n$, which contradicts the construction of $x_0$ and $y_0$. Hence, $x_\lambda = y_\lambda$ and $\pi$ is well defined. 
Note that $\pi$ is positively homogeneous. Furthermore, $\pi$ is additive. Indeed, let $x,y \in C$ and write $x = \alpha_1p_1+\dots +\alpha_np_n$ and $y = \beta_1p_1+\cdots+\beta_np_n$, where $p_k\in\mathrm{ext}(C)$ and $\alpha_k,\beta_k\ge 0$ for all $k$. With the set $M$ as defined above, we obtain
\[
\pi(x+y) = \sum_{k\in M}(\alpha_k+\beta_k)p_k = \pi(x) + \pi(y).
\]
It follows that we can extend $\pi$ to a positive linear operator $P \colon V \to \mathrm{Span}\, K_\lambda$ by writing any $x\in V$ as $x = y-z$ for $y,z \in C$ and defining $Px:=\pi(y)-\pi(z)$, see \cite[Theorem 1.2.5]{KalvanGaa2019}. Since we also have that $P \le I$, we obtain an order projection on $V$. Suppose there are $\lambda,\mu>0$ with $\lambda\neq \mu$ such that $K_\lambda,K_\mu\neq\emptyset$. Then $P\neq 0$ and $P\neq I$,  which contradicts the fact that $V$ is irreducible. We conclude that $\mathrm{ext}(C)=K_\lambda$ for some $\lambda>0$.  Then $S=\lambda I$ and therefore, $T=2(1-\lambda)I$. Hence, the map $t\mapsto tI$ is a bijection from $\mathbb{R}$ to the order theoretical centre of $V$, which clearly also is a linear order isomorphism.
\end{proof}

\begin{theorem}\label{T:order centre of pos spanned V}
Let $(V,C,u)$ be a complete order unit space such that $\mathrm{Span}_+\mathrm{ext}(C)=C$. There exist irreducible order unit spaces $(V_k,C_k,u_k)$, $k\in \{1,\ldots,n\}$, such that $V$ is isomorphic to the order direct sum $\bigoplus_{k=1}^nV_k$. Moreover, the order theoretical centre of $V$ is isomorphic to $\mathbb{R}^n$ as order unit spaces.
\end{theorem}

\begin{proof}
\Cref{L:finte decomposition} yields that $V$ is isomorphic to the order direct sum $\bigoplus_{k=1}^nV_k$ and also that it is unique up to possibly reordering the indices. If $T \in [0,I]$, then, by \Cref{L:centre is Banach}, we have that $S:=I-\frac{1}{2}T$ is a positive linear bijection with positive inverse and $S \in [0,I]$. By \Cref{lem.wegetlambda}, for any $p \in \mathrm{ext}(C)$, there is a $\lambda>0$ such that $Sp=\lambda p$. \Cref{lem.itsplits}$(iii)$ then yields that $S[\Phi_k[V_k]]\subseteq \Phi_k[V_k]$ for every $k$. Since $S$ is surjective, it follows for every $k$ that $S[\Phi_k[V_k]]= \Phi_k[V_k]$. For every $v\in V_k$, define $S_k(v):=S(\Phi_k(v))(k)$ and similarly, define $R_k(v):=S^{-1}(\Phi_k(v))(k)$. Then $R_k$ is the inverse of $S_k$, so that $S_k\colon V_k\to V_k$ is a positive bijection with a positive inverse. Also, $S_k\in [0,I_k]$, where $I_k$ is the identity on $V_k$. By \Cref{P: irreducible trivila centre}, there is a $0<\mu_k\le 1$ such that $S_k=\mu_k I_k$. It follows that $S=\bigoplus_{k=1}^n \mu_kI_k$, and so, $T=\bigoplus_{k=1}^n 2(1-\mu_k)I_k$. Conversely, for every $0\le \mu_1,\dots,\mu_n \le 1$, we have that the operator $T:=\bigoplus_{k=1}^n\mu_k I_k$ is in $[0,I]$, and we conclude that the order theoretical centre of $V$ and $\mathbb{R}^n$ are isomorphic as order unit spaces.
\end{proof}

For finite-dimensional spaces, either \Cref{pro.EVJB} or \Cref{T:order centre of pos spanned V} and \Cref{lem.finitedimhastheextremeproperty} together yield the following characterisation of the order theoretical centre.

\begin{corollary}\label{cor.finitedimensionalcentre}
Let $(V,C,u)$ be a finite-dimensional order unit space of dimension $d$. Then there exists $n\in\mathbb{N}$ with $n\le d$ such that the order theoretical centre of $V$ is isomorphic to $\mathbb{R}^n$ as order unit spaces.
\end{corollary}

According to \Cref{exa.antilatticeisirredicible}, anti-lattices are irreducible, so we have the following special case.

\begin{corollary}\label{cor.centreofantilattice}
Let $(V,C,u)$ be a finite-dimensional order unit space. If $V$ is an anti-lattice, then its order theoretical centre is isomorphic to $\mathbb{R}$ as order unit spaces.
\end{corollary}

\begin{example}
The order theoretical centre of the spin factor $\mathbb{R}^{n-1}\oplus\mathbb{R}$ is isomorphic to $\mathbb{R}$. Also the order theoretical centre of the space of symmetric $n\times n$-matrices with the cone of positive semi-definite matrices is isomorphic to $\mathbb{R}$. Indeed, both spaces are anti-lattices due to \Cref{P:disjoint_in_spin} and \Cref{P:disjoint_B(H)}, so the assertions follow from \Cref{cor.centreofantilattice}. Note that this conclusion also follows from \Cref{T: centre is centre}.
\end{example}

\appendix

\section{Appendix}\label{JBappendix}
This appendix contains a survey of the factors of atomic JBW-algebras. 

\subsection{Quaternionic Hilbert spaces and their operators}

By introducing the multiplication rules $i^2=j^2=k^2=ijk=-1$ on the symbols $i$, $j$, and $k$, we induce the structure of a unital associative algebra on the four dimensional real vector space 
\[
\mathbb{H}:=\bigl\{a1+bi+cj+dk\colon a,b,c,d\in\mathbb{R}\bigr\}
\]
where the general product distributes over the sum as usual, with unit $1$. This algebra $\mathbb{H}$ is referred to as the \emph{quaternions}, and $i$, $j$, and $k$ are called the imaginary units. The multiplication on $\mathbb{H}$ is not commutative as $ij=-ji=k$, $ki=-ik=j$, and $jk=-kj=i$. The algebraic centre of $\mathbb{H}$ equals $\mathbb{R}1$, and every non-zero quaternion is invertible as
\[
(a1+bi+cj+dk)(a1-bi-cj-dk)=(a^2+b^2+c^2+d^2)1.
\]
For $q=a1+bi+cj+dk$, the \emph{quaternionic conjugate} is defined to be $q^*:=a1-bi-cj-dk$, which defines an involution ${}^*\colon \mathbb{H}\to\mathbb{H}$ that reverses the order of multiplication. That is, for $q,r\in\mathbb{H}$, we have $(qr)^*=r^*q^*$. The \emph{real part} of $q$ is denoted by $\operatorname{Re}(q)$ and is given by $\operatorname{Re}(q):=\frac{1}{2}(q+q^*)=a1$. Furthermore, this conjugation gives rise to the multiplicative norm $|q|:=\sqrt{q^*q}$ on $\mathbb{H}$. Note that $\mathbb{R}1$ has been identified with $\mathbb{R}$ here.

An abelian group $(V,+)$ that admits a right action $\cdot\colon V\times\mathbb{H}\to V$ is called a \emph{quaternionic vector space} if the action distributes over $+$ in $V$ and the sum of quaternions. That is, we have $(v + w)\cdot q=v\cdot q + w\cdot q$, $v\cdot (q+r)=v\cdot q + v\cdot r$, $(v \cdot r)\cdot q=v\cdot (qr)$, and $v\cdot 1 =v$ for all $q,r\in\mathbb{H}$ and all $v,w\in V$. The reason for choosing a right action on $V$ is so that an $n\times n$ matrix $A$ acting on the left as usual on a vector $x\in \mathbb{H}^n$, that is, $x\mapsto Ax$, is now $\mathbb{H}$-linear. 

A \emph{quaternionic inner product} on $V$ is a $\mathbb{H}$-sesquilinear form $\langle\cdot,\cdot\rangle\colon V\times V\to \mathbb{H}$, so $\langle \cdot , \cdot \rangle$ satisfies $\langle u, v\cdot q +w\rangle=\langle u,v\rangle q+\langle u,w\rangle$ and $\langle v,w\rangle=\langle w,v\rangle^*$ for all $q\in\mathbb{H}$ and all $u,v,w\in V$, which in addition satisfies $\langle v,v\rangle\ge 0$ for all $v\in V$ and $\langle v,v \rangle=0$ if and only if $v=0$. It follows that $\|v\|:=\sqrt{\langle v,v \rangle}$ yields a norm on $V$, which we  prove next by using the quaternionic version of the Cauchy-Schwarz inequality.

\begin{lemma}[Cauchy-Schwarz]\label{L:Cauchy-Schwartz}
	Let $V$ be a quaternionic vector space equipped with a quaternionic inner product. Then, for any $v,w\in V$, it follows that $|\langle v,w \rangle|\le \|v\|\|w\|$ with equality if and only if $v$ and $w$ are $\mathbb{H}$-linearly dependent. 
\end{lemma}

\begin{proof}
	If $w=0$, the statement clearly holds, so we may assume that $w\neq 0$. Let $q:=\langle w,v \rangle\|w\|^{-2}$ and observe that
	\begin{align*}
		0 \le \|v-w\cdot q\|^2 & = \|v\|^2 - q^*\langle w,v \rangle  - \langle v,w \rangle q + q^* \|w\|^2 q \\
		& = \|v\|^2 - \frac{|\langle v,w \rangle|^2}{\|w\|^2} - \frac{|\langle v,w \rangle|^2}{\|w\|^2} + \frac{|\langle v,w \rangle|^2}{\|w\|^2}\\
		& = \|v\|^2 - \frac{|\langle v,w \rangle|^2}{\|w\|^2},
	\end{align*}
	and hence, we have $|\langle v,w \rangle |\le \|v\|\|w\|$. Moreover, in case of an equality above, it follows that $v=w\cdot q$ and if $v=w\cdot r$ for some $r\in\mathbb{H}$, then $|\langle v,w \rangle|=|\langle w,w\rangle ||r|=|r|\|w\|\|w\|=\|w\cdot r\|\|w\|=\|v\|\|w\|$. 
\end{proof}
\noindent
It now follows from \Cref{L:Cauchy-Schwartz} that, for $v,w\in V$, we have 
\[
\|v+w\|^2\le \|v\|^2+2|\langle v,w \rangle|+\|w\|^2\le (\|v\|+\|w\|)^2, 
\]
showing that the triangle inequality is satisfied. If $V$ is complete with respect to the norm $\|\cdot\|$, then $V$ is a \emph{quaternionic Hilbert space}. For more details, see \cite{Moretti-Oppio}. Quaternionic Hilbert spaces will from now on be denoted by $\mathcal{H}_q$. Most of the theory for Hilbert spaces passes over analogously to quaternionic Hilbert spaces as will be shown. Two vectors $v,w\in\mathcal{H}_q$ are said to be \emph{orthogonal} if $\langle v,w \rangle=0$, and similarly, the \emph{orthogonal complement} $S^\perp$ of a set is defined. A subset $B\subseteq\mathcal{H}_q$ is \emph{orthonormal} if the vectors in $B$ have norm one and are pairwise orthogonal. The Pythagorean theorem also holds for the quaternionic inner product. That is, if $v_1,\ldots,v_n\in\mathcal{H}_q$ are pairwise orthogonal, then 
\begin{align*}
	\|v_1 + \cdots + v_n\|^2=\|v_1\|^2 + \cdots + \|v_n\|^2\qquad\qquad \emph{\mbox{(Pythagorean identity)}}.
\end{align*}
The Pythagorean theorem can be used in turn to prove Bessel's inequality, stating that, for an orthonormal set $\{b_n \colon n\in \mathbb{N}\}$ and any $v\in\mathcal{H}_q$, we have

\begin{align*}
	\sum_{k=1}^\infty |\langle v,b_k \rangle|^2\le \|v\|^2\qquad\qquad \emph{\mbox{(Bessel's inequality)}}.
\end{align*}
We call a subset $B\subseteq\mathcal{H}_q$ an \emph{orthonormal basis} for $\mathcal{H}_q$ if it is a maximal orthonormal set. An application of Zorn's lemma tells us that every quaternionic Hilbert space $\mathcal{H}_q\neq\{0\}$ has an orthonormal basis. Indeed, the set
\[
\mathscr{B}:=\left\{B\subseteq\mathcal{H}_q\colon \mbox{$B$ is an orthonormal set}\right\}
\]
is non-empty and partially ordered by set inclusion. Let $(B_i)_i$ be a chain in $\mathscr{B}$. Then $B_0:=\bigcup_{i}B_i$ is an orthonormal set and Zorn's lemma implies that $\mathscr{B}$ contains a maximal element $B$. Hence, there is no non-zero vector $v\in\mathcal{H}_q$ such that $\langle v,b \rangle=0$ for all $b\in B$, showing that $B$ is an orthonormal basis for $\mathcal{H}_q$. By Bessel's inequality, for any $v\in\mathcal{H}_q$, there are at most countably many vectors in $B$ such that $\langle v,b \rangle\neq 0$, and it follows that 
\begin{align}\label{E:basis}
	v=\sum_{b\in B} b\cdot\langle v,b \rangle \qquad\mbox{and}\qquad\|v\|^2=\sum_{b\in B} |\langle v,b \rangle|^2.
\end{align}
For more details, see \cite[Lemma~I.4.12]{Conway} and \cite[Theorem~I.4.13]{Conway} (the arguments are the same for quaternionic Hilbert spaces).

An $\mathbb{H}$-linear operator $T \colon \mathcal{H}_q \to \mathcal{H}_q$ is bounded in the same way an operator between Hilbert spaces is bounded, so $T$ is \emph{bounded} if and only if $\sup\{\|Tv\| \colon \|v\|\le 1\}<\infty$, with this supremum denoted by $\|T\|$. Since multiplication in $\mathbb{H}$ is not commutative, the bounded operators on $\mathcal{H}_q$ can only be a equipped with the structure of a \emph{real} vector space. For $r\in \mathbb{R}$, the $\mathbb{H}$-linear operator $rT$ is defined by $rTx:=Tx\cdot r$, but linearity fails if $r$ is replaced by a general quaternion $q\in \mathbb{H}$. To illustrate this with an example, for an operator $T$ on $\mathcal{H}_q$, it follows that, for $S:=jT$, we have $S(v\cdot i)\neq Sv\cdot i$. The real vector space of bounded $\mathbb{H}$-linear operators is denoted by $\mathrm{B}(\mathcal{H}_q)$, and becomes a Banach space when equipped with the norm $T\mapsto \|T\|$. In order to define the quaternionic adjoint of an operator in $\mathrm{B}(\mathcal{H}_q)$, we need the quaternionic version of the Riesz representation theorem. Similarly, a $\mathbb{H}$-linear functional $\varphi \colon \mathcal{H}_q \to \mathbb{H}$ is \emph{bounded} if $\sup\{|\varphi(v)| \colon \|v\|\le 1\}<\infty$.

\begin{lemma}[Riesz representation theorem]\label{L:Riesz}
	If $\varphi \colon \mathcal{H}_q \to \mathbb{H}$ is a bounded $\mathbb{H}$-linear functional, then there exists a unique $v \in \mathcal{H}_q$ such that $\varphi(w)=\langle v,w \rangle$ for all $w \in \mathcal{H}_q$.
\end{lemma}

\begin{proof}
	Since $\varphi$ is continuous, it follows that $\ker \varphi$ is closed in $\mathcal{H}_q$. Clearly, if $\varphi=0$, then we can take $v=0$ to represent $\varphi$. So, suppose that $\varphi \neq 0$. By choosing an orthonormal basis for $\ker\varphi$, it follows that this basis is not maximal in $\mathcal{H}_q$, so there is a $v \in \ker\varphi^\perp$ such that $\varphi(v)=1$. For $w \in \mathcal{H}_q$, we have $\varphi(w-v\cdot \varphi(w))=0$, so $w-v\cdot \varphi(w) \in \ker\varphi$. Hence, $\langle v,w \rangle - \|v\|^2\varphi(w) = \langle v,w - v\cdot \varphi(w) \rangle = 0$ and so $\varphi(w) = \langle \|v\|^{-2}v,w \rangle$ for all $w\in \mathcal{H}_q$. Note that if $u \in \mathcal{H}_q$ is such that $\langle v,w \rangle = \langle u,w \rangle$ for all $w \in \mathcal{H}_q$, then, by choosing $w=u-v$, we get $\langle u-v, u-v \rangle = 0$, so $u=v$. 
\end{proof}

Given a bounded operator $T\in \mathrm{B}(\mathcal{H}_q)$, we have a well-defined linear operator $S \colon \mathcal{H}_q \to \mathcal{H}_q$ given by the relation $\langle w,Tv \rangle = \langle Sw,v \rangle$ as a consequence of \Cref{L:Riesz}, since, for any $v \in \mathcal{H}_q$, the map $w \mapsto \langle v,Tw \rangle$ is $\mathbb{H}$-linear and bounded by \Cref{L:Cauchy-Schwartz}. The properties of a quaternionic inner product also show that this operator $S$ is unique, and bounded with $\|S\|=\|T\|$. We say that $S$ is the \emph{quaternionic adjoint} of $T$, and will be denoted by $T^*$. It is again similar to the case of dealing with operators on a Hilbert space to find that $T^{**}=T$ and $\|T^*T\|=\|T\|^2$. An operator $T\in \mathrm{B}(\mathcal{H}_q)$ is \emph{self-adjoint} if $T^*=T$ and the subspace of self-adjoint operators will be denoted by $\mathrm{B}(\mathcal{H}_q)_\mathrm{sa}$. Since $\mathrm{B}(\mathcal{H}_q)_\mathrm{sa}$ is closed in $\mathrm{B}(\mathcal{H}_q)$, it is a Banach space as well. Note that $\langle Tv,v \rangle \in \mathbb{R}1$ for all $T\in \mathrm{B}(\mathcal{H}_q)_\mathrm{sa}$. An \emph{orthogonal projection} in $\mathrm{B}(\mathcal{H}_q)_\mathrm{sa}$ is an idempotent operator, and orthogonal projections are in bijection with the closed subspaces of $\mathcal{H}_q$. Indeed, for a closed subspace $\mathcal{K}\subseteq\mathcal{H}_q$, we have that $\mathcal{K}=\{0\}$ is precisely the range of the projection $P=0$. Otherwise, let $B$ be an orthonormal basis for $\mathcal{K}$, and, for any finite set $F\subseteq B$,  define $P_F \colon \mathcal{H}_q \to \mathcal{H}_q$ by 
\[
P_Fv:=\sum_{b\in F} b\cdot \langle b,v \rangle.
\]
Then, it follows that $P_F$ is linear, idempotent, and $\|P_F\|=1$. Furthermore, note that, for $v,w\in\mathcal{H}_q$, we have
\[
\langle P_Fv,w \rangle = \sum_{b\in F} \langle b\cdot \langle b,v \rangle ,w \rangle = \sum_{b\in F} (\langle w,b \rangle \langle b, v\rangle)^* = \sum_{b\in F} \langle v,b \rangle \langle b, w\rangle = \sum_{b\in F} \langle v,b\cdot \langle b,w \rangle \rangle = \langle v,P_Fw \rangle,
\]
hence, $P_F$ is self-adjoint. Since 
\[
Pv:=\sum_{b\in B} \langle v,b \rangle \cdot b
\]
is the limit of the net $\{P_F v \colon \mbox{$F\subseteq B$ finite}\}$, it follows that $P$ is an idempotent linear operator with $\|P\|=1$. Moreover, for $v,w\in\mathcal{H}_q$, we have
\[
\bigl| \langle Pv,w \rangle - \langle v,Pw \rangle \bigr| = \bigl| \langle (P-P_F)v,w \rangle - \langle v,(P-P_F)w \rangle \bigr| \le \|(P-P_F)v\|\|w\|+\|(P-P_F)w\|\|v\| 
\]
and hence, we find that $P$ is self-adjoint making it an orthogonal projection. By \eqref{E:basis}, the range of $P$ equals $\mathcal{K}$. The uniqueness of $P$ follows from the fact that any orthogonal projection $Q$ with range $\mathcal{K}$ must agree with $P$ as $\mathrm{ran}\, P = \mathrm{ran}\, Q = \mathcal{K}$ and $\ker P = \ker Q = \mathcal{K}^\perp$. Conversely, any orthogonal projection $P$ yields a closed subset $\mathrm{ran}\, P \subseteq \mathcal{H}_q$. 

The commutative bilinear product $\circ\colon \mathrm{B}(\mathcal{H}_q)_\mathrm{sa}\times \mathrm{B}(\mathcal{H}_q)_\mathrm{sa}\to \mathrm{B}(\mathcal{H}_q)_\mathrm{sa}$ defined by $T\circ S:=\frac{1}{2}(TS+ST)$ turns $\mathrm{B}(\mathcal{H}_q)_\mathrm{sa}$ into a real Jordan algebra and the norm $T\mapsto \|T\|$ satisfies $\|S\circ T\|\le \|S\|\|T\|$, and also $\|T^2\|=\|T^* T\|=\|T\|^2$ for all $S,T\in \mathrm{B}(\mathcal{H}_q)_\mathrm{sa}$. Furthermore, for $S,T\in \mathrm{B}(\mathcal{H}_q)_\mathrm{sa}$ and $v\in\mathcal{H}_q$ with $\|v\|\le 1$, it follows that 
\[
\|Tv\|^2=\langle  Tv,Tv\rangle \le  \langle  Tv,Tv\rangle + \langle  Sv,Sv\rangle = \langle T^2v,v\rangle + \langle S^2v,v\rangle = \langle (T^2+S^2)v,v\rangle \le \|T^2+S^2\| 
\]
by \Cref{L:Cauchy-Schwartz}, so $\|T^2\|\le \|T^2+S^2\|$, and we conclude that $\mathrm{B}(\mathcal{H}_q)_\mathrm{sa}$ is a JB-algebra with identity $I$. The \emph{spectrum} of an operator $T$ in $\mathrm{B}(\mathcal{H}_q)_\mathrm{sa}$ is denoted by $\sigma(T)$ and is defined to be 
\[
\sigma(T):=\{\lambda\in\mathbb{R}\colon T-\lambda I \mbox{ is not invertible in }\mathrm{B}(\mathcal{H}_q)_\mathrm{sa}\}.
\]
Note that real scalar multiples of $I$ are considered here as $\mathrm{B}(\mathcal{H}_q)_\mathrm{sa}$ is a real vector space. Note that by the functional calculus \cite[Corollary~1.19]{Alfsen}, the spectrum of an operator is never empty. The \emph{numerical range} of an operator $T \in \mathrm{B}(\mathcal{H}_q)_\mathrm{sa}$, which is defined to be $N(T):=\{\langle Tv,v \rangle\colon \|v\|=1\}$, is related to the spectrum of $T$ as follows.

\begin{lemma}\label{L:spectrum in numerical range}
	For $T\in \mathrm{B}(\mathcal{H}_q)_\mathrm{sa}$, we have $\sigma(T)\subseteq \overline{N(T)}$.
\end{lemma}

\begin{proof}
	If $\lambda$ is such that $T-\lambda I$ is not injective, then there is a normalised $v\in\mathcal{H}_q$ such that $Tv=\lambda v$, so $\langle Tv,v\rangle = \lambda$, and thus $\lambda \in N(T)$. If $T-\lambda I$ is not surjective, then there are two cases to distinguish, the range of $T-\lambda I$ is not dense in $\mathcal{H}_q$, and the range of $T-\lambda I$ is not closed but dense in $\mathcal{H}_q$. Firstly, suppose there is a normalised vector $v$ in the orthogonal complement of $\overline{\mathrm{ran}(T-\lambda I)}$. Then we find that $\langle Tv,v \rangle - \lambda = \langle (T-\lambda I)v,v \rangle = 0$, so $\lambda \in N(T)$. Secondly, if the range of $T-\lambda I$ is not closed but dense in $\mathcal{H}_q$, then there is no $\mu>0$ such that $\|(T-\lambda I)v\|\ge \mu\|v\|$ for all $v \in \mathcal{H}_q$ as the range is not closed, so there is a sequence of normalised vectors $(v_n)_{n\ge 1}$ in $\mathcal{H}_q$ such that $(T-\lambda I)v_n \to 0$. It follows that $\langle Tv_n,v_n \rangle - \lambda = \langle (T-\lambda I)v_n,v_n \rangle \to 0$, so $\lambda \in \overline{N(T)}$.
\end{proof}

The partial ordering on $\mathrm{B}(\mathcal{H}_q)_\mathrm{sa}$ can be formulated via the following equivalent properties.
\begin{lemma}\label{L:cone in B(H_q)_sa}
	For an operator $T\in \mathrm{B}(\mathcal{H}_q)_\mathrm{sa}$, the following statements are equivalent.
	\begin{itemize}
		\item[$(i)$] $\langle Tv,v \rangle \ge 0$ for all $v\in\mathcal{H}_q$.
		\item[$(ii)$] $T=S^2$ for some $S \in \mathrm{B}(\mathcal{H}_q)_\mathrm{sa}$.
		\item[$(iii)$] $\sigma(T)\subseteq [0,\infty)$.
	\end{itemize}
\end{lemma}

\begin{proof}
	$(i)\Longrightarrow (iii)$: If $\lambda < 0$ and $v \in \mathcal{H}_q$, then 
	\[
	\|(T-\lambda I)v\|^2 = \|Tv\|^2 -2\lambda \langle Tv,v \rangle + \lambda^2\|v\|^2 \ge \lambda^2\|v\|^2
	\]
	which implies that $T-\lambda I$ is injective. The same inequality shows that the range of $T$ is closed in $\mathcal{H}_q$, so we can define a left inverse $S \in \mathrm{B}(\mathcal{H}_q)$ of $T$ by $S := T'\oplus I_{\mathrm{ ran }\, T^\perp}$, where $T'$ is the inverse of $T$ restricted to the range of $T$. The fact that $T'$ is bounded follows from the inverse mapping theorem for bounded operators on Banach spaces \cite[Theorem~III.12.5]{Conway} where the argument also works for operators on $\mathcal{H}_q$. Since $T$ is self-adjoint, it follows that $T$ also has a right inverse. Note that $T^{-1}$ needs to be self-adjoint, as, for any $v \in \mathcal{H}_q$, we have 
	\[
	\langle T^{-1}v,v \rangle = \langle T^{-1}Tw,Tw \rangle = \langle w,Tw \rangle = \langle Tw,w \rangle = \langle v,T^{-1}v \rangle. 
	\]
	We conclude that $\lambda \notin \sigma(T)$. Hence, $\sigma(T)\subseteq [0,\infty)$.
	
	$(iii)\Longrightarrow (ii)$: The existence of an operator $S$ such that $T=S^2$ in this case follows from the functional calculus \cite[Corollary~1.19]{Alfsen}.
	
	$(ii)\Longrightarrow (i)$: If $S$ is such that $T=S^2$, then it follows directly that $\langle Tv,v \rangle = \langle S^2v,v \rangle = \langle Sv,Sv \rangle \ge 0$ for all $v\in \mathcal{H}_q$.
\end{proof}

\begin{lemma}\label{L: B(H_q) is monotone comp}
	The JB-algebra $\mathrm{B}(\mathcal{H}_q)_\mathrm{sa}$ is monotone complete.
\end{lemma}

\begin{proof}
	Without loss of generality, we may assume that $(T_i)_i$ is an increasing net such that $0\le T_i \le I$. By the functional calculus \cite[Proposition~1.12]{Alfsen}, we have $S^2\le S$ for all $0\le S\le I$, so that, for any $v\in\mathcal{H}_q$ and $i\le j$, we have
	\begin{align}\label{E:pw-limit T_i}
		\|(T_j-T_i)v\|^2=\langle (T_j -T_i)v,(T_j-T_i)v \rangle = \langle (T_j-T_i)^2 v,v \rangle \le \langle (T_j-T_i) v,v \rangle
	\end{align}
	by \Cref{L:cone in B(H_q)_sa}. Furthermore, as $\langle T_iv,v \rangle \le \|v\|^2$ for all $v\in \mathcal{H}_q $, the increasing net $(\langle T_iv,v \rangle)_i$ is Cauchy in $\mathbb{R}1$, which implies that $(T_iv)_i$ is a Cauchy net in $\mathcal{H}_q$ for every $v\in\mathcal{H}_q$ by \eqref{E:pw-limit T_i}. Hence, we can define a linear operator $T$ via the pointwise norm limits $Tv:=\lim_i T_iv$. It follows that $\|T\|\le 1$ and since we also have that $\langle T_iv,v \rangle \to \langle Tv,v \rangle$ for all $v\in \mathcal{H}_q$, we have $T\in \mathrm{B}(\mathcal{H}_q)_\mathrm{sa}$. Suppose that $S \in \mathrm{B}(\mathcal{H}_q)_\mathrm{sa}$ is such that $T_i\le S$ for all $i$. Then 
	\[
	\langle (S-T)v,v \rangle = \langle (S-T_i)v,v \rangle + \langle (T_i-T)v,v \rangle \ge \langle (T_i-T)v,v \rangle \to 0, 
	\]
	so $T\le S$ by \Cref{L:cone in B(H_q)_sa} and $T$ is the supremum of $(T_i)_i$.  
\end{proof}

Let $v \in \mathcal{H}_q$ be such that $\|v\|=1$. Then $\varphi_v\colon \mathrm{B}(\mathcal{H}_q)_\mathrm{sa}\to\mathbb{R}$ defined by $\varphi_v(T):=\langle Tv,v \rangle$ is a positive linear functional with $\varphi_v(I)=1$, so it is a state on $\mathrm{B}(\mathcal{H}_q)_\mathrm{sa}$. Furthermore, if $(T_i)_i$ is an increasing net with supremum $T$ in $\mathrm{B}(\mathcal{H}_q)_\mathrm{sa}$, then we saw in the proof of \Cref{L: B(H_q) is monotone comp} that $\varphi_v(T_i)=\langle T_iv,v \rangle \to \langle Tv,v \rangle = \varphi_v(T)$. Hence, we have that $\varphi_v$ is a normal state. These states are referred to as \emph{vector states}. 

\begin{lemma}\label{L:vector states are pure states}
	The vector states $\varphi_v$ on $\mathrm{B}(\mathcal{H}_q)_\mathrm{sa}$ are pure states.
\end{lemma}

\begin{proof}
	Let $v$ be a normalised vector and $\varphi_v$ be the corresponding vector state. Suppose that, for some $0<t<1$ and states $\psi$, $\eta$, we have $\varphi_v = t\psi+(1-t)\eta$. For the projection $P_vw:=v\cdot\langle v,w \rangle$, it follows that $\varphi_v(P_v)=1$ and as $0 \le \psi(P_v),\eta(P_v)\le 1$, since $0\le P_v\le I$, it follows that $\psi(P_v)=\eta(P_v)=1$. The symmetric bilinear form $(S,T) \mapsto \psi(S\circ T)$ is positive semi-definite, and so 
	\[
	|\psi(T\circ (I-P_v))|^2\le |\psi(T)||\psi(I-P_v)|=0 
	\]
	by the generalised Cauchy-Schwarz inequality, \cite[5.5.3]{Dudley}. Hence $\psi(T)=\psi(T\circ P_v)$ for all $T\in \mathrm{B}(\mathcal{H}_q)_\mathrm{sa}$. Using that $(P_v\circ T)\circ P_v = \frac{1}{2}P_vTP_v+\frac{1}{2}T\circ P_v$ and the fact that $P_vT P_v$ is self-adjoint, it follows that 
	\[
	\psi(T)=\psi((P_v\circ T)\circ P_v)={\textstyle\frac{1}{2}}\psi(P_vT P_v)+ {\textstyle\frac{1}{2}}\psi(T\circ P_v) = {\textstyle\frac{1}{2}}\psi(P_vT P_v)+ {\textstyle\frac{1}{2}}\psi(T),
	\]
	so that $\psi(T)=\psi (P_v T P_v)$. Since $(P_v T P_v)w=v\cdot \langle Tv,v \rangle\langle v,w \rangle = \varphi_v(T)P_v w$ for all $w\in\mathcal{H}_q$, we conclude that $\psi(T)=\psi(P_v T P_v)=\varphi_v(T)\psi(P_v)=\varphi_v(T)$ for all $T \in \mathrm{B}(\mathcal{H}_q)_\mathrm{sa}$, so $\psi=\varphi_v$. Hence, $\varphi_v = \psi = \eta$ and $\varphi_v$ is a pure state. 
\end{proof}

Note that if $T,S \in \mathrm{B}(\mathcal{H}_q)_\mathrm{sa}$ are such that $\langle (S-T)v,v \rangle = 0$ for all $v \in \mathcal{H}_q$, then $S\le T$ and $T\le S$ by \Cref{L:cone in B(H_q)_sa} and so $T=S$. We find that the vector states separate the points of $\mathrm{B}(\mathcal{H}_q)_\mathrm{sa}$ and hence, it follows that $\mathrm{B}(\mathcal{H}_q)_\mathrm{sa}$ is a JBW-algebra. 

We will show that $\mathrm{B}(\mathcal{H}_q)_\mathrm{sa}$ is an atomic JBW-algebra. Indeed, let $P$ be a non-zero orthogonal projection in $\mathrm{B}(\mathcal{H}_q)_\mathrm{sa}$. Then the range of $P$ is a closed subspace of $\mathcal{H}_q$, so we may choose an orthonormal basis for it. Let $v$ be an element of this orthonormal basis and note that the orthogonal projection $P_vw:=v\cdot \langle v,w \rangle $ satisfies $P_v=PP_v=P_vP$, so $P-P_v$ is idempotent and self-adjoint, so $\langle (P-P_v)w,w \rangle \ge 0$ for all $w\in \mathcal{H}_q$ and therefore $P_v\le P$. The following lemma shows that $P_v$ is an atom from which we can conclude that  $\mathrm{B}(\mathcal{H}_q)_\mathrm{sa}$ is atomic.

\begin{lemma}\label{L:P_v is minimal}
	For $v \in \mathcal{H}_q$ with $\|v\|=1$, the projection $P_v$ is an atom, and every atom in $\mathrm{B}(\mathcal{H}_q)_\mathrm{sa}$ is of this form. 
\end{lemma}

\begin{proof}
	Suppose that $P_v=Q+R$ for some orthogonal projections $Q$ and $R$. If $Qv$ and $Rv$ are non-zero, then they are $\mathbb{H}$-linearly independent and so the dimension over $\mathbb{H}$ of the range of $Q+R$ is at least two, which contradicts the fact that the range of $P_v$ equals $\mathbb{H}v:=\{v\cdot q\colon\,q\in\mathbb{H}\}$. Hence, we may assume without loss of generality that $Qv = v$ and $Rv=0$. But $RQR=0$ by \cite[Proposition~2.18]{Alfsen} and as 
	\[
	\|RQR\|=\|RQ^2R\|=\|(QR)^*QR\|=\|QR\|^2, 
	\]
	it follows that $QR=0$. Hence, $P_vw = Q(v\cdot \langle v,w \rangle )= Qw+QRw = Qw$, so $Q=P_v$ and $P_v$ is an atom.
	
	Let $P$ be an atom in $\mathrm{B}(\mathcal{H}_q)_\mathrm{sa}$. Then the range of $P$ must be of the form $\mathbb{H}v$ for some $v\in\mathcal{H}_q$ with $\|v\|=1$. If $w\in \mathcal{H}_q$ and $Pw=v\cdot q$, then it follows from $w=Pw+(I-P)w$ that $q=\langle v,w \rangle$ by taking the inner product with $v$ as $v\in \ker P^\perp$. Hence, $P=P_v$.
\end{proof}

\begin{lemma}
	The JBW-algebra $\mathrm{B}(\mathcal{H}_q)_\mathrm{sa}$ is a factor.
\end{lemma}

\begin{proof}
	Suppose $T\in \mathrm{B}(\mathcal{H}_q)_\mathrm{sa}$ operator commutes with all $S \in \mathrm{B}(\mathcal{H}_q)_\mathrm{sa}$. For any normalised $v\in\mathcal{H}_q$, the atom $P_v$ yields
	\[
	{\textstyle\frac{1}{2}}(Tv+v\cdot \langle  Tv,v\rangle )=T\circ(P_v\circ P_v)v= P_v\circ(T\circ P_v)v=v\cdot{\textstyle\frac{3}{4}}\langle Tv,v \rangle +{\textstyle\frac{1}{4}}Tv,
	\]
	hence, $Tv=v\cdot \langle Tv,v \rangle $. If $v$ and $w$ are two linearly independent normalised vectors in $\mathcal{H}_q$, then the linearity of $T$ implies that 
	\[
	(v+w)\cdot\langle T(v+w),v+w \rangle\|v+w\|^{-2}= T(v+w)=Tv+Tw=\langle Tv,v \rangle v+\langle Tw,w \rangle w
	\]
	and so $\langle Tv,v\rangle = \langle Tw,w\rangle$. Hence, it follows that $T=\lambda I$ for some $\lambda \in \mathbb{R}$ showing that the algebraic centre of $\mathrm{B}(\mathcal{H}_q)_\mathrm{sa}$ equals $\mathbb{R}I$.
\end{proof}

\subsection{Spin factors}
Let $H$ be a real Hilbert space of dimension at least two. If we equip the direct sum $H\oplus\mathbb{R}$ with the product
\begin{equation}\label{eq.spinproduct}
(x,\lambda)\circ(y,\mu):=(\mu x+\lambda y,\langle x,y \rangle +\lambda\mu),
\end{equation}
then $H\oplus\mathbb{R}$ becomes a Jordan algebra with unit $(0,1)$. The Cauchy-Schwarz inequality implies that we can define the norm $\|(x,\lambda)\|:=\sqrt{\langle x,x \rangle}+|\lambda|$ on $H\oplus\mathbb{R}$ and it follows that this defines a JB-algebra norm.

\begin{lemma}\label{L:spin factor norm}
	The norm on $H \oplus \mathbb{R}$ is a JB-algebra norm.
\end{lemma}

\begin{proof}
	Let $(x,\lambda)$ and $(y,\mu)$ be in $H \oplus \mathbb{R}$. Then 
	\begin{align*}
		\|(x,\lambda)\circ(y,\mu)\| &= \|(\mu x+\lambda y,\langle x,y \rangle + \lambda\mu)\| = \sqrt{ \mu^2\langle x,x  \rangle + 2\lambda\mu\langle x,y \rangle + \lambda^2\langle y,y \rangle} + |\langle x,y \rangle + \lambda\mu| \\&\le \sqrt{\mu^2\langle x,x  \rangle + 2|\lambda||\mu|\sqrt{\langle x,x \rangle} \sqrt{\langle y,y \rangle} + \lambda^2\langle y,y \rangle} + \sqrt{\langle x,x \rangle }\sqrt{\langle  y,y\rangle} +|\lambda||\mu| \\&= |\mu|\sqrt{\langle x,x \rangle} + |\lambda|\sqrt{\langle y,y \rangle} + \sqrt{\langle x,x \rangle }\sqrt{\langle  y,y\rangle} +|\lambda||\mu| \\&= (\sqrt{\langle x,x \rangle} + |\lambda|)(\sqrt{\langle y,y \rangle} + |\mu|) \\& = \|(x,\lambda)\|\|(y,\mu)\|,
	\end{align*}
	by the Cauchy-Schwarz inequality, and we can explicitly check the identity 
	\[
	\|(x,\lambda)^2\| = \|(2\lambda x,\langle x,x\rangle+\lambda^2)\| = 2|\lambda|\sqrt{\langle x,x \rangle}+\langle x,x \rangle + \lambda^2 = (\sqrt{\langle x,x \rangle}+|\lambda|)^2 = \|(x,\lambda)\|^2. \]
	Lastly, by using the Cauchy-Schwarz inequality once more, it follows that 
	\begin{align*}
		\|(x,\lambda)^2 + (y,\mu)^2\| &= \sqrt{ 4\lambda^2 \langle x,x \rangle +8\lambda\mu \langle x,y \rangle + 4\mu^2\langle y,y \rangle} + \langle x,x \rangle + \langle y,y \rangle + \lambda^2+ \mu^2
		\\&\ge \sqrt{ 4\lambda^2 \langle x,x \rangle -8|\lambda||\mu |\sqrt{\langle x,x \rangle}\sqrt{\langle y,y \rangle} + 4\mu^2\langle y,y \rangle} + \langle x,x \rangle + \langle y,y \rangle + \lambda^2+ \mu^2 \\&= 2|\lambda|\sqrt{\langle x,x \rangle} - 2|\mu|\sqrt{\langle y,y \rangle}  + \langle x,x \rangle + \langle y,y \rangle + \lambda^2+ \mu^2 \\&= (\sqrt{\langle x,x \rangle}+|\lambda|)^2 + (\sqrt{\langle y,y \rangle} - |\mu|)^2 \\&\ge \|(x,\lambda)^2\|.
	\end{align*}
	Hence this norm satisfies the properties of a JB-algebra norm.
\end{proof}

Furthermore, note that the Hilbert space direct sum norm $\|(x,\lambda)\|_2:=\sqrt{\langle x,x \rangle + \lambda^2}$ on $H\oplus \mathbb{R}$ is equivalent to $\left\|\cdot\right\|$. Indeed, we have 
\[
\|(x,\lambda)\|_2 \le \sqrt{\langle x,x \rangle} + |\lambda| = \|(x,\lambda)\|, 
\]
and by the concavity of the square root function, we find that
\[
\|(x,\lambda)\| = \sqrt{\langle x,x \rangle} + |\lambda| \le \sqrt{2}\sqrt{\langle x,x \rangle + \lambda^2} = \sqrt{2}\|(x,\lambda)\|_2.
\]
Hence, $H\oplus\mathbb{R}$ is reflexive by \cite[Proposition~1.11.8]{Megginson}. It follows that $H\oplus\mathbb{R}$ is a JB-algebra that is a dual space, so by \cite[Theorem~2.55]{Alfsen} it is a JBW-algebra. These JBW-algebras are called a \emph{spin factor}.

\begin{lemma}\label{L: projections in spin}
	The projections in $H\oplus\mathbb{R}$ are precisely $(0,0)$, $(0,1)$, and $(x,\frac{1}{2})$ such that $\sqrt{\langle x,x \rangle} = \frac{1}{2}$. Moreover, the latter are precisely the atoms.
\end{lemma}

\begin{proof}
	The equation $(x,\lambda)^2=(2\lambda x,\langle x,x \rangle +\lambda^2)$ yields $(2\lambda-1)x=0$ and $\lambda^2-\lambda + \langle x,x \rangle = 0$. So, if $x=0$, then $\lambda = 0$ or $\lambda = 1$ gives the idempotents $(0,0)$ and $(0,1)$, and if $x\neq 0$, then $\lambda = \frac{1}{2}$ and $\langle x,x \rangle = \frac{1}{4}$ as required. If $x$ is such that $\sqrt{\langle x,x \rangle} = \frac{1}{2}$, then $(x,\frac{1}{2}) + (-x,\frac{1}{2}) = (0,1)$, and it is clear that $(x,\frac{1}{2})$ can not be written as the sum of two projections of the form $(y,\frac{1}{2})$ or $(0,1)$. Hence, the atoms are precisely of the form $(x,\frac{1}{2})$.   
\end{proof}

By \cite[Lemma~1.10]{Alfsen}, spin factors are partially ordered by the cone of squares.

\begin{lemma}\label{L: cone of squares in spin id}
	The cone of squares $C$ in a spin factor $H \oplus \mathbb{R}$ equals $\Lambda:=\{(x,\lambda)\colon \sqrt{\langle x,x \rangle} \le \lambda\}$.
\end{lemma}

\begin{proof}
	Observe that $(0,\lambda)$ is a square if and only if $\lambda\ge 0$. Suppose that  $(x,\lambda)$ is not a multiple of $(0,1)$. Then $(x,\lambda)^2=(2\lambda x,\langle x,x \rangle +\lambda^2)$ satisfies $\sqrt{\langle 2\lambda x,2\lambda x \rangle} = 2|\lambda|\sqrt{\langle x,x \rangle} \le \langle x,x \rangle +\lambda^2$ as $(\sqrt{\langle x,x \rangle}-|\lambda|)^2\ge 0$ and it follows that $C\subseteq\Lambda$. Conversely, note that the two atoms $(x,\frac{1}{2})$ and $(-x,\frac{1}{2})$ are orthogonal since $(x,\frac{1}{2})\circ(-x,\frac{1}{2})=(0,0)$. So, if $(x,\lambda)$ is an element of $\Lambda$, then $(x,\lambda)$ is the square of 
	\[
	\sigma_1 \left(\frac{x}{2\sqrt{\langle x,x \rangle}},\frac{1}{2}\right) + \sigma_2\left(-\frac{x}{2\sqrt{\langle x,x \rangle}},\frac{1}{2}\right)
	\]
	where $\sigma_k:=\sqrt{\lambda + (-1)^{k+1}\sqrt{\langle x,x \rangle}}$ for $k=1,2$. Hence $\Lambda\subseteq C$.
\end{proof}

To see that the spin factor $H\oplus\mathbb{R}$ is an atomic JBW-algebra, by \Cref{L: projections in spin}, it remains to show that there is an atom below $(0,1)$. Actually, for every atom $(x,\frac{1}{2})$, we have that $(x,\frac{1}{2}) \le (0,1)$ by \Cref{L: cone of squares in spin id}. Thus, $H\oplus \mathbb{R}$ is an atomic JBW-algebra. 

Next, we consider the spectrum of elements of $H\oplus\mathbb{R}$. Suppose that $(x,\lambda)$ is not a multiple of $(0,1)$. Then we can write 
\[
(x,\lambda) = \lambda(0,1) + \sqrt{\langle x,x \rangle}\left(\frac{x}{\sqrt{\langle x,x \rangle}},0\right)
\]
where $\langle x,x \rangle^{-1/2}(x,0)$ squares to $(0,1)$. By the functional calculus \cite[Corollary~1.19]{Alfsen}, the spectrum of the element $\langle x,x \rangle^{-1/2}(x,0)$ must be $\{\pm 1\}$, so the spectrum of the element $(x,\lambda)$ must, therefore, equal $\{\lambda + \sqrt{\langle x,x \rangle},\lambda - \sqrt{\langle x,x \rangle}\}$. On the other hand, any multiple of $(0,1)$ has a spectrum containing at most two numbers, so any element of $H\oplus\mathbb{R}$ has a spectrum consisting of at most two numbers. Furthermore, 
\[
(x,\lambda) = (\lambda + \sqrt{\langle x,x \rangle})\left(\frac{x}{2\sqrt{\langle x,x \rangle}},\frac{1}{2}\right) + (\lambda - \sqrt{\langle x,x \rangle})\left(-\frac{x}{2\sqrt{\langle x,x \rangle}},\frac{1}{2}\right)
\]
is the \emph{spectral decomposition} of $(x,\lambda)$. Note that the spectrum of $(x,\lambda)$ is positive if and only if $(x,\lambda)$ is an element of $C$. 

The cone $\Lambda:=\{(x,\lambda)\colon \sqrt{\langle x,x \rangle} \le \lambda\}$ in the vector space $H \oplus \mathbb{R}$ is called the \emph{Lorentz cone}. Clearly, $(H \oplus \mathbb{R}, \Lambda)$ is an Archimedean partially ordered vector space and $(0,1)$ is an order unit. In particular, for every $n\ge 3$, the vector space $\mathbb{R}^n$ can be endowed with a Lorentz cone, by viewing $\mathbb{R}^n$ as $\mathbb{R}\times \mathbb{R}^{n-1}$ and considering the Euclidean inner product on $\mathbb{R}^{n-1}$. 

\begin{lemma}\label{L: spin is a factor}
	A spin factor $H\oplus\mathbb{R}$ is in fact a factor.
\end{lemma}

\begin{proof}
	Let $x,y\in H$ be orthogonal unit vectors. Then $(x,0)\circ(y,0)^2 = (x,0)\circ(0,1) = (x,0)$ and $(y,0)\circ \bigl((x,0)\circ(y,0)\bigr) = (0,0)$, so $(x,0)$ does not operator commute with $(y,0)$, thus $(x,0)$ can not be in the algebraic centre of $H\oplus\mathbb{R}$. This implies that if $(x,\lambda)=(x,0)+\lambda (0,1)$ is an element of the algebraic centre, then $(x,0)$ must be an element of the algebraic centre as $\lambda(0,1)$ is. We conclude that $x=0$ and, therefore, the algebraic centre equals $\mathbb{R}(0,1)$ and $H\oplus\mathbb{R}$ is a factor. 
\end{proof}

We show that a functional $\varphi$ is a state of $H\oplus\mathbb{R}$ if and only if there exists $y\in H$ with $\langle y,y\rangle =1$ such that $\varphi((x,\lambda))=\langle (x,\lambda),(y,1)\rangle$ for every $(x,\lambda)\in H\oplus\mathbb{R}$. Indeed, let $\varphi$ be a state of $H\oplus\mathbb{R}$. By the Riesz representation theorem, it follows that there is a $(y,\mu)$ such that $\varphi((x,\lambda))=\langle (x,\lambda),(y,\mu) \rangle = \langle  x,y\rangle + \lambda\mu$ for all $(x,\lambda) \in H\oplus\mathbb{R}$. Since $\varphi(0,1) = 1$, we must have $\mu=1$. If $y\neq 0$, then  
\[
\varphi\left(\left(\frac{-y}{2\sqrt{\langle y,y \rangle}},\frac{1}{2}\right)\right) = \left\langle \left(\frac{-y}{2\sqrt{\langle y,y \rangle}},\frac{1}{2}\right),(y,1)\right\rangle = \frac{1}{2}\left(-\sqrt{\langle y,y \rangle} + 1\right) \ge 0
\]
since $\varphi$ is positive, so $\sqrt{\langle y,y \rangle} \le 1$. On the other hand, if $y\in H$ is such that $\sqrt{\langle y,y \rangle}\le 1$ and we define the linear functional $\psi$ by $\psi((x,\lambda)):=\langle x,y \rangle + \lambda$, then $\psi((0,1))=1$ and 
\begin{align*}
	\psi((x,\lambda)^2) &= 2\lambda \langle x,y \rangle + \langle x,x \rangle + \lambda^2 \ge
	-2|\lambda| |\langle x,y \rangle| + \langle x,x \rangle + \lambda^2
	\ge -2|\lambda|\sqrt{\langle x,x \rangle} + \langle x,x \rangle + \lambda^2 \\&= \left(\sqrt{\langle x,x \rangle} - |\lambda|\right)^2 \ge 0
\end{align*}
by the Cauchy-Schwarz inequality, so $\psi$ is positive. Hence, $\psi$ is a state. 

\begin{lemma}\label{L:pure states in spin}
	Let $(y,1)$ represent the state $\varphi$ on $H\oplus\mathbb{R}$. Then $\varphi$ is a pure state if and only if $\sqrt{\langle y,y \rangle}=1$.
\end{lemma}

\begin{proof}
	Suppose that $y \in H$ is such that $\sqrt{\langle y,y \rangle}<1$. If $y=0$, then we can write $(0,1) = \frac{1}{2}(x,1)+\frac{1}{2}(-x,1)$ for a unit vector  $x\in H$, and the states represented by $(x,1)$ and $(-x,1)$ are distinct since 
	\[
	\langle (x,1),(-x,1) \rangle = -\langle x,x \rangle + 1 = 0\qquad\mbox{and}\qquad\langle (-x,1),(-x,1) \rangle = \langle x,x \rangle + 1 = 2.
	\]
	Hence, the state represented by $(0,1)$ is not a pure state. If $y\neq 0$, then we can write 
	\[
	(y,1) = t\left(\frac{y}{\sqrt{\langle y,y \rangle}},1\right) + (1-t)\left(\frac{-y}{\sqrt{\langle y,y \rangle}},1\right)
	\]
	for some $0<t<1$, and the states represented by $(\pm\ \langle y,y \rangle^{-1/2}y,1)$ are again distinct since 
	\[
	\left\langle \left(\frac{y}{\sqrt{\langle y,y \rangle}},1\right),\left(\frac{-y}{\sqrt{\langle y,y \rangle}},1\right) \right\rangle = 0\qquad\mbox{and}\qquad\left\langle \left(\frac{-y}{\sqrt{\langle y,y \rangle}},1\right),\left(\frac{-y}{\sqrt{\langle y,y \rangle}},1\right) \right\rangle = 2.
	\]
	We conclude that the state represented by $(y,1)$ can, therefore, not be a pure state. Conversely, suppose that $y$ is a unit vector in $H$. If $(x,1)$ and $(z,1)$ represent states such that a non-trivial convex combination of them equal the state represented by $(y,1)$, then $(y,1)=t(x,1)+(1-t)(z,1)$ for some $0<t<1$ and so $y=tx+(1-t)z$. Since the unit sphere in $H$ is strictly convex, it follows that $x=z=y$, so $(y,1)$ represents a pure state.  
\end{proof}

\subsection{Matrices with octonionic entries}
We introduce the multiplication rules on $\{e_1,\ldots, e_7\}$ as follows. Set $e_i^2=-1$ for all $1\le i\le 7$ and determine the product of any two $e_i$ and $e_j$ via the so called \emph{Fano plane} below.
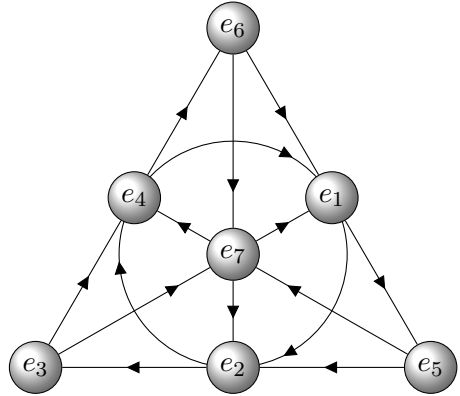
\[
\begin{tikzpicture}[scale=1.5]
\tikzstyle{point}=[ball color=white, circle, draw=black, inner sep=0.1cm]
\node (v7) at (0,0) [point] {$e_7$};
\draw (0.866,0.5) arc (30:160:1cm) node[currarrow, pos=0.25, xscale=1, sloped, scale=1]{};
\draw (-0.9397,0.342) arc (160:270:1cm) node[currarrow, pos=0.25, xscale=-1, sloped, 
scale=1]{}; \draw (0,-1) arc (270:390:1cm) node[currarrow, pos=0.25, xscale=-1, sloped, scale=1]{};
\node (v1) at (90:2cm) [point] {$e_6$};
\node (v2) at (210:2cm) [point] {$e_3$};
\node (v4) at (330:2cm) [point] {$e_5$};
\node (v3) at (150:1cm) [point] {$e_4$};
\node (v6) at (270:1cm) [point] {$e_2$};
\node (v5) at (30:1cm) [point] {$e_1$};
\draw (v1) -- (v3) node[currarrow, pos=0.5, xscale=1, sloped, scale=1]{} -- (v2) node[currarrow, pos=0.5, xscale=1, sloped, scale=1]{};
\draw (v2) -- (v6) node[currarrow, pos=0.5, xscale=-1, sloped, scale=1]{} -- (v4) node[currarrow, pos=0.5, xscale=-1, sloped, scale=1]{};
\draw (v4) -- (v5) node[currarrow, pos=0.5, xscale=1, sloped, scale=1]{} -- (v1) node[currarrow, pos=0.5, xscale=1, sloped, scale=1]{};
\draw (v3) -- (v7) node[currarrow, pos=0.5, xscale=-1, sloped, scale=1]{} -- (v4) node[currarrow, pos=0.25, xscale=-1, sloped, scale=1]{};
\draw (v5) -- (v7) node[currarrow, pos=0.5, xscale=1, sloped, scale=1]{} -- (v2) node[currarrow, pos=0.25, xscale=1, sloped, scale=1]{};
\draw (v6) -- (v7) node[currarrow, pos=0.5, xscale=-1, sloped, scale=1]{} -- (v1) node[currarrow, pos=0.25, xscale=-1, sloped, scale=1]{};
\end{tikzpicture}
\]
The elements $e_i$ and $e_j$ lie on a unique line consisting of three elements, including the circle. The product is defined by following the arrow and using cyclic permutations (which preserve the directions of the arrows). For example, the elements $e_1$ and $e_2$ lie on the line $(e_1,e_2,e_4)$, so $e_1 e_2=e_4$. Furthermore, we also have that $e_1$ and $e_4$ lie on the line $(e_1,e_2,e_4)$, which yields the same line $(e_4,e_1,e_2)$ by applying cyclic permutations, so $e_4e_1=e_2$. Transversing in the opposite direction of the indicated arrow yields a minus sign, that is, $e_1e_4=-e_2$. The 8-dimensional real vector space
\[
\mathbb{O}:=\bigl\{a_01 + a_1 e_1 +\dots+ a_7 e_7 \colon a_0,\dots,a_7 \in \mathbb{R}\bigr\}
\]
equipped with the multiplication rules described above and unit $1$, where a general product distributes over the sums, forms the so called \emph{octonions}. The product on the octonions is not commutative as we have seen, and it also fails to be associative. Indeed, note that $(e_1e_2)e_3=e_4e_3=-e_6$ and $e_1(e_2e_3)=e_1e_5=e_6$. The octonions are \emph{alternative}, meaning that the subalgebra generated by two elements in $\mathbb{O}$ is associative. The real multiples of the identity $1$ commute with all octonions. For $x=a_01+\sum_{k=1}^7 a_ke_k$, the \emph{octonionic conjugate} of $x$ is defined to be $x^*:= a_01-\sum_{k=1}^7a_ke_k$. The octonionic conjugate is an involution on $\mathbb{O}$ that reverses the order of multiplication, that is, for $x,y\in\mathbb{O}$, we have $(xy)^*=y^*x^*$. The \emph{real part} of $x$ is denoted by $\operatorname{Re}(x)$ and is given by $\operatorname{Re}(x):=\frac{1}{2}(x+x^*)=a_01$. Note that every non-zero octonion is invertible since $x^*x=(a_0^2+\dots+a_7^2)1$. Furthermore, the octonionic conjugation induces a norm on $\mathbb{O}$ given by $\|x\|:=\sqrt{x^*x}$, where $\mathbb{R}1$ has been identified with $\mathbb{R}$. Similar to the norm on the quaternions, the norm is multiplicative, that is, $\|xy\|=\|x\|\|y\|$ for all $x,y\in\mathbb{O}$. This implies that
\[
\|xy\|^2=(xy)^*(xy)=(y^*x^*)(xy)=y^*(x^*x)y=\|x\|^2\|y\|^2.
\]
The octonions can be equipped with the real inner product $\langle x,y \rangle := \operatorname{Re}(xy^*)= \frac{1}{2}(xy^*+yx^*)$, where again the real multiples of the identity $1$ are identified with the real numbers. The inner product coincides with the standard inner product on $\mathbb{R}^8$, that is, for $x:=a_01+\sum_{k=1}^7a_ke_k$ and $y:=b_01+\sum_{k=1}^7b_ke_k$ it follows that $\langle x,y \rangle = \sum_{k=0}^7a_kb_k$. Furthermore, note that the norm relates to the inner product as usual, $\|x\|=\sqrt{\langle x,x \rangle}$. For the reader interested in studying properties of the octonions in more detail, we recommend the well written and extensive exposition on the subject \cite{Baez}. 

In view of the theory of JB-algebras, let $\mathrm{M}_n(\mathbb{O})$ denote the $n\times n$ matrices over the octonions which form a non-associative unital real algebra. Similar to the Hermitian adjoint, an involution can be defined on $\mathrm{M}_n(\mathbb{O})$ given by $(A^*)_{ij}:=(A_{ji})^*$. Since every JB-algebra is formally real, the subspace of self-adjoint matrices $\mathrm{M}_n(\mathbb{O})_\mathrm{sa}$ are considered instead of $\mathrm{M}_n(\mathbb{O})$, equipped with the commutative product $A\circ B:=\frac{1}{2}(AB+BA)$ (note that squares coincide for both products). It was shown by Jordan, von Neumann, and Wigner in \cite{Jordan-von-Neumann-Wigner} that $\mathrm{M}_n(\mathbb{O})_\mathrm{sa}$ is a Jordan algebra for $1\le n \le 3$ and not for $n\ge 4$, see also \cite[Theorem~2.7.6, Theorem~2.7.8]{Hanche-Olsen-Stormer}. In particular, it turns out that $\mathrm{M}_2(\mathbb{O})_\mathrm{sa}$ is a spin factor, see \cite[p.~28]{Baez}. 

\begin{lemma}\label{L:M_2(O) is spin}
	$\mathrm{M}_2(\mathbb{O})_\mathrm{sa}$ is a spin factor.
\end{lemma}

\begin{proof}
	Define the map $f \colon \mathrm{M}_2(\mathbb{O})_\mathrm{sa} \to \mathbb{R}^9 \oplus \mathbb{R}$ by 
	\[
	\begin{pmatrix}
		\alpha+\beta & x \\
		x^* & \alpha-\beta
	\end{pmatrix} \mapsto ((x,\beta),\alpha).
	\]
	It is a straightforward verification that $f$ is a linear bijection that maps the identity matrix $I_2$ to the unit $(0,1)$. Let 
	\[
	A:=\begin{pmatrix}
		\alpha_1+\beta_1 & x \\
		x^* & \alpha_1-\beta_1
	\end{pmatrix} \qquad\mbox{and}\qquad B:=\begin{pmatrix}
		\alpha_2+\beta_2 & y \\
		y^* & \alpha_2-\beta_2
	\end{pmatrix}.
	\]
	Using that $\langle x,y \rangle = \langle x^*,y^* \rangle$, it follows that 
	\[
	A\circ B = \begin{pmatrix}
		a+b & \alpha_2 x+\alpha_1 y \\
		\alpha_2 x^*+\alpha_1 y^* & a-b
	\end{pmatrix}
	\]
	where $a=\alpha_1\alpha_2+\beta_1\beta_2+\langle x,y \rangle$ and $b=\alpha_1\beta_2+\alpha_2\beta_1$, so 
	\begin{align*}
		f(A\circ B)&=((\alpha_2 x+\alpha_1 y,\alpha_1\beta_2+\alpha_2\beta_1),\alpha_1\alpha_2+\beta_1\beta_2+\langle x,y \rangle)=((x,\beta_1),\alpha_1)\circ((y,\beta_2),\alpha_2)\\&=f(A)\circ f(B)
	\end{align*}
	showing that $f$ is a Jordan homomorphism. Hence $\mathrm{M}_2(\mathbb{O})_\mathrm{sa}$ is isomorphic to the spin factor $\mathbb{R}^9\oplus\mathbb{R}$. 
\end{proof}

\begin{remark}\label{R:spin}
	A similar argument proves that $\mathrm{M}_2(\mathbb{H})_\mathrm{sa}$ is isomorphic to the spin factor $\mathbb{R}^5\oplus\mathbb{R}$, that $\mathrm{M}_2(\mathbb{C})_\mathrm{sa}$ is isomorphic to the spin factor $\mathbb{R}^3\oplus\mathbb{R}$, and that $\mathrm{M}_2(\mathbb{R})_\mathrm{sa}$ is isomorphic to the spin factor $\mathbb{R}^2\oplus\mathbb{R}$. Therefore, by \Cref{L: projections in spin}, all the minimal projections in $\mathrm{M}_2(\mathbb{R})_\mathrm{sa}$ are of the form 
	\[
	\begin{pmatrix}
		\frac{1}{2}+x_2&x_1\\x_1&\frac{1}{2}-x_2
	\end{pmatrix}
	\]
	where $x_1^2+x_2^2=\frac{1}{4}$.
\end{remark}

The self-adjoint $3\times 3$ matrices over the octonions is called the \emph{Albert algebra} and is an exceptional Jordan algebra, as it is not Jordan isomorphic to a subalgebra of an associative real algebra $A$ with the product $a\circ b:=\frac{1}{2}(ab+ba)$, see \cite[Theorem~4.6]{Alfsen} and \cite[Corollary~2.8.5]{Hanche-Olsen-Stormer}. Furthermore, by \cite[Theorem~3.32]{Alfsen}, the Albert algebra is a JBW-algebra and even a factor. Hence, it follows from \cite[Lemma~1.10]{Alfsen} that $\mathrm{M}_3(\mathbb{O})_\mathrm{sa}$ is partially ordered by the cone of squares. The minimal projections (or atoms) in $\mathrm{M}_3(\mathbb{O})_\mathrm{sa}$ can be characterised as follows.

\begin{proposition}\label{L:projections in Albert algebra}
	The minimal projections $P$ in $\mathrm{M}_3(\mathbb{O})_\mathrm{sa}$ are of the form 
	\[P:=
	\begin{pmatrix}
		\|x_1\|^2 & x_1x_2^* & x_1x_3^*\\
		x_2x_1^* & \|x_2\|^2 & x_2x_3^*\\
		x_3x_1^* & x_3x_2^* & \|x_3\|^2
	\end{pmatrix}
	\]
	where $x_1,x_2,x_3\in\mathbb{O}$ associate, that is, $(x_1x_2)x_3=x_1(x_2x_3)$, and $\|x_1\|^2+\|x_2\|^2+\|x_3\|^2=1$.
\end{proposition}

\begin{proof}
	Let 
	\[A:=
	\begin{pmatrix}
		r_1 & y_1^* & y_2^*\\
		y_1 & r_2 & y_3^*\\
		y_2 & y_3 & r_3
	\end{pmatrix}
	\]
	be so that $A^2=A$. Then, as 
	\begin{align}\label{E:A^2}
		A^2=
		\begin{pmatrix}
			r_1^2+\|y_1\|^2+\|y_2\|^2 & (r_1+r_2)y_1^*+y_2^*y_3 & (r_1+r_3)y^*_2+y^*_1y_3^*\\
			(r_1+r_2)y_1+y_3^*y_2 & r_2^2+\|y_1\|^2+\|y_3\|^2 &  (r_2+r_3)y^*_3+y_1y_2^*\\
			(r_1+r_3)y_2+y_3y_1 & (r_2+r_3)y_3+y_2y_1^* & r_3^2+\|y_2\|^2+\|y_3\|^2
		\end{pmatrix},
	\end{align}
	it follows that $0\le r_1,r_2,r_3\le 1$ and not all the $r_i$ are zero, as otherwise $A=0$. Furthermore, from the system of equations
	\[
	\begin{cases}
		(1-r_1-r_2)y_1=y_3^*y_2\\
		(1-r_1-r_3)y_2=y_3y_1\\
		(1-r_2-r_3)y_3=y_2y_1^*
	\end{cases}
	\]
	we see that $y_1,y_2$, and $y_3$ are in a subalgebra $N\subseteq\mathbb{O}$ generated by two elements (and 1). Since $\mathbb{O}$ is alternative, we must have that $N$ is associative. Let $x\in N$ be non-zero. Since $\mathbb{O}$ has no zero divisors, the $\mathbb{R}$-linear map $L_x(y):=xy$ is injective on $\mathbb{O}$. Hence, the restriction of $L_x$ to $N$ is injective as well and as $N$ is finite-dimensional, it is a bijection. Let $z\in N$ be such that $xz=1$. It follows that $(zx)^2=(zx)(zx)=z(xz)x=zx$ as $N$ is associative, and so $zx=1$, again since $\mathbb{O}$ has no zero divisors. This shows that $x$ has an inverse $z$ and, therefore, $N$ is a real division algebra. By Hurwitz's theorem \cite{Hurwitz}, we have that $N$ is isomorphic to $\mathbb{R}$, $\mathbb{C}$, or $\mathbb{H}$. It follows that the entries of $A$ are elements of the algebra $N$ which is isomorphic to $\mathbb{H}$. Under this isomorphism, the inner product on $N^3$ induced by the inner product of $\mathbb{O}^3$ coincides with the inner product of $\mathbb{H}^3$.
	%every projection in $M_3(\mathbb{O})_{sa}$ can be viewed as a projection in $M_3(\mathbb{H})_{sa}$. 
	Hence, by \Cref{L:P_v is minimal}, there is a unital vector $x\in  N^3$ such that $Ay=x \cdot \langle x,y \rangle $. It follows that  $x:=(x_1,x_2,x_3)$ is a unital vector in $\mathbb{O}^3$ with $(x_1x_2)x_3=x_1(x_2x_3)$ and 
	\begin{equation}\label{eq.theformofA}
		A=\begin{pmatrix}
			\|x_1\|^2 & x_1x_2^* & x_1x_3^*\\
			x_2x_1^* & \|x_2\|^2 & x_2x_3^*\\
			x_3x_1^* & x_3x_2^* & \|x_3\|^2
		\end{pmatrix}.
	\end{equation}
	
	Conversely, suppose $A$ is as in \eqref{eq.theformofA} for some unit vector $(x_1,x_2,x_3)\in\mathbb{O}^3$ such that $(x_1x_2)x_3=x_1(x_2x_3)$. Then it follows that $A^2=A$, and as the subalgebra $M\subseteq\mathbb{O}$ generated by $x_1$, $x_2$, and $x_3$ (and 1) is associative, it follows that $M$ is isomorphic to $\mathbb{R}$, $\mathbb{C}$, or $\mathbb{H}$ by Hurwitz's theorem once more. Therefore, the matrix $A$ satisfies $Ay=x\cdot \langle x,y \rangle$ for $x=(x_1,x_2,x_3)$, so it is a minimal projection by \Cref{L:P_v is minimal}. 
\end{proof}

The \emph{trace} of $A\in\mathrm{M}_3(\mathbb{O})_\mathrm{sa}$ is defined as usual for matrices by $\mathrm{trace}(A):=A_{11}+A_{22}+A_{33}$, where $A_{ii}$ are the diagonal entries of $A$. It follows that $\langle A,B \rangle := \mathrm{trace}(A\circ B)$ is a real inner product on $\mathrm{M}_3(\mathbb{O})_\mathrm{sa}$. Indeed, note that, by \eqref{E:A^2}, it follows that $\langle A,A \rangle \ge 0$, and $\langle A,A \rangle = 0$ if and only if $A=0$. Furthermore, $\mathrm{trace}(A\circ B)=\operatorname{Re}(\mathrm{trace}(AB))$ and, by \cite[Proposition~V.2.2]{Faraut-Koranyi}, we have that $\langle A\circ B, C \rangle = \langle A, B\circ  C \rangle$ for all $A,B,C\in \mathrm{M}_3(\mathbb{O})_\mathrm{sa}$. Hence, with this inner product $\mathrm{M}_3(\mathbb{O})_\mathrm{sa}$ is a Euclidean Jordan algebra. For any $A\in \mathrm{M}_3(\mathbb{O})_\mathrm{sa}$, there are unique $\lambda_1,\dots,\lambda_m$ and unique pairwise orthogonal projections $P_1,\dots,P_m$ such that $A=\lambda_1 P_1+\dots+\lambda_mP_m$ by \cite[Theorem~III.1.1]{Faraut-Koranyi}. This is the \emph{spectral decomposition} of $A$. The \emph{spectrum} of $A$, denoted by $\sigma(A)$, consists of the eigenvalues that occur in the spectral decomposition of $A$, that is, $\sigma(A)=\{\lambda_1,\ldots,\lambda_m\}$.

\begin{lemma}\label{L: positive equivalence M_3(O)_sa}
	Let $A\in \mathrm{M}_3(\mathbb{O})_\mathrm{sa}$. Then the following statements are equivalent.
	\begin{itemize}
		\item [$(i)$] $A\ge 0$.
		\item[$(ii)$] $\sigma(A)\subseteq [0,\infty)$.
		\item[$(iii)$] $\langle A,B \rangle\ge 0$ for all $B\ge 0$.
	\end{itemize}
\end{lemma}

\begin{proof}
	$(i)\Longleftrightarrow (ii)$: If $A\ge 0$, then $A=B^2$ for some $B\in\mathrm{M}_3(\mathbb{O})_\mathrm{sa}$, and the spectral decomposition of $B=\lambda_1P_1+\dots+\lambda_mP_m$ now yields $A=B^2=\lambda_1^2P_1+\dots+\lambda_m^2P_m$, so $\sigma(A)\subseteq [0,\infty)$. On the other hand, if $\sigma(A)\subseteq [0,\infty)$, then the spectral decomposition of $A=\mu_1Q_1+\dots+\mu_nQ_n$ yields $A=B^2$ for $B:=\sqrt{\mu_1}Q_1+\dots+\sqrt{\mu_n}Q_n$.
	
	$(i)\Longleftrightarrow (iii)$: This equivalence follows from the fact that the cone of squares in a Euclidean Jordan algebra yields a symmetric cone by \cite[Theorem~III.2.1]{Faraut-Koranyi}. In particular, $\langle A,B\rangle \ge 0$ for all $B\ge 0$ if and only if $A\ge 0$.
\end{proof}

By the Riesz representation theorem, for every functional $\varphi\colon \mathrm{M}_3(\mathbb{O})_\mathrm{sa} \to \mathbb{R}$, there is a unique $B\in \mathrm{M}_3(\mathbb{O})_\mathrm{sa}$ such that $\varphi(A)=\langle A, B\rangle$. Furthermore, it follows from \Cref{L: positive equivalence M_3(O)_sa} that $\varphi:=\langle \cdot,B\rangle$ is a state if and only if $B\ge 0$ and $\mathrm{trace}(B)=1$.

\begin{lemma}\label{L:pure_states_in_Albert}
	A state $\langle \cdot,B\rangle$ on $\mathrm{M}_3(\mathbb{O})_\mathrm{sa}$ is pure if and only if $B$ is a minimal projection, i.e., an atom. 
\end{lemma}

\begin{proof} 
	Suppose that $\langle \cdot,B\rangle$ is a pure state on $\mathrm{M}_3(\mathbb{O})_\mathrm{sa}$. Let $B=\sum_{k=1}^m \lambda_k P_k$ be the spectral decomposition of $B$ such that $\lambda_k\neq 0$ for all $k$. Since every $P_k$ can be written as the sum of minimal projections, we may assume that each $P_k$ is minimal. Suppose that there are two distinct minimal projections $P_i$ and $P_j$ in this decomposition. Then we can write 
	\begin{align}\label{E:pure_state}
		\langle \cdot,B\rangle=\lambda_i\langle \cdot,P_i\rangle+\textstyle{\sum_{k\ne j}\lambda_k\left((\sum_{k\ne j}\lambda_k)^{-1}\left\langle \cdot,\sum_{k\ne j}\lambda_kP_k\right\rangle\right)}
	\end{align}
	and since $\mathrm{trace}(P_k)=1$ for all $k$ by \Cref{L:projections in Albert algebra}, it follows that \eqref{E:pure_state} writes $\langle \cdot,B\rangle$ as a non-trivial convex combination of two distinct states, which is impossible. Hence, we have that $\langle \cdot,B\rangle = \lambda_k\langle \cdot,P_k\rangle$ for some $k$ and as $\mathrm{trace}(B)=1$, we find that $B=P_k$.   
	
	On the other hand, let $P$ be a minimal projection and suppose $\langle \cdot,P\rangle=t\langle\cdot ,C\rangle+(1-t)\langle \cdot, D\rangle$ for $0<t<1$, $C,D\ge 0$ with $\mathrm{trace}(C)=\mathrm{trace}(D)=1$. Then we must have $P=tC+(1-t)D$ and so $C=\lambda P$ and $D=\mu P$ by \cite[Proposition~2.15]{Alfsen} and \cite[Lemma~3.29]{Alfsen} since $P$ is a minimal projection. Because $\mathrm{trace}(C)=\mathrm{trace}(D)=1$, it follows that $\lambda=\mu =1$ and so $\langle \cdot,P\rangle$ is a pure state. 
\end{proof}

\subsection{The pre-duals of atomic JBW-algebra factors}

We conclude this appendix by determining the pre-duals of all atomic JBW-algebra factors $M$. For this we need the notion of so called \emph{trace class} elements. These are elements $x \in M$ that can be written as 
\[
x = \sum_{k=1}^\infty \lambda_k p_k,
\]
where $(p_k)_k$ is a sequence of pairwise orthogonal atoms in $M$ and $(\lambda_k)_k \subseteq \R$ satisfies $\sum_{k=1}^\infty |\lambda_k| < \infty$. The set of trace class elements will be denoted by $M_{\mathrm{tr}}$, and a \emph{trace} can be define on $M_{\mathrm{tr}}$ by 
\[
\mathrm{tr}(x) := \sum_{k=1}^\infty \lambda_k.
\]
The trace does not depend on the representation of $x$, so it is well defined on $M_\mathrm{tr}$; see \cite[Definition~5.65]{Alfsen} and the paragraph below for more details. Given $x \in M$, we consider the JB-subalgebra JB$(x,e)$ of $M$ generated by $x$ and $e$, which is isomorphic to a space of continuous functions. Hence, in JB$(x,e)$ the modulus $|x|$ of $x$ exists. The \emph{trace norm} of $x \in M_{\mathrm{tr}}$ is defined by 
\[ 
\|x\|_\mathrm{tr}:=\mathrm{tr}(|x|) = \sum_{k=1}^\infty |\lambda_k|.
\] 
It follows from \cite[Proposition~5.66]{Alfsen} that $M_{\mathrm{tr}}$ equipped with the trace norm is a Banach space. The pre-dual $M_*$ of $M$ is isometrically isomorphic to $M_{\mathrm{tr}}$. In particular, for all finite-dimensional factors and spin factors, the pre-dual is the same space equipped with the trace norm. For the self-adjoint bounded operators on a real, complex, or quaternionic Hilbert space, the pre-dual is identified with the three analogues of trace class operators. 

\begin{acknowledgement}
    The authors are grateful to A.W.\ Wickstead for making them aware of reference \cite{Buck}.
\end{acknowledgement}

\bibliographystyle{alpha}
\bibliography{bandsJBW_nov2022}

\end{document}